\documentclass[11pt]{amsart}
\usepackage{amsmath}
\usepackage{amsxtra}
\usepackage{amscd}
\usepackage{amsthm}
\usepackage{amsfonts}
\usepackage{amssymb}
\usepackage{eucal}

\newtheorem{thm}{Theorem}[section]
\newtheorem{conj}{Conjecture}

\newtheorem{lem}[thm]{Lemma}
\newtheorem{prop}[thm]{Proposition}

\theoremstyle{remark}
\newtheorem{remark}[thm]{Remark}

\theoremstyle{definition}
\newtheorem{definition}[thm]{Definition}

\numberwithin{equation}{section}

\newcommand{\bean}{\begin{eqnarray}}
\newcommand{\eean}{\end{eqnarray}}
\newcommand{\be}{\begin{displaymath}}
\newcommand{\ee}{\end{displaymath}}
\newcommand{\bea}{\begin{eqnarray*}}
\newcommand{\eea}{\end{eqnarray*}}

\newcommand{\thmref}[1]{Theorem~\ref{#1}}
\newcommand{\secref}[1]{Section~\ref{#1}}
\newcommand{\lemref}[1]{Lemma~\ref{#1}}
\newcommand{\propref}[1]{Proposition~\ref{#1}}

\newcommand{\conjref}[1]{Conjecture~\ref{#1}}
\newcommand{\nc}{\newcommand}
\nc{\on}{\operatorname}
\nc{\ch}{\mbox{ch}}
\nc{\Z}{{\mathbb Z}}
\nc{\C}{{\mathbb C}}
\nc{\pone}{{\mathbb P}^1}
\nc{\pa}{\partial}
\nc{\F}{{\mathcal F}}
\nc{\arr}{\rightarrow}
\nc{\larr}{\longrightarrow}
\nc{\al}{\alpha}
\nc{\ri}{\rangle}
\nc{\lef}{\langle}
\nc{\W}{{\mathcal W}}
\nc{\la}{\lambda}
\nc{\ep}{\epsilon}
\nc{\su}{\widehat{{\mathfrak s}{\mathfrak l}}_2}
\nc{\sw}{{\mathfrak s}{\mathfrak l}}
\nc{\g}{{\mathfrak g}}
\nc{\h}{{\mathfrak h}}
\nc{\n}{{\mathfrak n}}
\nc{\N}{\widehat{\n}}
\nc{\G}{\widehat{\g}}
\nc{\De}{\Delta}
\nc{\gt}{\widetilde{\g}}
\nc{\Ga}{\Gamma}
\nc{\one}{{\mathbf 1}}
\nc{\z}{{\mathfrak Z}}
\nc{\La}{\Lambda}
\nc{\wt}{\widetilde}
\nc{\wh}{\widehat}
\nc{\cri}{_{\kappa_c}}
\nc{\kk}{h^\vee}
\nc{\sun}{\widehat{\sw}_N}
\nc{\si}{\sigma}
\nc{\el}{\ell}
\nc{\bi}{\bibitem}
\nc{\om}{\omega}
\nc{\ol}{\overline}
\nc{\ds}{\displaystyle}
\nc{\dzz}{\frac{dz}{z}}
\nc{\Res}{\on{Res}}
\nc{\mc}{\mathcal}
\nc{\Cal}{\mathcal}
\nc{\bb}{{\mathfrak b}}
\nc{\ot}{\otimes}
\nc{\R}{{\mc R}}
\nc{\yy}{{\mc Y}}
\nc{\ga}{\gamma}

\nc{\us}{\underset}
\nc{\opl}{\oplus}
\nc{\beq}{\begin{equation}}
\nc{\Fq}{{\mathcal F}}
\nc{\Mq}{{\mathcal M}}
\nc{\Rep}{\on{Rep}}
\nc{\sssec}{\subsubsection}
\nc{\ssec}{\subsection}
\nc{\lan}{\langle}
\nc{\ran}{\rangle}

\nc{\D}{\mathcal D}
\nc{\Vect}{\on{Vect}}
\nc{\ghat}{\G}
\nc{\T}{\mc T}
\nc{\Tloc}{\T^\g_{\on{loc}}}
\nc{\vac}{|0\ran}
\nc{\Wick}{{\mb :}}
\nc{\mb}{\mathbf}
\nc{\delz}{\partial_z}
\nc{\K}{{\cali K}}
\nc{\cali}{\mathcal}
\nc{\li}{\mathfrak l}
\nc{\lt}{\widetilde{\li}}
\nc{\astar}{a^*}
\nc{\cA}{{\mc A}}
\nc{\ka}{\kappa}

\nc{\OO}{{\mc O}}
\nc{\AutO}{\on{Aut}\OO}
\nc{\DerO}{\on{Der}\OO}
\nc{\DerpO}{\on{Der}_+\OO}
\nc{\Au}{{\mc A}ut}
\nc{\mf}{\mathfrak}
\nc{\V}{{\mathbb V}}
\nc{\hh}{\wh{\h}}

\nc{\pp}{{\mathfrak p}}
\nc{\mm}{{\mathfrak m}}
\nc{\rr}{{\mathfrak r}}
\nc{\ket}{\rangle}
\nc{\zz}{{\mathfrak z}}
\nc{\gr}{\on{gr}}
\nc{\Spe}{\on{Spec}}
\nc{\rv}{\crho}
\nc{\can}{\on{can}}
\nc{\CC}{\on{Op}_G(D))}
\nc{\Op}{\on{Op}_G(D)}
\nc{\MOp}{\on{MOp}_G(D)}
\nc{\Db}{{\mathbb D}}
\nc{\ww}{w}

\nc{\af}{{\mathbb A}^1}
\nc{\bs}{\backslash}
\nc{\laa}{(\la_i)}
\nc{\zn}{(z_i)}

\nc{\cla}{\check{\la}}
\nc{\cmu}{\check{\mu}}
\nc{\crho}{\check{\rho}}
\nc{\chal}{\check{\al}}
\nc{\cc}{{\mathfrak c}}

\nc{\M}{{\mathbb M}}

\nc{\ZZ}{{\mc Z}}

\nc{\UU}{{\mathbb U}}

\nc{\Conn}{\on{Conn}(\Omega^{\crho})}
\nc{\Con}{\on{Conn}(\Omega^{-\rho})}
\nc{\Co}{\on{Conn}(\Omega^{\rho})}

\nc{\ppart}{(\!(t)\!)}
\nc{\zpart}{(\!(z)\!)}
\nc{\ppzi}{(\!(t-z_i)\!)}
\nc{\ppinf}{(\!(t^{-1})\!)}
\nc{\Ind}{\on{Ind}}
\nc{\I}{{\mathbb I}}

\nc{\Bun}{\on{Bun}}

\newcommand {\bfone}{\mathbf{1}}

\newcommand {\IC}{\mathbb{C}}

\newcommand {\IM}{\mathbb{M}}
\newcommand {\IN}{\mathbb{N}}                          
\newcommand {\IP}{\mathbb{P}}

\newcommand {\IZ}{\mathbb{Z}}

\newcommand {\VV}{\mathbf{V}}

\renewcommand {\l}{\mathfrak l}
\newcommand {\m}{\mathfrak m}

\newcommand {\A}{{\mathcal A}}

\renewcommand {\sl}[1]{\mathfrak{sl}_{#1}}

\newcommand {\ad}{\operatorname{ad}}

\newcommand {\End}{\operatorname{End}}

\newcommand {\Ker}{\operatorname{Ker}}

\newcommand {\id}{\operatorname{id}}

\newcommand {\rk}{\operatorname{rk}}

\newcommand {\nablak}{\nabla_h}

\newcommand {\reg}{_{\operatorname{reg}}}
\newcommand {\hreg}{\h\reg}

\newcommand {\fd}{finite-dimensional }

\newcommand {\wrt}{with respect to }

\newcommand {\aand}{\qquad\text{and}\qquad}

\newcommand {\odots}[1]{#1\cdots #1}

\renewcommand {\>}{\rangle}

\newcommand {\hatg}{\wh{\g}}
\newcommand {\IV}{\mathbb{V}}
\newcommand {\Der}{\operatorname{Der}}
\newcommand {\tildeg}{\wt{\g}}
\newcommand {\ind}{\operatorname{ind}}
\newcommand {\opp}{^{\scriptscriptstyle{\operatorname{opp}}}}
\newcommand {\ord}{\operatorname{ord}}

\nc{\MM}{\mathcal M}

\begin{document}

\title{Gaudin models with irregular singularities}

\author{B. Feigin}

\address{Landau Institute for Theoretical Physics, Kosygina St 2,
Moscow 117940, Russia}

\author{E. Frenkel}\thanks{E.F. and V.T.L. were supported by the DARPA
grant HR0011-04-1-0031 and the NSF grant DMS-0303529, and B.F. was
supported by the grants RFBR 05-01-01007, RFBR 05-01-02934 and
NSh-6358.2006.2.}

\address{Department of Mathematics, University of California,
  Berkeley, CA 94720, USA}

\author{V. Toledano Laredo}

\address{Universit\'e Pierre et Marie Curie--Paris 6, UMR 7586,
Institut de Math\'ematiques de Jussieu, Case 191, 16 rue Clisson,
F--75013 Paris and Department of Mathematics, Northeastern University,
360 Huntington Avenue, Boston, MA 02115, USA}


\begin{abstract}
We introduce a class of quantum integrable systems generalizing the
Gaudin model. The corresponding algebras of quantum Hamiltonians are
obtained as quotients of the center of the enveloping algebra of an
affine Kac--Moody algebra at the critical level, extending the
construction of higher Gaudin Hamiltonians from \cite{FFR} to the case
of non-highest weight representations of affine algebras. We show that
these algebras are isomorphic to algebras of functions on the spaces
of opers on $\pone$ with regular as well as irregular singularities at
finitely many points. We construct eigenvectors of these Hamiltonians,
using Wakimoto modules of critical level, and show that their spectra
on finite-dimensional representations are given by opers with trivial
monodromy. We also comment on the connection between the generalized
Gaudin models and the geometric Langlands correspondence with
ramification.
\end{abstract}

\maketitle

\tableofcontents

\section{Introduction}

Quantum integrable systems associated to simple Lie algebras are
sometimes best understood in terms of the corresponding affine
Kac--Moody algebras. A case in point is the Gaudin model
\cite{G}. Let us recall the setup (see, e.g., \cite{FFR,F:faro}): for
each simple Lie algebra $\g$ there is a collection of commuting
quadratic Gaudin Hamiltonians $\Xi_i, i=1,\ldots,N$, in
$U(\g)^{\otimes N}$, defined for any set of $N$ distinct complex
numbers $z_1,\ldots,z_N$. A natural question is to find a maximal
commutative subalgebra of $U(\g)^{\otimes N}$ containing $\Xi_i,
i=1,\ldots,N$. It turns out that such a subalgebra may be
constructed with the help of the affine Kac--Moody algebra $\ghat$,
which is the universal central extension of the formal loop algebra
$\g\ppart$. The completed universal enveloping algebra
$\wt{U}_{\ka_c}(\ghat)$ of the affine Kac--Moody algebra at the {\em
critical level} contains a large center \cite{FF:gd,F:wak}. It was
shown in \cite{FFR} (see also \cite{F:faro}) that the sought-after
commutative algebra, called the {\em Gaudin algebra} in \cite{F:faro},
may be obtained as a quotient of the center of
$\wt{U}_{\ka_c}(\ghat)$, using the spaces of conformal blocks of
$\ghat$-modules. In particular, the quadratic Segal-Sugawara central
elements of $\wt{U}_{\ka_c}(\ghat)$ give rise to the quadratic Gaudin
Hamiltonians, whereas higher order central elements of
$\wt{U}_{\ka_c}(\ghat)$ give rise to higher order generalized Gaudin
Hamiltonians.

The center $Z(\ghat)$ of $\wt{U}_{\ka_c}(\ghat)$ has been described in
\cite{FF:gd,F:wak}, where it was shown that $Z(\ghat)$ is isomorphic
to the algebra of functions on the space of $^L G$-opers on the
punctured disc (for an introduction to this subject, see
\cite{F:book}). Here $^L G$ is the group of inner automorphisms of the
Langlands dual Lie algebra of $^L \g$ (whose Cartan matrix is the
transpose of the Cartan matrix of $\g$). The notion of opers has been
introduced by A. Beilinson and V. Drinfeld in \cite{BD} (following an
earlier work \cite{DS}). Roughly speaking, a $^L G$-oper on $X$ is a
principal $^L G$-bundle on the punctured disc (or a smooth complex
curve), equipped with a connection and a reduction to a Borel subgroup
of $^L G$, satisfying a certain transversality condition with respect
to the connection (we recall the precise definition in \secref{opers}
below).

This description of the center was used in \cite{F:faro} to show that
the Gaudin algebra is in fact isomorphic to the algebra of functions
on the space of $^L G$-opers on $\pone$ with regular singularities at
the points $z_1,\ldots,z_N$ and $\infty$. This implies that the joint
eigenvalues of the Gaudin algebra on any module $M_1 \otimes \ldots
\otimes M_N$ over $U(\g)^{\otimes N}$ are encoded by opers on $\pone$
with regular singularities.

\subsection{Non-highest weight representations and irregular
  singularities}

In this paper we pursue further the connection between affine
Kac--Moody algebras and quantum integrable systems. Recall from
\cite{FFR,F:faro} that the action of the Gaudin algebra on the tensor
product $M_1 \otimes \ldots \otimes M_N$ of $\g$-modules comes about
via its realization as the {\em space of coinvariants} of the tensor
product of the induced modules $\M_1,\ldots,\M_N$ over $\ghat_{\ka_c}$
of critical level, where we use the notation $\M = \Ind_{\g[[t]]
\oplus \C{\mb 1}}^{\ghat_{\ka_c}} M$.

The notion of the space of coinvariants (which is the dual space to
the so-called space of conformal blocks) comes from conformal field
theory (see, e.g., \cite{FB} for a general definition). More
precisely, we consider the curve $\pone$ with the marked points
$z_1,\ldots,z_N$ and $\infty$. We attach the $\ghat_{\ka_c}$-modules
$\M_1,\ldots,\M_N$ to the points $z_1,\ldots,z_N$ and another
$\ghat_{\ka_c}$-module $\M_\infty = \Ind_{\g \otimes t\C[[t]] \oplus
\C {\mb 1}} \C$ to the point $\infty$. The corresponding space of
coinvariants is isomorphic to $M_1 \otimes \ldots \otimes M_N$. It is
shown in \cite{FFR,F:faro} that the center $Z(\ghat)$ acts on this
space by functoriality, and its action factors through the Gaudin
subalgebra of $U(\g)^{\otimes N}$.

In this construction the $\ghat_{\ka_c}$-modules $\M_i$ satisfy an
important property: they are {\em highest weight modules} (provided
that the $M_i$'s are highest weight $\g$-modules; in general, $\M_i$
is generated by vectors annihilated by the Lie subalgebra $\g \otimes
t \C[[t]] \subset \ghat_{\ka_c}$). The representation theory of affine
Kac--Moody algebras has up to now almost exclusively been concerned
with highest weight representations, such as Verma modules, Weyl
modules and their irreducible quotients.

But these are by far not the most general representations of
$\ghat_{\ka_c}$. It is more meaningful to consider a larger category
of all {\em smooth} representations. Those are generated by vectors
annihilated by the Lie subalgebra $\g \otimes t^N\C[[t]]$ for some $N
\in \Z_+$.

The main question that we address in this paper is the following:

\medskip


\noindent{\em What quantum integrable systems correspond to
non-highest weight representations?}


\medskip

It is instructive to consider the following analogy: highest weight
representations are like differential equations with the mildest
possible singularities, namely, regular singularities, while
non-highest weight representations are like equations with irregular
singularities. Actually, this is much more than an analogy. The point
is that the action of the center $Z(\ghat)$ on the spaces of
coinvariants corresponding to non-highest weight representations of
$\ghat_{\ka_c}$ factors through the algebra of functions on the space
of {\em opers on $\pone$ with irregular singularities at
$z_1,\ldots,z_N$ and $\infty$}. The order of the pole at each point is
determined by the "depth" of the corresponding $\ghat_{\ka_c}$-module:
if it is generated by vectors annihilated by the Lie subalgebra $\g
\otimes t^N\C[[t]]$, then the order of the pole at this point is less
than or equal to $N$. This may be summarized by the following diagram:

$$
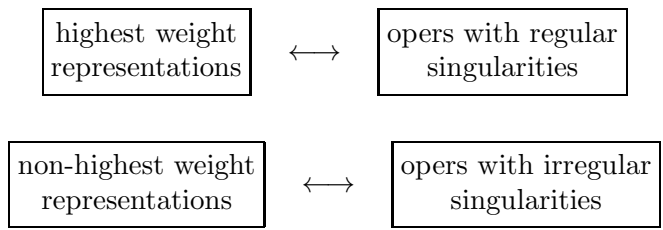
\boxed{\begin{matrix} \text{highest weight} \\
\text{representations} \end{matrix}} \quad \longleftrightarrow \quad
\boxed{\begin{matrix} \text{opers with regular} \\
    \text{singularities} \end{matrix}}
$$

\bigskip

$$
\boxed{\begin{matrix} \text{non-highest weight} \\
\text{representations} \end{matrix}} \quad \longleftrightarrow \quad
\boxed{\begin{matrix} \text{opers with irregular} \\
    \text{singularities} \end{matrix}}
$$

\bigskip

Motivated by this picture, we call the corresponding quantum
integrable systems the generalized {\em Gaudin models with irregular
singularities}. Their classical limits may be identified with a class
of Hitchin integrable systems on the moduli spaces of bundles on
$\pone$ with level structures at the points $z_1,\ldots,z_N,\infty$
and with Higgs fields having singularities at those points of orders
equal to the orders of the level structures. Classical integrable
systems of this type have been considered in
\cite{Beauville,Bot,Markman,DM,ER}.

Thus, using non-highest weight representations, we obtain more general
spaces of coinvariants (and conformal blocks) than those considered
previously. At the critical level this yields interesting commutative
subalgebras and quantum integrable systems, and away from the critical
level, interesting systems of differential equations.

In particular, the Knizhnik-Zamolodchikov (KZ) equations may be
obtained as the differential equations on the conformal blocks
associated to highest weight representations of $\ghat$ away from the
critical level (see, e.g., \cite{FFR}). On the other hand, there is a
flat connection constructed by J. Millson and V. Toledano Laredo
\cite{MTL,TL}, and, independently, C. De Concini (unpublished; a
closely related connection was also considered in
\cite{FMTV}). Following \cite{Boalch}, we will call it the DMT
connection. We expect that this connection may also be obtained from a
system of differential equations on the conformal blocks, but
associated to non-highest weight representations of $\ghat$, away from
the critical level. An indication that this is the case comes from the
work of P. Boalch \cite{Boalch}, in which the quasi-classical limit of
the DMT connection was related to certain isomonodromy equations, and
the paper \cite{BF} where it was shown that isomonodromy equations of
this type often arise as quasi-classical limits of equations on
conformal blocks (for example, the so-called Schlesinger isomonodromy
equations arise in the quasi-classical limit of the KZ equations, see
\cite{Res}). We hope that our results on the generalized Gaudin models
will help elucidate the relation between the DMT connection and
conformal field theory.

\subsection{The shift of argument subalgebra and the DMT Hamiltonians}

What are the simplest Gaudin models with irregular singularities? They
are obtained by allowing regular singularities at all but one point,
where we allow a pole of order two. It is convenient to take $\infty
\in \pone$ as this special point. The corresponding Gaudin algebra is
then a commutative subalgebra of $U(\g)^{\otimes N} \otimes S(\g)$,
where $S(\g)$ is the symmetric algebra of $\g$. Here $S(\g)$ arises as
the universal enveloping algebra of the commutative Lie algebra $\g
\otimes t\C[[t]]/t^2\C[[t]]$. This algebra naturally acts on the
simplest non-highest weight $\ghat$-module $\Ind_{\g \otimes
t^2\C[[t]] \oplus \C{\mb 1}}^{\ghat_{\ka_c}} \C$ attached to the point
$\infty$ (by endomorphisms commuting with the action of
$\ghat_{\ka_c}$). We may then specialize this algebra at a point $\chi
\in \g^* = \on{Spec} S(\g)$. As the result, we obtain a commutative
subalgebra of $U(\g)^{\otimes N}$ depending on $\chi \in \g^*$.

Recently, this algebra has been constructed by L. Rybnikov
\cite{Ryb2}, using the results of \cite{FF:gd,F:wak} on the center at
the critical level and the method of \cite{FFR}. Thus, the
construction of \cite{Ryb2} is essentially equivalent to our
construction, applied in this special case.

Consider in particular the case when $N=1$. Then the corresponding
Gaudin algebra, which we denote by ${\mc A}_\chi$, is a subalgebra of
$U(\g)$. It was shown in \cite{Ryb2} that for regular semi-simple
$\chi$ the algebra ${\mc A}_\chi$ is a quantization of the ``shift of
argument'' subalgebra $\ol{\mc A}_\chi$ of the associated graded
algebra $S(\g) = \on{gr} U(\g)$ (we show below that this is true for
all regular $\chi$). The subalgebra $\ol{\mc A}_\chi$, introduced by
A.S. Mishchenko and A.T. Fomenko in \cite{MF} (see also \cite{Man}),
is the Poisson commutative subalgebra of\footnote{Here and below, for
an affine algebraic variety $Z$ we denote by $\on{Fun} Z$ the algebra
of regular functions on $Z$.} $S(\g) = \on{Fun} \g^*$ generated by the
derivatives of all orders in the direction of $\chi$ of all elements
of $\on{Inv} \g^* = (\on{Fun} \g^*)^{\g}$, the algebra of invariants
in $\on{Fun} \g^*$. The algebra ${\mc A}_\chi$ is the quantization of
$\ol{\mc A}_\chi$ in the sense that $\on{gr} {\mc A}_\chi = \ol{\mc
A}_\chi$.

We note that the problem of quantization of $\ol{\mc A}_\chi$ was
posed by E.B. Vinberg \cite{Vin}. Such a quantization has been
previously constructed for $\g$ of classical types in \cite{NO} (using
twisted Yangians) and for $\g=\sw_n$ in \cite{Tar} (using the
symmetrization map) and \cite{CT} (using explicit
formulas).\footnote{in \cite{Ryb1} it was shown that, if it exists, a
quantization ${\mc A}_\chi$ is unique for generic $\chi$}

Our general results on the structure of the Gaudin algebras identify
${\mc A}_\chi$ with the algebra of functions on the space $\on{Op}_{^L
G}(\pone)_{\pi(-\chi)}$ of $^L G$-opers on $\pone$ with regular
singularity at $0 \in \pone$ and irregular singularity, of order $2$,
at $\infty \in \pone$, with the most singular term $-\chi$ (where
$\chi$ is any regular element of $\g^*$). This means that a joint
generalized eigenvalue of the quantum shift of argument subalgebra
${\mc A}_\chi \subset U(\g)$, for regular semi-simple $\chi$, on any
$\g$-module is encoded by a point in $\on{Op}_{^L
G}(\pone)_{\pi(-\chi)}$.

Note that $\ol{\mc A}_\chi$ is a graded subalgebra of $\on{Fun} \g^*$
with respect to the standard grading, so we have $\ol{\mc A}_\chi =
\bigoplus_{i\geq 0} \ol{\mc A}_{\chi,i}$. It is easy to see that the
degree $1$ piece $\ol{\mc A}_{\chi,1}$ is the Cartan subalgebra $\h
\subset \g$ (realized as degree one polynomial functions on $\g^*$)
which is the centralizer of $\chi$. The degree two piece $\ol{\mc
A}_{\chi,2}$ was determined in \cite{Vin}: it is spanned by elements
of the form
$$
\ol{T}_\ga(\chi) = \sum_{\al \in \De_+}
\frac{(\al,\gamma)(\al,\al)}{(\al,\chi)} e_\al f_\al,
\qquad \ga \in \h^*,
$$
where $\al \in \De_+$ is the set of positive roots of $\g$ with
respect to a Borel subalgebra $\bb$ containing $\h$, and $e_\al \in
\g_\al, f_\al \in \g_{-\al}$ are generators of an $\sw_2$ triple
corresponding to $\al$ (we also use in this formula an inner product
on $\h^*$ corresponding to a non-degenerate invariant inner product on
$\g$). Since we could not find the proof of this result in the
literature, we show that the elements $\ol{T}_\ga(\chi)$ indeed belong
to $\ol{\mc A}_\chi$ in \secref{DMT}.

Let now $T_\ga(\chi)$ be the element of $U(\g)$ given by the formula
\begin{equation}    \label{T ham}
T_\ga(\chi) = \frac{1}{2} \sum_{\al \in \De_+}
\frac{(\al,\gamma)(\al,\al)}{(\al,\chi)} (e_\al f_\al + f_\al e_\al),
\qquad \ga \in \h^*,
\end{equation}
so that the symbol of $T_\ga(\chi)$ is equal to
$\ol{T}_\ga(\chi)$. These operators coincide with the connection
operators of the DMT flat connection on $\h_{\on{reg}} \simeq
\h^*_{\on{reg}}$ mentioned above. The flatness of this connection
implies that these operators commute in $U(\g)$. Thus, we obtain that
the quadratic part of the Poisson commutative algebra $\ol{\mc
A}_\chi$ may be quantized by the DMT operators. It then follows from
the results of \cite{Tarasov} that ${\mc A}_\chi$ is a maximal
commutative subalgebra of $U(\g)$ containing $\h$ and the operators
$T_\ga(\chi)$.

To summarize, we find that the commutative algebra corresponding to
the ordinary Gaudin model, and the quantum shift of argument
subalgebra ${\mc A}_\chi$, arise as special cases of the general
construction of Gaudin models with irregular singularities that we
propose in this paper. In particular, the Gaudin operators and the DMT
operators arise as the quadratic generators of the corresponding
generalized Gaudin algebras.

Note that they both come from flat holomorphic connections, the KZ
connection and the DMT connection, respectively. As shown in
\cite{TL}, the KZ and DMT connections are dual to each other in the
case of ${\mathfrak g}{\mathfrak l}_N$. Therefore it is natural to
expect that the corresponding Gaudin algebras are also dual in this
case (some results in this direction are obtained in \cite{MTV}). We
hope to discuss this question elsewhere.

\subsection{Diagonalizing quantum Hamiltonians}
\label{diagonalizing}

The next step is to consider the problem of diagonalization of the
generalized Gaudin algebras on various representations.

For example, given a $\g$-module $M$, we can try to find joint
eigenvectors and eigenvalues of the algebra ${\mc A}_\chi$ acting on
$M$, where $\chi$ is a regular semi-simple. We know from the above
description of the spectrum that each joint eigenvalue of ${\mc
A}_\chi$ on any $\g$-module is encoded by a $^L G$-oper on $\pone$
which belongs to $\on{Op}_{^L G}(\pone)_{\pi(-\chi)}$. However, we
would like to know which $^L G$-opers may be realized on a given
$\g$-module. For example, we show that if $M$ is generated by a
highest weight vector (for example, if $M$ is a Verma module), then
the $^L G$-opers arising from the joint eigenvalues on $M$ have a
fixed residue at the point $0 \in \pone$ which is determined by the
highest weight (this is similar to what happens in the Gaudin model,
see \cite{F:faro}).

The most interesting case is when $M$ is an irreducible
finite-dimensional $\g$-module. It was proved in \cite{F:faro} that
in the case of the Gaudin model the joint eigenvalues of the Gaudin
algebra on the tensor product of finite-dimensional irreducible
$\g$-modules are encoded by the $^L G$-opers on $\pone$ with regular
singularities at $z_1,\ldots,z_N$ and $\infty$, with prescribed
residues at those points (determined by the highest weights of our
modules), and with {\em trivial monodromy}. Moreover, it was
conjectured in \cite{F:faro} that this sets up a bijection between the
joint spectrum of the Gaudin algebra on this tensor product and this
space of opers. This conjecture was motivated by the geometric
Langlands correspondence (see below).

By applying the same argument in the irregular case, we show that the
eigenvalues of ${\mc A}_\chi$ (for a regular semi-simple $\chi$) on an
irreducible finite-dimensional $\g$-module $V_\la$ are encoded by $^L
G$-opers on $\pone$ with irregular singularity of order $2$ with the
most singular term $-\chi$ at $\infty$, with regular singularity with
residue $\la$ at $0$ and trivial monodromy.

Let $\on{Op}_{^L G}(\pone)^\la_{\pi(-\chi)}$ be the set of such
opers. Then we obtain an injective map from the joint spectrum of
${\mc A}_\chi$ on $V_\la$ to $\on{Op}_{^L
G}(\pone)^\la_{\pi(-\chi)}$. We conjecture that this map is a
bijection for any regular semi-simple $\chi$ (we also give a
multi-point generalization of this conjecture). We expect that for
generic $\chi$ the algebra ${\mc A}_\chi$ is diagonalizable on $V_\la$
and has simple spectrum (this has been proved in \cite{Ryb2} for
$\g=\sw_n$). If this is so, then our conjecture would imply that there
exists an eigenbasis of ${\mc A}_\chi$ in a $\g$-module $V_\la$
parameterized by the monodromy-free $^L G$-opers on $\pone$ with
prescribed singularities at two points.

In the case of the ordinary Gaudin model, there is a procedure for
diagonalization of the Gaudin Hamiltonians called {\em Bethe
Ansatz}. In \cite{FFR,F:faro} it was shown that this procedure can
also be understood in the framework of coinvariants of $\ghat$-modules
of critical level. We need to use a particular class of
$\ghat$-modules, called the {\em Wakimoto modules}.

Let us recall that the Wakimoto modules of critical level are
naturally parameterized by objects closely related to opers, which are
called {\em Miura opers} \cite{F:wak}. They may also be described more
explicitly as certain connections on a particular $^L H$-bundle
$\Omega^{-\rho}$ on the punctured disc, where $^L H$ is the Cartan
subgroup of $^L G$. The center acts on the Wakimoto module
corresponding to a Cartan connection by the Miura transformation of
this connection (see \cite{F:wak}). The idea of \cite{FFR} was to use
the spaces of conformal blocks of the tensor product of the Wakimoto
modules to construct eigenvectors of the generalized Gaudin
Hamiltonians. It was found in \cite{FFR} that the eigenvalues of the
Gaudin Hamiltonians on these vectors are encoded by the $^L G$-opers
which are obtained by applying the Miura transformation to certain
very simple Cartan connections on $\pone$.

In this paper we apply the methods of \cite{FFR,F:faro} to define an
analogue of Bethe Ansatz for the Gaudin models with irregular
singularities. In particular, we use Wakimoto modules of critical
level to construct eigenvectors of the generalized Gaudin algebras
(such as the algebra ${\mc A}_\chi$) on Verma modules and
finite-dimensional irreducible representations of $\g$.

\subsection{Connection to the geometric Langlands correspondence}
\label{connection to glc}

One of the motivations for studying the generalized Gaudin systems and
their connection to opers comes from the geometric Langlands
correspondence. Here we give a very rough outline of this connection.

The ramified geometric Langlands correspondence proposed in
\cite{FG:local} (see \cite{F:ram} and the last chapter of
\cite{F:book} for an exposition) assigns to holomorphic $^L G$-bundles
with meromorphic connections on a Riemann surface $X$, certain
categories of Hecke eigensheaves on the moduli stacks of $G$-bundles
on $X$ with level structures\footnote{Here by a level structure on a
$G$-bundle on a curve $X$ of order $m$ at a point $x \in X$ we
understand a trivialization of the bundle on the $(m-1)$st
infinitesimal neighborhood of $x$.} at the positions of the poles of
the connection of the orders equal to the orders of the poles of the
connection (in the case of regular singularities, we may choose
instead a parabolic structure, which is a reduction of the fiber of
the bundle to a Borel subgroup of $G$). We note that recently the
geometric Langlands correspondence with ramification has been related
in \cite{GW} to the S-duality of four-dimensional supersymmetric
Yang-Mills theory.

Now, the point is that the generalized Gaudin algebra
$\ZZ^{(m_i),m_\infty}_{(z_i),\infty}(\g)$ (and its versions with
non-trivial characters $(\chi_i),\chi_\infty$) introduced in this
paper gives rise to a commutative algebra of global (twisted)
differential operators on moduli stack
$\on{Bun}_{G,(z_i),\infty}^{(m_i),m_\infty}$ of $G$-bundles on $\pone$
with level structures of orders $m_i$ at the points $z_i,
i=1,\ldots,N$, and $m_\infty$ at $\infty \in \pone$. In the case when
all $m_i=1$ this has been explained in detail in \cite{F:icmp},
following the seminal work \cite{BD} in the unramified case (which
applies to curves of arbitrary genus), and in general the construction
is similar.

In this paper we identify the Gaudin algebra
$\ZZ^{(m_i),m_\infty}_{(z_i),\infty}(\g)$ with the algebra of
functions on the space $\on{Op}_{^L
G}(\pone)^{(m_i),m_\infty}_{(z_i),\infty}$ of $^L G$-opers on $\pone$
with singularities of orders $m_i$ at the points $z_i, i=1,\ldots,N$,
and $m_\infty$ at $\infty$. Thus, each $f \in \on{Fun} \on{Op}_{^L
G}(\pone)^{(m_i),m_\infty}_{(z_i),\infty}$ gives rise to an element
of the Gaudin algebra $\ZZ^{(m_i),m_\infty}_{(z_i),\infty}(\g)$ and
hence a differential operator on
$\on{Bun}_{G,(z_i),\infty}^{(m_i),m_\infty}$, which we denote by
$D_f$.

For each point $\tau \in \on{Op}_{^L
G}(\pone)^{(m_i),m_\infty}_{(z_i),\infty}$ we may now write down the
following system of differential equations on
$\on{Bun}_{G,(z_i),\infty}^{(m_i),m_\infty}$:
\begin{equation}    \label{hecke}
D_f \cdot \Psi = f(\tau) \Psi, \qquad f \in \on{Fun}
\on{Op}_{^L G}(\pone)^{(m_i),m_\infty}_{(z_i),\infty}.
\end{equation}
This system defines a ${\mc D}$-module $\Delta_\tau$ on
$\on{Bun}_{G,(z_i),\infty}^{(m_i),m_\infty}$, and by adapting an
argument from \cite{BD} (which treated the unramified case), we obtain
that $\Delta_\tau$ is a Hecke eigensheaf, whose ``eigenvalue'' is the
$^L G$-local system on $\pone$ underlying the $^L G$-oper $\tau$ (see
\cite{F:icmp,F:rev,F:ram} for more details).

Thus, one can construct explicitly examples of Hecke eigensheaves on
$\on{Bun}_{G,(z_i),\infty}^{(m_i),m_\infty}$ corresponding to
irregular connections by using the generalized Gaudin algebras
introduced in this paper. These ${\mc D}$-modules provide us with a
useful testing ground for the geometric Langlands correspondence.

The philosophy of the geometric Langlands correspondence also gives us
insights into the structure of the spectra of the generalized Gaudin
algebras on tensor products of finite-dimensional modules. Namely, the
existence of such an eigenvector for a particular eigenvalue $\tau$
implies that the ${\mc D}$-module $\Delta_\tau$ should be in some
sense trivial. Therefore the $^L G$-local system underlying the $^L
G$-oper $\tau$ should also be trivial (i.e., monodromy-free). For
connections with regular singularities this is explained in
\cite{F:icmp,F:faro}, and in this paper we formulate conjectures to
this effect in the case of connections with irregular
singularities. This motivates, in particular, our conjecture, already
discussed above in \secref{diagonalizing}, that the spectrum of the
Gaudin algebra ${\mc A}_\chi$ on an irreducible finite-dimensional
$\g$-module is described by monodromy-free $^L G$-opers on $\pone$
with singularities at $0$ and $\infty$.

Another useful perspective is provided by the {\em separation of
variables} method pioneered by E. Sklyanin \cite{Skl}. The connection
between the geometric Langlands correspondence and the separation of
variables in the Gaudin model is discussed in detail in
\cite{F:icmp}. The point is that, roughly speaking, the Hecke property
of the ${\mc D}$-module corresponding to the system \eqref{hecke} is
reflected in the existence of separated variables for the system
\eqref{hecke}. These separated variables have been discovered by
Sklyanin \cite{Skl} in the case of $\g=\sw_2$ and regular
singularities. One can approach more general Gaudin models with
irregular singularities in a similar way. However, this subject is
beyond the scope of the present paper.

\subsection{Plan of the paper}

The paper is organized as follows. In \secref{construction} we present
a general construction of quantum Hamiltonians using spaces of
coinvariants and the center $\zz(\ghat)$ of the vertex algebra $\V_0$
associated to the affine Kac--Moody algebra $\ghat$ at the critical
level. We define the universal Gaudin algebra as a commutative
subalgebra of $U(\g[[t]])^{\otimes N} \otimes U(t\g[[t]])$ and its
various quotients. This generalizes the construction of the Gaudin
algebra from \cite{FFR,F:faro} which appears as a special
case. Another special case is a commutative subalgebra ${\mc A}_\chi$
of $U(\g)$. In \secref{graded} we describe the associated graded
algebras of the Gaudin algebras constructed in \secref{construction}
and relate them to the generalized Hitchin systems. In particular, we
identify the associated graded algebra of ${\mc A}_\chi$ and the
shift of argument subalgebra $\ol{\mc A}_\chi$ of $S(\g) =
\on{Fun} \g^*$ for regular $\chi$ (this has previously been done in
\cite{Ryb2} for regular semi-simple $\chi$).

Next, we wish to identify the spectra of the Gaudin algebras with the
appropriate spaces of $^L G$-opers on $\pone$. We start by collecting
in \secref{opers} various results on opers. Then we recall in
\secref{spectrum} the isomorphism between $\zz(\ghat)$ and the algebra
of functions on the space of $^L G$-opers on the disc from
\cite{FF:gd,F:wak}. After that we identify in \secref{spectrum} the
spectra of the universal Gaudin algebra and its quotients with various
spaces of opers on $\pone$ with singularities at finitely many
points. In particular, we identify the spectrum of the algebra ${\mc
A}_\chi$ with the space of $^L G$-opers on $\pone$ with singularities
of orders $1$ and $2$ at two marked points. We show, in the same way
as in \cite{F:faro}, that joint eigenvalues of the algebra ${\mc
A}_\chi$ (and its generalizations) on finite-dimensional
representations correspond to opers with trivial monodromy. Finally,
in \secref{bethe1} we study the problem of diagonalization of the
generalized Gaudin algebras. Applying the results of
\cite{FFR,F:faro}, we develop the Bethe Ansatz construction of
eigenvectors of the Gaudin algebras, such as ${\mc A}_\chi$ and its
multi-point generalizations, on tensor products of Verma modules and
finite-dimensional representations of $\g$.

\bigskip

\noindent{\bf Acknowledgments.} We thank L. Rybnikov for useful
comments on a draft of this paper. We are also grateful to
B. Kostant for providing a copy of his unpublished manuscript
\cite{Ko2}.

\section{Construction of generalized Gaudin algebras}
\label{construction}

In this section we introduce the universal Gaudin algebra and its
various quotients by using the coinvariants construction for affine
Kac--Moody algebras from \cite{FFR,F:faro}.

\subsection{Affine Kac--Moody algebra and its modules}\label{aff km}

Let $\g$ be a simple Lie algebra. Recall that the space of invariant
inner products on $\g$ is one-dimensional. Choose a non-zero
element $\ka$ in this space.

The affine Kac--Moody algebra $\ghat_{\ka}$ is the (universal)
extension of the Lie algebra $\g\ppart=\g \otimes \C\ppart$ by the
one-dimensional center $\C {\mb 1}$:
\begin{equation}\label{central extension}
0 \to \C {\mb 1} \to \ghat_{\ka} \to \g \ppart \to 0.
\end{equation}
The commutation relations in $\ghat_{\ka}$ read
\begin{equation}
\label{comm rel}
[A \otimes f(t),B \otimes g(t)] = [A,B] \otimes fg - \ka(A,B)
\on{Res}_{t=0} f dg \cdot {\mb 1}.
\end{equation}
The $\ghat_{\ka}$-modules on which ${\mb 1}$ acts as the identity will
be referred to as modules of {\em level} $\ka$.

We introduce the following notation: for a Lie algebra $\l$,
a subalgebra $\m\subset\l$ and an $\m$-module $M$,
we denote by
$$\Ind_{\m}^\l M = U(\l) \underset{U(\m)}\otimes M$$
the $\l$-module induced from $M$.

Let $\G_+$ be the Lie subalgebra $\g [[t]] \oplus \C {\mb
1}$ of $\ghat_{\ka}$. Given a $\g [[t]]$-module $M$, we
extend it to a $\G_+$-module by making ${\mb 1}$ act as the
identity. Denote by ${\mathbb M}_\ka$ the corresponding induced
$\ghat_{\ka}$-module
$${\mathbb M}_\ka = \Ind_{\G_+}^{\ghat_{\ka}} M$$ of level $\ka$.

For example, let $\bb\subset\g$ be a Borel subalgebra, $\n=[\bb,\bb]$
its nilpotent radical and $\h=\bb/\n$. For any $\la\in\h^*$, let
$\C_\la$ be the one-dimensional $\bb$-module on which $\h$ acts by the
character $\la: \h \to \C$ and $M_{\la}$ be the Verma module over $\g$
of highest weight $\la$,
$$
M_\la = \Ind_{\bb}^\g \C_\la.
$$
The corresponding induced $\ghat_{\ka}$-module $\M_{\la,\ka}$ is the
{\em Verma module} over $\ghat_{\ka}$ of level $\ka$ with highest
weight $\la$.

For a dominant integral weight $\la \in \h^*$, denote by $V_\la$ the
irreducible finite-dimensional $\g$-module of highest weight $\la$.
The corresponding induced module
$$
\V_{\la,\ka} = \Ind_{\G_+}^{\ghat_{\ka}} V_\la
$$
is called the {\em Weyl module} of level $\ka$ with highest weight
$\la$ over $\ghat_{\ka}$.

A $\ghat_{\ka}$-module $R$ is called a {\em highest weight module}
if it is generated by a {\em highest weight vector}, that is a $v\in R$ such
that $\wh\n_+ \cdot v = 0$, where
$$
\wh\n_+ = (\n\otimes 1) \oplus t\g[[t]],
$$
and $\h \otimes 1$ acts on $v$ through a linear functional $\la: \h
\to \C$. In this case we say that $v$ (or $M$) has highest weight $\la$.

For example, both $\M_{\la,\ka}$ and $\V_{\la,\ka}$ (for a dominant
integral weight $\la$) are highest weight modules, with highest weight
$\la$.

Here is an example of a non-highest weight module, which we will use
in this paper. Let $\chi: \g \to \C$ be a linear functional. Observe
that the composition
$$
t\g[[t]] \to t\g[[t]]/t^2\g[[t]] \simeq \g
\overset\chi\longrightarrow \C
$$
defines a one-dimensional representation $\C_\chi$ of the Lie algebra
$t\g[[t]]$. We extend it to the direct sum $t\g[[t]]\oplus \C{\mb 1}$ by
making ${\mb 1}$ act as the identity. Now set
\begin{equation}    \label{I1}
\I_{1,\chi,\ka} = \Ind_{t\g[[t]] \oplus \C{\mb 1}}^{\ghat_{\ka}}
\C_\chi.
\end{equation}
Note that $\I_{1,\chi,\ka}$ may also be realized as
$$
\Ind_{\G_+}^{\ghat_{\ka}} I_{1,\chi}, \qquad \on{where} \quad
I_{1,\chi} = \Ind_{t\g[[t]]}^{\g[[t]]} \C_\chi.
$$

\subsection{Spaces of coinvariants} \label{coinvariants}

We define spaces of coinvariants of the tensor products of
$\ghat_{\ka}$-modules, associated to the projective line and a
collection of marked points. These spaces and their duals, called
spaces of conformal blocks, arise naturally in conformal field theory
(for more on this, see, e.g., the book \cite{FB}).

Consider the projective line $\pone$ with a global coordinate $t$
and $N$ distinct finite points $z_1,\ldots,z_N \in \pone$. In the
neighborhood of each point $z_i$ we have the local coordinate
$t-z_i$ and in the neighborhood of the point $\infty$ we have the
local coordinate $t^{-1}$. Set
$$\widetilde{\g}(z_i)=\g\ppzi,\qquad\wt\g(\infty)=\g\ppinf$$ and let
$\hatg_\ka(z_i),\hatg_\ka(\infty)$ be the corresponding central
extensions \eqref{central extension} respectively. Let $\G_N$ be the
extension of the Lie algebra $\bigoplus_{i=1}^N \widetilde{\g}(z_i)
\oplus \wt\g(\infty)$ by a one-dimensional center $\C {\mb 1}$ whose
restriction to each summand $\widetilde{\g}(z_i)$ or $\wt\g(\infty)$
coincides with $\hatg_\ka(z_i)$ or $\hatg_\ka(\infty)$.  Thus, $\G_N$
is the quotient of $\bigoplus_{i=1}^N\hatg_\ka(z_i)\oplus
\hatg_\ka(\infty)$ by the subspace spanned by
$\bfone_{z_i}-\bfone_{z_j}$, $1\leq i<j\leq N$ and
$\bfone_{z_i}-\bfone_\infty$.

In what follows, we will consider exclusively {\it smooth}
$\g[[t]]$-modules. By definition, a $\g[[t]]$-module $M$ is smooth if
for any $v\in M$ we have $t^k \g[[t]] v=0$ for sufficiently large
$k\in\IZ_+$. Equivalently, these are the $\g[[t]]$-modules such that
for any $v \in M$ the map $\g \to M, x \mapsto x \cdot v$ is
continuous with respect to the $t$-adic topology on $\g[[t]]$ and the
discrete topology on $M$. Note that if $M$ is a smooth finitely
generated $\g[[t]]$-module, the action of $\g[[t]]$ on $M$ factors
through the Lie algebra $\g[[t]]/t^k\g[[t]]$ for some $k\in\IZ_+$.

Suppose we are given a collection $M_1,\ldots,M_N$ and $M_\infty$ of
$\g[[t]]$-modules. Then the Lie algebra $\G_N$ naturally acts on the
tensor product of the induced $\ghat_{\ka}$-modules
$\bigotimes_{i=1}^N {\mathbb M}_i \otimes {\mathbb M}_\infty$ (in
particular, ${\mb 1}$ acts as the identity).

Let $\g_{\zn}=\g_{z_1,\ldots,z_N}$ be the Lie algebra of $\g$-valued
regular functions on $\pone\backslash\{ z_1, \ldots,$ $z_N,\infty \}$
(i.e., rational functions on $\pone$, which may have poles only at the
points $z_1,\ldots,z_N$ and $\infty$). Clearly, such a function can be
expanded into a Laurent power series in the corresponding local
coordinates at each point $z_i$ and at $\infty$. Thus, we obtain an
embedding
\begin{equation}\label{embedding}
\g_{\zn} \hookrightarrow \bigoplus_{i=1}^N
\widetilde{\g}(z_i) \oplus \wt\g(\infty).
\end{equation}
It follows from the residue theorem and formula \eqref{comm rel} that
the restriction of the central extension to the image of this
embedding is trivial. Hence \eqref{embedding} lifts to an embedding
$\g_{\zn}\to\G_N$.

Denote by $H(\M_{1,\ka},\ldots,\M_{N,\ka},\M_{\infty,\ka})$ the space
of {\em coinvariants} of $\bigotimes_{i=1}^N {\mathbb M}_{i,\ka}
\otimes {\mathbb M}_{\infty,\ka}$ with respect to the action of the
Lie algebra $\g_{\zn}$:
$$
H(\M_{1,\ka},\ldots,\M_{N,\ka},\M_{\infty,\ka}) = (\bigotimes_{i=1}^N
{\mathbb M}_{i,\ka} \otimes {\mathbb M}_{\infty,\ka})/\g_{\zn}.
$$

By construction, we have a canonical embedding of a $\g[[t]]$-module
$M$ into the induced $\ghat_{\ka}$-module $\M_\ka$:
$$x \in M \arr 1 \otimes x \in \M_\ka,$$ which commutes with the action
of $\g[[t]]$ on both spaces. Thus, we have an embedding
$$\bigotimes_{i=1}^N M_i \otimes M_\infty \hookrightarrow
\bigotimes_{i=1}^N \M_{i,\ka} \otimes \M_{\infty,\ka}.$$

The following result gives a description of
$H(\M_1,\ldots,\M_N,\M_\infty)$ in terms of coinvariants for the \fd
Lie algebra $\g\subset\g_{\zn}$. It will allow us to construct
quantum Hamiltonians acting on tensor products of $\g$-modules in
Section \ref{action} below.

\begin{lem}[\cite{FFR}, Lemma 1] \label{iso}
The composition of the above embedding and the projection
$$\bigotimes_{i=1}^N \M_{i,\ka} \otimes \M_{\infty,\ka}
\twoheadrightarrow H(\M_{1,\ka},\ldots,\M_{N,\ka},\M_{\infty,\ka})$$
gives rise to an isomorphism
$$
H(\M_{1,\ka},\ldots,\M_{N,\ka},\M_{\infty,\ka}) \simeq
(\bigotimes_{i=1}^N M_{i,\ka} \otimes M_{\infty,\ka})/\g.
$$
\end{lem}

\proof Fix a point $u\in\pone\setminus\{z_1,\ldots,z_N,\infty\}$ and
let $\g_{\zn}^0\subset\g_{\zn}$ be the ideal of $\g$-valued functions
vanishing at $u$. Then, $\g_{\zn}/\g_{\zn}^0\cong \g$, so that
$$H(\M_{1,\ka},\ldots,\M_{N,\ka},\M_{\infty,\ka})=
(\bigotimes_{i=1}^N\IM_{i,\ka} \otimes\IM_{\infty,\ka})/\g_{\zn}=
\left. \left((\bigotimes_{i=1}^N\IM_{i,\ka}
\otimes\IM_{\infty,\ka})/\g^0_{\zn}\right) \right/\g.$$
Since
$$\bigoplus_{i=1}^N\wt\g(z_i)\oplus\wt\g(\infty)=
\left(\bigoplus_{i=1}^N\g[[t-z_i]]\oplus\g[[t^{-1}]]\right)
\bigoplus \g_{\zn}^0$$
the $\hatg_N$-module
$$\bigotimes_{i=1}^N\IM_{i,\ka} \otimes\IM_{\infty,\ka} =
\ind^{\hatg_N}_
{\bigoplus_{i=1}^N\wt\g(z_i)\oplus\wt\g(\infty)\oplus\IC{\mathbf 1}}
\left(\bigotimes_{i=1}^N M_i \otimes M_\infty\right)$$
is freely generated over $U\g_{\zn}^0$ by $\bigotimes_{i=1}^N
M_i \otimes M_\infty$. The conclusion now follows. \qed

\subsection{The algebra of endomorphisms of $\IV_{0,\ka}$}

Let
$$\V_{0,\ka}=\ind_{\hatg_+}^{\hatg_\ka}\IC$$
be the {\em vacuum module} of level $\ka$, that is the Weyl module
corresponding to the highest weight $0$, and let $v_0$ be its generating
vector.

We will be concerned with the algebra $\End_{\hatg_\ka}(\IV_{0,\ka})$
of endomorphisms of $\IV_{0,\ka}$. Note that we may identify this
algebra with the space
$$\zz_\ka(\G) = \V_{0,\ka}^{\g[[t]]} $$ of $\g[[t]]$-invariant vectors
in $\V_{0,\ka}$. Indeed, a $\g[[t]]$-invariant vector $v$ gives rise
to an endomorphism of $\V_{0,\ka}$ commuting with the action of
$\ghat_{\ka}$ which sends $v_0$ to $v$. Conversely, any
$\ghat_{\ka}$-endomorphism of $\IV_{0,\ka}$ is uniquely determined
by the image of $v_0$, which necessarily belongs to $\zz_\ka(\G)$.
Thus, we obtain an isomorphism
$\zz_\ka(\G)\simeq\on{End}_{\ghat_{\ka}} (\V_{0,\ka})$ which gives
$\zz_\ka(\G)$ an algebra structure.

The space $\IV_{0,\ka}$ has the structure of a vertex algebra and it
is easy to see that $\zz_\ka(\G)$ is its center (see, e.g.,
\cite{FB}). This immediately implies the following:

\begin{prop}\label{commut}
The algebra $\zz_\ka(\G)$ is commutative.
\end{prop}

Let $\hatg_-$ be the Lie subalgebra $t^{-1}\g[t^{-1}] \subset
\ghat_\ka$. Let us identify $\IV_{0,\ka}$ with $U(\hatg_-)$ as
$\hatg_-$-modules, with $\hatg_-$ acting on $U(\hatg_-)$ by
left multiplication.  This yields an embedding
\begin{equation}\label{eq:center}
\End_{\hatg_\ka}(\IV_{0,\ka}) \subset \End_{\hatg_-}(\IV_{0,\ka})\cong
U(\hatg_-)\opp
\end{equation}
where the latter acts on $\IV_{0,\ka}\cong U(\hatg_-)$ by right
multiplication.  Thus, $\zz_\ka(\G)$ may be realized as a subalgebra
of $U(\hatg_-)\opp$. This realization will be useful below.

Note that the Lie algebra $\Der(\IC[[t]])=\IC[[t]]\partial_t$ of
continuous derivations of $\IC[[t]]$ acts on $\g\ppart$ and leaves
$\g[[t]]$ invariant.  This action lifts to one on $\ghat_\ka$ since
the cocycle appearing in \eqref{comm rel} is invariant under
coordinate changes, and induces one on $\IV_{0,\ka}$ which leaves
$\IV_{0,\ka}^{\g[[t]]}$ invariant. Of particular relevance to us will
be the translation operator $T=-\partial_t$ which acts on
$\IV_{0,\ka}$ by
\begin{equation}\label{translation}
Tv_0=0,\qquad\qquad
[T,J_n]=-nJ_{n-1}.
\end{equation}
where, for $J\in\g$ and $n\in\IZ$ we denote $J\otimes t^n$ by $J_n$.

According to \thmref{center},(1), $\zz_\ka(\G) = \C$, unless $\ka$
takes a special value $\ka_c$, called the {\em critical level}, and so
a meaningful theory involving $\zz_\ka(\G)$ exists only for the
critical level. However, since the proofs of the results of the next
section do not require it, we will postpone specializing $\ka$ to
$\ka_c$ until \secref{universal}.

\subsection{Action of $\zz_\ka(\hatg)$ on coinvariants}\label{action}

In order to simplify our notation, from now on we will omit the
subscript $\ka$ in our formulas.

Fix a point $u\in\pone\setminus\{z_1,\ldots,z_N,\infty\}$. Recall that
$\IV_{0}$ denotes the vacuum module of $\hatg_\ka(u)$ and $v_0$ its
generating vector. We define below a canonical action of the
commutative algebra $\zz_\ka(\hatg)\cong\End_{\hatg_\ka(u)}(\IV_{0})$
on spaces of coinvariants.  This relies on the following well-known
result.

\begin{prop}\label{pr:propagation}
For any modules $\{\IM_{i}\}_{i=1}^N$ and $\IM_{\infty}$ of
$\hatg_\ka(z_i)$, $i=1,\ldots,N$ and $\hatg_\ka(\infty)$ respectively,
the map
$$\bigotimes_{i=1}^N \IM_{i} \otimes
\IM_{\infty}\longrightarrow \bigotimes_{i=1}^N
\IM_{i}\otimes\IM_{\infty}\otimes \IV_{0,\ka},\qquad
\bigotimes_{i=1}^N v_i\otimes v_\infty\longrightarrow \bigotimes
_{i=1}^N v_i\otimes v_\infty\otimes v_0$$ induces an
isomorphism
$$\jmath_u:H(\IM_{1},\ldots,\IM_{N},\IM_{\infty})
\longrightarrow
H(\IM_{1},\ldots,\IM_{N},\IM_{\infty},\IV_{0}).$$
\end{prop}
\proof Let $\g_{(z_i,u)}$ be the Lie algebra of $\g$-valued regular
functions on $\IC\setminus\{z_1,\ldots,z_N,u\}$. Then,
$\g_{(z_i,u)}=\g_{\zn}\oplus(t-u)^{-1}\g[(t-u)^{-1}]$, as a vector
space, so that
\begin{equation*}
\begin{split}
H(M_1,\ldots,M_N,M_\infty,\IV_{0})
&=
(\bigotimes_{i=1}^N\IM_i \otimes\IM_\infty\otimes\IV_{0})/
\g_{(z_i,u)}\\
&\simeq
(\bigotimes_{i=1}^N\IM_i)/\g_{\zn} \otimes \IV_{0}/
(t-u)^{-1}\g[(t-u)^{-1}] \\
&\cong
(\bigotimes_{i=1}^N\IM_i \otimes\IM_\infty)/\g_{\zn}
\end{split}
\end{equation*}
where the last isomorphism is due to the fact that $\IV_{0}$ is a free
module of rank $1$ over $U((t-u)^{-1}\g[(t-u)^{-1}])$. \qed

\smallskip

Proposition \ref{pr:propagation} yields a canonical action of the
algebra $\End_{\hatg_\ka(u)}(\IV_{0})$ on the space $H(\IM_1,
\ldots,\IM_N,\IM_\infty)$ given by
\begin{equation}\label{eq:action}
X \cdot [\bigotimes_{i=1}^N v_i\otimes v_\infty]=
\jmath_u^{-1}[\bigotimes_{i=1}^N v_i\otimes v_\infty\otimes X \cdot
  v_0],
\end{equation}
where $[\bigotimes_{i=1}^N v_i\otimes v_\infty]$ denotes the image of
$$\bigotimes_{i=1}^N v_i\otimes v_\infty\in\bigotimes_{i=1}^N\IM
_i\otimes\IM_\infty$$ in $H(\IM_1,\ldots,\IM_\infty)$.

Let us work out this action explicitly. Set
$$\hatg(u)_-=(t-u)^{-1}\g[(t-u)^{-1}]\subset \wt\g(u) \subset
\hatg_\ka(u)$$ and identify $\IV_{0}$ with $U(\hatg(u)_-)$ as
$\hatg(u)_-$-modules, where $\hatg(u)_-$ acts on $U(\hatg(u)_-)$ by
left multiplication.  As explained in the previous section, this
yields an embedding
\begin{equation}\label{eq:center1}
\End_{\hatg_\ka(u)}(\IV_{0})\subset
\End_{\hatg(u)_-}(\IV_{0})\cong
U(\hatg(u)_-)\opp
\end{equation}
where the latter acts on $\IV_{0}\cong U(\hatg(u)_-)$ by right
multiplication.

Notice next that the $\tildeg(u)$-component of the embedding
$$\g_{(z_i,u)}\hookrightarrow\bigoplus_{i=1}^N
\wt\g(z_i)\oplus\wt\g(u)\oplus\wt\g(\infty)$$ has a section over
$\ghat(u)_-$ given by regarding an element of $\ghat(u)_-$ as a
regular function on $\IP^1\setminus\{u,\infty\}$. Composing with
Laurent expansion at the points $z_i$ and $\infty$ yields a Lie
algebra homomorphism
\begin{equation}\label{eq:sigmau}
\tau_{u,(z_i)}:
\ghat(u)_-\hookrightarrow
\bigoplus_{i=1}^N\tildeg(z_i)\oplus\tildeg(\infty)
\end{equation}

Let $\Phi_{u,(z_i)}$ be the extension of $-\tau_{u,(z_i)}$ to an
{\it anti}-homomorphism
$$\Phi_{u,(z_i)}:U(\ghat(u)_-)\longrightarrow U\hatg_N.$$

\begin{lem}\label{le:swap}
Let $\IM_i$, $\IM_\infty$ and $\IM_u$ be representations of
$\hatg_\ka(z_i)$, $\hatg_\ka(\infty)$ and $\hatg_\ka(u)$
respectively. The following holds in
$H(\IM_1,\ldots,\IM_N,\IM_\infty,\IM_u)$ for any
$v_i\in\IM_i,v_\infty\in\IM _\infty$, $v_u\in\IM_u$ and $X\in
U(\ghat(u)_-)$ we have
$$[\bigotimes_{i=1}^N v_i\otimes v_\infty\otimes Xv_u]=
[\Phi_{u,(z_i)}(X)(\bigotimes_{i=1}^N v_i\otimes v_\infty)\otimes
  v_u].$$
\end{lem}

\proof It suffices to prove this for $X\in\ghat(u)_-$. In that case it
follows because $X+\tau_{u,(z_i)}(X)\in\g_{(z_i,u)}$, so that
$$[\tau_{u,(z_i)}(X)\cdot (\bigotimes_{i=1}^N v_i\otimes
v_\infty)\otimes v_u] + [\bigotimes_{i=1}^N v_i\otimes v_\infty\otimes
X \cdot v_u]=0.$$ \qed

\smallskip

\begin{prop}\label{pr:swap}
The action of $$X\in\End_{\hatg_\ka(u)}(\IV_{0})\subset
U(\hatg(u)_-)\opp$$ on $H(\IM_1,\ldots,\IM_N,\IM_\infty)$ defined by
\eqref{eq:action} is given by
$$X[\bigotimes_{i=1}^N v_i\otimes v_\infty]=
[\Phi_{u,(z_i)}(X)(\bigotimes_{i=1}^N v_i\otimes v_\infty)]$$
\end{prop}
\proof This follows at once from \eqref{eq:action} and Lemma
\ref{le:swap}.
\qed

\smallskip

For later use, we work out the homomorphism \eqref{eq:sigmau}.

\begin{prop}\label{pr:sigmau}
The following holds for any $J\in\g$, $m\geq 0$ and $n>0$
\begin{align*}
\tau_{u,(z_i)}(J_m)
&=
\sum_{i=1}^N \sum_{p=0}^m (-1)^{m-p}
\begin{pmatrix}m\\p\end{pmatrix} (u-z_i)^{m-p}J^{(i)}_p\\
&+
\sum_{p=0}^m (-1)^{m-p}
\begin{pmatrix}m\\p\end{pmatrix} u^{m-p}J^{(\infty)}_{-p}, \\
\tau_{u,(z_i)}(J_{-n})
&=
\frac{\partial_{u}^{n-1}}{(n-1)!}\left(
-\sum_{i=1}^N \sum_{p\geq 0} \frac{J^{(i)}_p}{(u-z_i)^{p+1}}
+\sum_{p\geq 0}u^p J^{(\infty)}_{p+1} \right),
\end{align*}
where $J^{(i)}_p=J\otimes (t-z_i)^p\in\wt\g(z_i)$ and $J^{(\infty)}_p=
J\otimes t^{-p}\in\wt\g(\infty)$.
\end{prop}
\proof The first identity follows from the fact that
$$(t-u)^m=
\sum_{p=0}^m\begin{pmatrix}m\\p\end{pmatrix}(t-z_i)^p(z_i-u)^{m-p}$$
and
$$(t-u)^m=
\sum_{p=0}^m (-1)^{m-p}\begin{pmatrix}m\\p\end{pmatrix}t^p u^{m-p}$$
The second one follows from
$$\frac{1}{t-u}=
-\frac{1}{u-z_i}\frac{1}{1-\frac{t-z_i}{u-z_i}}=
-\sum_{m\geq 0}\frac{(t-z_i)^m}{(u-z_i)^{m+1}}$$
and 
$$\frac{1}{t-u}=
\frac{1}{t}\sum_{m\geq 0}\left(\frac{u}{t}\right)^m.$$
\qed

\subsection{The universal Gaudin algebra and its quotients}
\label{universal}

Now we identify $\wt\g(u)$ and $\bigoplus_{i=1}^N\wt\g(z_i)
\oplus\wt\g(\infty)$ with $\g\ppart$ and $\g\ppart^{\oplus (N+1)}$,
respectively, by putting $\wt\g(\infty)$ as the last component in the
direct sum. We also identify
$$
\hatg_- = t^{-1} \g[t^{-1}] \simeq
\hatg(u)_-=(t-u)^{-1}\g[(t-u)^{-1}]
$$
Then we obtain an anti-homomorphism
$$\Phi_{u,(z_i)}:
U(\ghat_-)\to
U(\g[[t]])^{\otimes N}\otimes U(t\g[[t]]).$$

By Proposition \ref{pr:sigmau}, $\Phi_{u,(z_i)}$ is given by the
formula
\begin{equation} \label{Psi}
\Phi_{u,(z_i)}(J^{a_1}_{-n_1} \ldots J^{a_m}_{-n_m})=
{\mb J}^{a_m}_{-n_m}(u) \ldots {\mb J}^{a_1}_{-n_1}(u),
\end{equation}
where
\begin{equation}\label{J u}
{\mb J}_{-n}(u)=
\frac{\partial_{u}^{n-1}}{(n-1)!}\left(
\sum_{i=1}^N \sum_{p\geq 0} \frac{J^{(i)}_p}{(u-z_i)^{p+1}}-
\sum_{p\geq 0}u^p J^{(\infty)}_{p+1}
\right).
\end{equation}

It is clear from the construction that $\Phi_{u,(z_i)}$ is independent
of $\ka$. What does depend on $\ka$ is the algebra $\zz_\ka(\ghat)$
which we realize as a subalgebra of $U(\ghat_-)\opp$, and to which we
then restrict the anti-homomorphism $\Phi_{u,(z_i)}$.

Note that $\g$ acts on $\hatg_-$ by adjoint action and that the
induced action on $U(\hatg_-)$ coincides with that on $\IV_{0}$
under the identification $\IV_{0}\cong U(\hatg_-)$. Since any
$v\in\zz(\G)\subset\V_{0}$ satisfies $\g\cdot v=0$, we find that
$\Phi_{u,(z_i)}(v)$ is invariant under the the diagonal action of $\g$
on $U(\g[[t]])^{\otimes N} \otimes U(t\g[[t]])$. Thus,
$\Phi_{u,(z_i)}$ restricts to an algebra homomorphism
\begin{equation}\label{univ hom}
\zz_\ka(\ghat) \to \left(U(\g[[t]])^{\otimes N}\otimes
U(t\g[[t]])\right)^\g
\end{equation}

At this point we quote \thmref{center},(1), according to which
$\zz_\ka(\ghat) = \C$, if $\ka \neq \ka_c$, where $\ka_c$ is the {\em
critical} invariant inner product on $\g$ defined by the
formula\footnote{Note that $\ka_c = - h^\vee \ka_0$, where $\ka_0$ is
the inner product normalized as in \cite{Kac} (so that the square
length of the maximal root is equal to $2$) and $h^\vee$ is the dual
Coxeter number.}
$$\ka_c(A,B) = - \frac{1}{2} \on{Tr}_\g \on{ad} A \on{ad} B.$$ On the
other hand, the center $\zz_{\ka_c}(\ghat)$ is non-trivial (see
\thmref{center},(2)).

Thus, the homomorphism \eqref{univ hom} is non-trivial only for
$\ka=\ka_c$. So from now on we will specialize $\ka$ to $\ka_c$ and
omit the symbol $\ka$ from most of our formulas. In particular, we
will write $\zz(\ghat)$ for $\zz_{\ka_c}(\ghat)$, $\IM$ for
$\IM_{\ka_c}$, and so on.

\begin{definition}
The {\em universal Gaudin algebra}
$$\ZZ_{(z_i),\infty}(\g)\subset
\left(U(\g[[t]])^{\otimes N}\otimes U(t\g[[t]])\right)^\g$$
associated to $\g$ and the set of points $z_1,\ldots,z_N\in\IC$ is the
image of $\zz(\hatg)$ under $\Phi_{u,(z_i)}$.
\end{definition}

Thus, $\ZZ_{(z_i),\infty}(\g)$ is a commutative algebra, and the
action of $\zz(\hatg)$ on coinvariants
$H(\IM_1,\ldots,\IM_N,\IM_\infty)$ factors through
$\ZZ_{(z_i),\infty}(\g)$ by Proposition \ref{pr:swap}.

We will show below that the algebra $\ZZ_{(z_i),\infty}(\g)$ is
independent of the chosen point
$u\in\pone\setminus\{z_1,\ldots,z_N,\infty\}$, which is the reason why
we suppress its notational dependence on $u$. This is easier to
establish in terms of the quotients
$\ZZ_{(z_i),\infty}^{(m_i),m_\infty} (\g)$ of $\ZZ_{(z_i),\infty}$
that we introduce presently.

Note first that the algebra $U(\g[[t]])^{\otimes N}\otimes
U(t\g[[t]])$ is a complete topological algebra, whose topology is
given as follows.  For each positive integer $m$, the Lie subalgebra
$t^m\g[[t]]\subset\g[[t]]$ is an ideal. The left ideal ${\mc I}_m$ it
generates in $U(\g[[t]])$ is therefore a two-sided ideal and the
quotient $U(\g[[t]])/{\mc I}_m$ is isomorphic to the universal
enveloping algebra $U(\g_m)$, where $\g_m=
\g[[t]]/t^m\g[[t]]$. Similarly, the subalgebra $t^m\g[[t]]\subset
t\g[[t]]$ generates a two-sided ideal $\ol{\mc I}_m$ of $U(t\g[[t]])$
and the quotient $U(t\g[[t]])/\ol{\mc I}_m$ is isomorphic to
$U(\ol\g_m)$, where $\ol\g_m = t\g[[t]]/t^m\g[[t]]$.

The topology on $U(\g[[t]])^{\otimes N}\otimes U(t\g[[t]])$ is defined
by declaring
$$
{\mc I}_{(m_i)} = \bigoplus_{i=1}^N U(\g[[t]])^{\otimes (i-1)}
\otimes{\mc I}_{m_i}\otimes U(\g[[t]])^{\otimes (N-i)}\otimes
U(t\g[[t]]) \bigoplus U(\g[[t]])^{\otimes N}\otimes \ol{\mc
I}_{m_\infty}
$$
to be the base of open neighborhoods of $0$. This
algebra is complete in this topology, and it is the inverse
limit
\begin{equation*}
\begin{split}
U(\g[[t]])^{\otimes N}\otimes U(t\g[[t]])
&=
\lim_{\longleftarrow} U(\g[[t]])^{\otimes N}\otimes U(t\g[[t]])/
{\mc I}_{(m_i)} \\
&=
\lim_{\longleftarrow}\bigotimes_{i=1}^N U(\g_{m_i})\otimes
U(\ol\g_{m_\infty}).
\end{split}
\end{equation*}

For any collection of positive integers $m_1,\ldots,m_N,m_\infty$, let
\begin{equation}    \label{Phi}
\Phi_{u,(z_i)}^{(m_i),m_\infty}:
U(\ghat_-)\to
\bigotimes_{i=1}^N U(\g_{m_i}) \otimes U(\ol\g_{m_\infty})
\end{equation}
be the composition of
$\Phi_{u,(z_i)}$ with the natural surjection
\begin{equation}\label{surjection}
U(\g[[t]])^{\otimes N} \otimes U(t\g[[t]])\twoheadrightarrow
\bigotimes_{i=1}^N U(\g_{m_i}) \otimes U(\ol\g_{m_\infty}).
\end{equation}
Thus, $\Phi_{u,(z_i)}^{(m_i),m_\infty}$ is obtained from the formulas
\eqref{Psi} and \eqref{J u} by setting $J^{(i)}_m=0$ for $m>m_i$ and
$J^{(\infty)}_m = 0$ for $m>m_\infty$.

Let
$$\ZZ_{(z_i),\infty}^{(m_i),m_\infty}(\g)\subset
\left(\bigotimes_{i=1}^N U(\g_{m_i})\otimes
U(\ol\g_{m_\infty})\right)^\g$$ be the image of $\zz(\hatg)$ under
$\Phi_{u,(z_i)}^{(m_i),m_\infty}$. The algebras $\ZZ_{(z_i),\infty}$
form an inverse system of commutative algebras, and
\begin{equation}\label{eq:inv limit}
\ZZ_{(z_i),\infty}(\g)=\lim_{\longleftarrow}
\ZZ_{(z_i),\infty}^{(m_i),m_\infty}(\g).
\end{equation}

Formulas \eqref{Psi}--\eqref{J u} show that if $u$ is considered as a
variable, $\Phi_{u,(z_i)}^{(m_i),m_ \infty}$ may be regarded as an
anti-homomorphism
$$\Phi_{u,(z_i)}^{(m_i),m_\infty}:U(\ghat_-)
\to
(\bigotimes_{i=1}^N U(\g_{m_i}) \otimes U(\ol\g_{m_\infty}))
\otimes\C[(u-z_i)^{-1}]_{i=1,\ldots,N}\otimes\C[u]$$

The following result gives an alternative description of
$\ZZ_{(z_i),\infty} ^{(m_i),m_\infty}(\g)$, which shows in particular
that it is independent of $u$.

\begin{prop}\label{pr:modes}
The algebra $\ZZ_{(z_i),\infty}^{(m_i),m_\infty}(\g)$ is equal to the
span of the coefficients of
$\Phi_{u,(z_i)}^{(m_i),m_\infty}(A),A\in\zz(\ghat)$, appearing in
front of the monomials of the form $\prod_{i=1}^N
(u-z_i)^{-n_i}u^{n_\infty}$.
\end{prop}

\proof Since the translation operator acts as $-\partial_t=\partial_u$
on $\ghat(t-u)_-$ and $\Phi_{u,(z_i)}^{(m_i),m_\infty}$ is given by
Laurent expansion at the points $z_1,\ldots,z_N$, it follows that, for
any $A\in\zz(\hatg)$
$$\Phi_{u,(z_i)}^{(m_i),m_\infty}(TA) =
\pa_u\Phi_{u,(z_i)}^{(m_i),m_\infty}(A).$$ where $T$ is given by
\eqref{translation}.  For any $B\in\zz(\ghat)$, the expression
$\Phi_{u,(z_i)}^{(m_i),m_\infty} (B)$, viewed as a rational function
of $u$, is a finite linear combination of monomials of the form
$\prod_{i=1}^N (u-z_i)^{-n_i} u^{n_\infty}$, where $n_i>0,
n_\infty\geq 0$. Such a linear combination is uniquely determined by
its expansion in a power series at a fixed point (or equivalently, by
the values of its derivatives at this point). Therefore the span of
these expressions for a fixed $u$ and all $B \in\zz(\ghat)$ of the
form $T^n A, n\geq 0$, is the same as the span of the coefficients of
$\Phi_{u,(z_i)}^{(m_i),m_\infty}(A)$, considered as a rational
function in $u$, appearing in front of the monomials of the form
$\prod_{i=1}^N (u-z_i)^{-n_i} u^{n_\infty}$. \qed

\smallskip

By \eqref{eq:inv limit} and Proposition \ref{pr:modes},
$\ZZ_{(z_i),\infty} (\g)$ is also independent of $u$. Moreover, if $u$
is regarded as a variable, and $\Phi_{u,(z_i)}$ as an
anti-homomorphism
$$\Phi_{u,(z_i)}:U(\ghat_-) \to \left(U(\g[[t]])^{\otimes N}\otimes
U(t\g[[t]])\right) \wh{\otimes} \left(\C[(u-z_i)^{-1}]_{i=1,\ldots,N}
\otimes \C[u]\right),$$ where the tensor product is completed with
respect to the natural topology on the Lie algebra $\g[[t]]^{\oplus
(N+1)}$, then $\ZZ_{(z_i),\infty}(\g)$ may also be obtained as the
completion of the span of the coefficients of the series $\Phi_
{u,(z_i)}(A), A \in \zz(\ghat)$, appearing in front of the monomials
of the form $\prod_{i=1}^N (u-z_i)^{-n_i} u^{n_\infty}$.

It is instructive to think of the algebras $\ZZ_{(z_i),\infty}^{(m_i),
m_\infty}(\g)$ and their generalizations considered below as quotients
of the universal Gaudin algebra. However, note that while
$\ZZ_{(z_i),\infty}(\g)$ is a complete topological algebra, the
algebras $\ZZ_{(z_i),\infty}^{(m_i),m_\infty}(\g)$ are its discrete
quotients. For this reason in practice it is easier to work with the
latter. In fact, for any collection of smooth finitely generated
$\g[[t]]$-modules $M_1,\ldots,M_N$, $M_\infty$, the action of
$\g[[t]]$ on $M_i$ (resp., $M_\infty$) factors through $\g_{m_i}$
(resp., $\g_{m_\infty}$). Therefore the action of
$\ZZ_{(z_i),\infty}(\g)$ on the tensor product $\bigotimes_{i=1}^N M_i
\otimes M_\infty$, or on the corresponding space of coinvariants
$$H(\M_1,\ldots,\M_N,\M_\infty) \cong (\bigotimes_{i=1}^N M_i\otimes
M_\infty)/\g$$ factors through that of
$\ZZ_{(z_i),\infty}^{(m_i),m_\infty}(\g)$. Since the only
$\g[[t]]$-modules that we consider are smooth finitely generated
modules, we do not lose any generality by working with the quotients
$\ZZ_{(z_i),\infty}^{(m_i),m_\infty}(\g)$ rather than with
the algebra $\ZZ_{(z_i),\infty}(\g)$ itself.

\subsection{Example: the Gaudin model}    \label{gaudin model}

Consider the simplest case when all $m_i$ and $m_\infty$ are equal to
$1$. The corresponding algebra $\ZZ_{(z_i),\infty}^{(1),1}(\g)$ is a
commutative subalgebra of $(U\g^{\otimes N})^{\g}$. The homomorphism
$\Phi_{u,(z_i)}^{(1),1}(S)$ was defined in \cite{FFR} (see also
\cite{F:faro}).

Consider the Segal-Sugawara vector in $\V_0$:
\begin{equation} \label{sugawara}
S = \frac{1}{2} \sum_{a=1}^d J_{a,-1} J^a_{-1} v_0,
\end{equation}
where $d=\dim\g$ and $\{J_a\},\{J^a\}$ are dual bases of $\g$ \wrt a
non-zero invariant inner product. It is easy to show (see,
e.g., \cite{FB}, Sect. 3.4.8) that this vector belongs to
$\zz(\G)$. Let us compute $\Phi_{u,(z_i)}^{(1),1}(S)$. Let
$$
\Delta = \dfrac{1}{2} \sum_a J_a J^a \in U\g
$$
be the Casimir operator. Formulas \eqref{Psi}--\eqref{J u} readily
yield the following:

\begin{lem}[\cite{FFR}, Prop. 1] \label{coincide}
We have
$$
\Phi_{u,(z_i)}^{(1),1}(S) =
\sum_{i=1}^N \frac{\Delta^{(i)}}{(u-z_i)^2}+
\sum_{i=1}^N \frac{\Xi_i}{u-z_i} ,
$$
where the $\Xi_i$'s are the Gaudin Hamiltonians
\begin{equation}    \label{Gaudin ham}
\Xi_i = \sum_{j\neq i} \sum_{a=1}^d \frac{J_a^{(i)}
J^{a(j)}}{z_i-z_j}, \qquad i=1,\ldots,N.
\end{equation}
\end{lem}

By Proposition \ref{pr:modes}, $\ZZ_{(z_i),\infty}^{(1),1}(\g)$
contains each $\Delta^{(i)}$ and the Gaudin Hamiltonians $\Xi_i$,
$i=1,\ldots,N$. Elements of $\ZZ_{(z_i),\infty}^{(1),1}(\g)$ are
generalized Gaudin Hamiltonians which act on the tensor product
$\bigotimes_{i=1}^N M_i$ of any $N$-tuple of $\g$-modules or on the
space of coinvariants
$$
H(\M_1,\ldots,\M_N,\M_\infty) = (\bigotimes_{i=1}^N M_i
\otimes M_\infty)/\g,
$$
where each $M_i$ is regarded as a $\g[[t]]$-module by letting
$t\g[[t]]$ act by $0$. If all $M_i$'s and $M_\infty$ are highest
weight $\g$-modules, then the corresponding induced modules $\M_i$ and
$\M_\infty$ are highest weight $\ghat_{\ka_c}$-modules. Thus, the
choice $m_i = 1, m_\infty = 1$ corresponds to highest weight
modules. The spectrum of the corresponding algebra
$\ZZ_{(z_i),\infty}^{(1),1}(\g)$ is the space of $^L G$-opers with
regular singularities at the points $z_1,\ldots,z_N$ and $\infty$ (see
\cite{F:faro} and \secref{spectrum} below).

However, if we choose some of the $m_i$'s (or $m_\infty$) to be
greater than $1$, then for a general $\g_{m_i}$-module $M_i$ (or
$M_\infty$) the corresponding induced module $\M_i$ (or $\M_\infty$)
will not be a highest weight module. The corresponding algebra
$\ZZ_{(z_i),\infty}^{(m_i),m_\infty}(\g)$ is isomorphic to the algebra
of functions on the space $^L G$-opers with irregular singularities at
the points $z_i$ (and $\infty$) of orders $m_i$ (resp., $\m_\infty$),
see \thmref{ZZ descr}. Therefore we refer to the corresponding
integrable quantum models as {\em Gaudin models with irregular
singularities}.

\subsection{Another example: the algebra ${\mc A}_\chi$}    \label{Achi}

Consider the case when there are two points: $z_1 \in \C$ and
$\infty$, with $m_1 = 1$ and $m_\infty = 2$. This is the simplest
model with an irregular singularity. The corresponding Gaudin algebra
$\ZZ_{z_1,\infty}^{1,2}(\g)$ is a commutative subalgebra of
$(U(\g)\otimes U(\ol\g_2))^{\g}$. It is easy to see that it is
independent of $z_1$, so we will set $z_1=0$.

The Lie algebra
$$\ol\g_2 = t\g[[t]]/t^2\g[[t]] $$ is abelian and isomorphic to $\g$
as a vector space. Any linear functional $\chi:\g\to\C$ on $\g$
therefore defines an algebra homomorphism $U(\ol \g_2)\cong
S\g\to\C$. Let
\begin{equation}\label{eq:A_chi}
{\mc A}_\chi =
\id\otimes\chi\left(\ZZ_{0,\infty}^{1,2}(\g)\right)\subset U\g
\end{equation}
be the image of $\ZZ_{0,\infty}^{1,2}(\g)$ under the homomorphism
$(U(\g) \otimes U(\ol\g_2))^{\g} \to U(\g)$ given by applying $\chi$
to the second factor. It is clear that ${\mc A}_\chi$ is a commutative
subalgebra of the centralizer $U(\g)^{\g_\chi}\subset U(\g)$, where
$\g_\chi\subseteq\g$ is the stabilizer of $\chi$.

As a subalgebra of $U(\g)$, ${\mc A}_\chi$ acts on any $\g$-module
$M$. From the point of view of the coinvariants construction of
\secref{coinvariants}, this action comes about as follows: we consider
the $\ghat_{\ka_c}$-module $\M$ induced from $M$, attached to $0 \in
\pone$, and the non-highest weight $\ghat_{\ka_c}$-module
$\I_{1,\chi}$ introduced in formula \eqref{I1}. Then the Lie algebra
$\g_{\zn}$ is $\g[t,t^{-1}]$ and the space of coinvariants is
$$H(\M,\I_{1,\chi})=(\M \otimes \I_{1,\chi})/\g[t,t^{-1}] \simeq (M
\otimes I_\chi)/\g \simeq M,$$ since
$I_\chi=\ind_{t\g[[t]]}^{\g[[t]]}\IC_\chi$ is isomorphic to $U\g$ as a
$\g$-module.

The action of $\zz(\ghat)$ on $H(\M,\I_{1,\chi}) \simeq M$ then
factors through that of ${\mc A}_\chi$.

We note that it follows from the definition that ${\mc A}_\chi = {\mc
A}_{c\chi}$ for any non-zero $c \in \C$ and that $\on{Ad}_g({\mc
A}_\chi) = {\mc A}_{\on{Ad}_g(\chi)}$ for any $g \in G$.

\subsection{Multi-point generalization}    \label{multi-point}

The algebra ${\mc A}_\chi$ has a natural multi-point generalization.
Namely, let $m_i, i=1,\ldots,N$, and $m_\infty$ be a collection of
positive integers as in \secref{universal}. Let us also fix characters
$\chi_i: t^{m_i}\g[[t]]\to\C$ and $\chi_\infty:t^{m_\infty} \g[[t]]
\to \C$. We will attach to this data a commutative algebra which
simultaneously generalizes both $\ZZ_{(z_i),\infty}^{(m_i),m_\infty}$
and ${\mc A}_\chi$. Note that a character
$$\chi: t^m\g[[t]]\to\C$$ has to vanish on the derived subalgebra of
$t^m\g[[t]]$, that is on $t^{2m}\g[[t]]$. Hence, defining
$\chi$ is equivalent to choosing an arbitrary linear functional on
the abelian Lie algebra
$$t^m\g[[t]]/t^{2m}\g[[t]]\simeq\g\otimes\C^m.$$
Let
$$
I_{m,\chi} = \Ind_{t^m\g[[t]]}^{\g[[t]]} \C_\chi, \qquad
\ol{I}_{m,\chi} = \Ind_{t^m\g[[t]]}^{t\g[[t]]} \C_\chi.
$$
Note that the action of $\g[[t]]$ on $I_{m,\chi}$ factors through
$\g_{2m}=\g[[t]]/t^{2m}\g[[t]]$, and as a $\g_{2m}$-module,
$I_{m,\chi}$ is isomorphic to
$$I_{m,\chi} \simeq \Ind_{t^m\g_{2m}}^{\g_{2m}}\C_\chi$$ (and
similarly for $\ol{I}_{m,\chi}$). Denote by
${\mc I}_{m,\chi}$ (resp., $\ol{\mc I}_{m,\chi}$) the annihilator of
$I_{m,\chi}$ in $U(\g_{2m})$ (resp., of $\ol{I}_{m,\chi}$ in
$U(\ol\g_{2m})$) and by $U_{m,\chi}$ (resp., $\ol{U}_{m,\chi}$) the
quotient of $U(\g_{2m})$ (resp., $U(\ol\g_{2m})$) by ${\mc
I}_{m,\chi}$ (resp., $\ol{\mc I}_{m,\chi}$). Thus, $U_{m,\chi}$ is the
image of $U(\g_{2m})$ and $U(\g[[t]])$ in $\on{End}_{\C}I_{m,\chi}$,
and $\ol{U}_{m,\chi}$ is the image of $U(\ol\g_{2m})$ and $U(t\g[[t]])$
in $\on{End}_{\C}I_{m,\chi}$.

Now, given the data of $(m_i), m_\infty,(\chi_i),\chi_\infty$, we
obtain the algebra $\bigotimes_{i=1}^N U_{m_i,\chi_i}
\otimes\ol{U}_{m_\infty,\chi_\infty}$, which is isomorphic to the
quotient of $\bigotimes_{i=1}^N U(\g_{2m_i})\otimes
U(\ol\g_{2m_\infty})$ by a two-sided ideal.  Let
$$\A_{(z_i),\infty}^{(m_i),m_\infty}(\g)_{(\chi_i),\chi_\infty}
\subset
\left( U_{m_i,\chi_i} \otimes \ol{U}_{m_\infty,\chi_\infty}
\right)^{\g}$$
be the image of
$$
\ZZ_{(z_i),\infty}^{(2m_i),2m_\infty}(\g) \subset \left( U_{2m_i}
\otimes \ol{U}_{2m_\infty} \right)^{\g}.
$$
Equivalently,
$\A_{(z_i),\infty}^{(m_i),m_\infty}(\g)_{(\chi_i),\chi_\infty}$ is
the image of the universal Gaudin algebra $\ZZ_{(z_i),\infty}(\g)$
under the homomorphism
\begin{equation}\label{surjection1}
U(\g[[t]])^{\otimes N} \otimes U(t\g[[t]]) \to \bigotimes_{i=1}^N
U_{m_i,\chi_i} \otimes \ol{U}_{m_\infty,\chi_\infty}.
\end{equation}

Note in particular that
$${\mc A}_\chi = \A_{0,\infty}^{1,1}(\g)_{0,\chi}.$$

The algebra
$\A_{(z_i),\infty}^{(m_i),m_\infty}(\g)_{(\chi_i),\chi_\infty}$
naturally acts on the tensor product $\bigotimes_{i=1}^N
I_{m_i,\chi_i} \otimes I_{m_\infty,\chi_\infty}$, or on the
corresponding space of coinvariants
$$
H(\I_{m_1,\chi_1},\ldots,\I_{m_N,\chi_N},\I_{m_\infty,\chi_\infty})
\cong
(\bigotimes_{i=1}^N I_{m_i,\chi_i} \otimes
I_{m_\infty,\chi_\infty})/\g \simeq (\bigotimes_{i=1}^N I_{m_i,\chi_i}
\otimes \ol{I}_{m_\infty,\chi_\infty}).
$$

For $\chi=0$, $U_{m,0}\simeq U(\g_m)$ and $\ol{U}_{m,0} \simeq
U(\ol\g_m)$. But for a non-zero character $\chi$ the algebras
$U_{m,\chi}$ and $\ol{U}_{m,\chi}$ are not in general isomorphic to
universal enveloping algebras.

There is however one exception. If $\chi_\infty: t^m\g[[t]] \to
\C$ factors as
$$
t^m\g[[t]] \to \g \otimes t^m \overset{\chi}\longrightarrow
\C
$$
for some linear functional $\chi$ on $\g$, then $\ol{U}_{m,\chi}
\simeq \ol{U}(\g_{m-1})$ for any $\chi$. Therefore the algebra
$\A_{(z_i),\infty}^{(m_i),m_\infty}(\g)_{(0),\chi_\infty}$ may
be realized in $\bigotimes_{i=1}^N U_{m_i} \otimes
\ol{U}_{m_\infty-1}$ in this case.

In particular, if $m=1$ then $\ol{U}_{1,\chi}\simeq \C$ regardless of
$\chi$. Therefore in the case when $\chi_i = 0$ for $i=1,\ldots,N$,
and $m_\infty=1$ the corresponding algebra
$\A_{(z_i),\infty}^{(m_i),1}(\g)_{(0),\chi_\infty}$ is a subalgebra
of the universal enveloping algebra $\bigotimes_{i=1}^N U(\g_{m_i})$.

In the case when $m_i=1, i=1,\ldots,N$, we obtain a subalgebra of
$U(\g)^{\otimes m}$ that has been previously constructed in
\cite{Ryb2}. This subalgebra
$\A_{(z_i),\infty}^{(1),1}(\g)_{(0),\chi}$ is obtained by applying
the homomorphism
$$
\id \otimes \chi: U(\g)^{\otimes m} \otimes S(\g) \to
U(\g)^{\otimes m}
$$
to the algebra $\ZZ_{(z_i),\infty}^{(1),2}(\g)$. Let us apply $\id
\otimes \chi$ to $\Phi_{u,(z_i)}^{(1),2}(S)$, where $S$ is the
Segal-Sugawara element in $\V_0$. Then we find, in the same way as in 
\lemref{coincide}, that
$$
(\id \otimes \chi) \circ \Phi_{u,(z_i)}^{(1),2}(S) =
\sum_{i=1}^N \frac{\Delta^{(i)}}{(u-z_i)^2}+
\sum_{i=1}^N \frac{\Xi_{i,\chi}}{u-z_i} + (\chi,\chi),
$$
where
\begin{equation}    \label{Gaudin ham1}
\Xi_{i,\chi} = \sum_{j\neq i} \sum_{a=1}^d \frac{J_a^{(i)}
J^{a(j)}}{z_i-z_j} + \chi^{(i)}, \qquad
  i=1,\ldots,N,
\end{equation}
where we identity $\h \simeq \h^*$ using the invariant inner product
used in the definition of $S$. Thus, $\Xi_{i,\chi}, i=1,\ldots,N$, are
elements of $\A_{(z_i),\infty}^{(1),1}(\g)_{(0),\chi}$. These
operators appeared in \cite{FMTV} in the study of generalized KZ
equations. We note that in the case of $\g=\sw_2$ they were probably
first considered in \cite{Skl}.

\section{Associated graded algebras and Hitchin systems}
\label{graded}

In this section we consider generalized Hitchin systems on the moduli
spaces of Higgs bundles, where the Higgs fields are allowed to have
poles. Integrable systems of this type, generalizing the original
Hitchin systems from \cite{Hitch} (which correspond to Higgs fields
without poles), have been previously considered in
\cite{Beauville,Bot,Markman,DM,ER}. Here we will focus on the case
when the underlying curve is $\pone$. We will construct algebras of
Poisson commuting Hamiltonians in these systems in the standard
way. We will then show that the generalized Gaudin algebras of
commuting quantum Hamiltonians constructed in the previous section are
quantizations of the Poisson commutative algebras of Hitchin
Hamiltonians. The proof of this result relies on a local statement,
due to \cite{FF:gd,F:wak}, that the center $\zz(\ghat)$ of
$\V_{0,\ka_c}$ is a quantization of
$S(\g\ppart/\g[[t]])^{\g[[t]]}$. By definition, the Gaudin algebras
are quotients of $\zz(\ghat)$, whereas the algebras of Hitchin
Hamiltonians are quotients of $S(\g\ppart/\g[[t]])^{\g[[t]]}$, hence
the result.

The idea of using the center at the critical level for quantizing the
Hitchin systems is due to Beilinson and Drinfeld \cite{BD}, who used
it to quantize the Hitchin systems defined on arbitrary smooth
projective curves, without ramification. Here we develop this theory
in a ``transversal'' direction: quantizing the Hitchin systems on
$\pone$, but with arbitrary ramification. The two scenarios may
certainly be combined, giving rise to quantum integrable systems
corresponding to arbitrary curves, with ramification. However, these
systems are hard to analyze explicitly unless the underlying curve is
$\pone$ or an elliptic curve. In the general setting the focus shifts
instead to the investigation of the salient features of the ${\mc
D}$-modules on the moduli stacks of bundles with level (or parabolic)
structures defined by the corresponding algebras of quantum
Hamiltonians (see \secref{connection to glc}). It follows from the
results of \cite{BD} that these ${\mc D}$-modules are Hecke
eigensheaves, and hence they play an important role in the geometric
Langlands correspondence.

On the other hand, in the case of $\pone$ the generalized Gaudin
algebras may be analyzed explicitly. In the case of regular
singularities, corresponding to the ordinary Gaudin models, this has
been done in \cite{FFR,F:faro} (see also \cite{ER,EFR} for a
generalization to the case of elliptic curves). Here we look closely
at another special case corresponding to irregular singularity of
order $2$ at one point of $\pone$. The corresponding Poisson
commutative algebra of Hitchin Hamiltonians may be identified with the
shift of argument algebra $\ol{\A}_\chi$ introduced in \cite{MF}. We
show that the quantum commutative algebra ${\mc A}_\chi$ introduced in
\secref{Achi} is the quantization of $\ol{\A}_\chi$ for all regular
$\chi \in \g^*$. This has been previously proved in \cite{Ryb2} in the
case when $\chi$ is regular semi-simple. In addition, we show that
when $\chi$ is regular semi-simple the algebra ${\mc A}_\chi$
contains the DMT Hamiltonians \eqref{T ham}.

\subsection{Hitchin systems with singularities}

Let us first recall the definition of the unramified Hitchin
system. Denote by $G$ the connected simply-connected simple Lie group
with Lie algebra $\g$. Let $\MM_G$ be the moduli space of stable
$G$-bundles on a smooth projective curve $X$. The tangent space to
$\MM_G$ at ${\mc F} \in \MM_G$ is isomorphic to $H^1(X,\g_{{\mc F}})$,
where $\g_{{\mc F}} = {\mc F} \underset{G}\times \g$. Hence, by the
Serre duality, the cotangent space at ${\mc F}$ is isomorphic to
$H^0(X,\g^*_{{\mc F}} \otimes \Omega)$, where $\Omega$ is the
canonical line bundle on $X$, by Serre duality. A vector $\eta \in
H^0(X,\g^*_{{\mc F}} \otimes \Omega)$ is referred to as a {\em Higgs
field}. We construct the Hitchin map $p: T^* \MM_G \to {\mc H}_G$,
where ${\mc H}_G$ is the {\em Hitchin space}
\begin{equation}    \label{Hitchin space}
{\mc H}_G(X) = \bigoplus_{i=1}^\ell H^0(X,\Omega^{\otimes(d_i+1)}),
\end{equation}
as follows. We use the following result of C. Chevalley.

\begin{thm}    \label{chevalley}
The algebra $\on{Inv} \g^*$ of $G$-invariant polynomial functions on
$\g^*$ is isomorphic to the graded polynomial algebra
$\C[\ol{P}_1,\ldots,\ol{P}_\ell]$, where $\deg \ol{P}_i = d_i+1$, and
$d_1,\ldots,d_\ell$ are the exponents of $\g$.
\end{thm}

For $\eta \in H^0(X,\g^*_{{\mc F}} \otimes \Omega)$, $\ol{P}_i(\eta)$
is well-defined as an element of $H^0(X,\Omega^{\otimes(d_i+1)})$.  By
definition, the Hitchin map $p$ takes $({\mc F},\eta) \in T^* \MM_G$
to
$$
(\ol{P}_1(\eta),\ldots,\ol{P}_\ell(\eta)) \in {\mc H}_G.
$$

\begin{remark}
This definition depends on the choice of generators of $\on{Inv}
\g^*$, which is not canonical. To give a more canonical definition,
let
$$
{\mc P} = \on{Spec} \on{Inv} \g^*.
$$
Since $\on{Inv} \g^*$ is a graded algebra, we obtain a canonical
$\C^\times$-action on ${\mc P}$. Now let $\Omega^\times$ be the
$\C^\times$-bundle on $X$ corresponding to the line bundle
$\Omega$. We then have a vector bundle
$$
{\mc P}_\Omega = \Omega^\times \underset{\C^\times}\times {\mc P},
$$
and we set
$$
{\mc H}_G = H^0(X,{\mc P}_\Omega).
$$
A choice of homogeneous generators $\ol{P}_i, i=1,\ldots,\ell$, of
$\on{Inv} \g^*$ gives rise to a set of coordinates on ${\mc P}$ which
enable us to identify ${\mc H}_G$ with \eqref{Hitchin space}. However,
using this definition of ${\mc H}_G$, we obtain a definition of the
Hitchin map that is independent of the choice of generators. Likewise,
the generalized Hitchin map considered below may also be defined in a
generator-independent way. But in order to simplify the exposition we
will define them by using a particular choice of generators $\ol{P}_i,
i=1,\ldots,\ell$, of $\on{Inv} \g^*$.\qed
\end{remark}

Now, given a linear functional $\phi: {\mc H}_G \to \C$, we obtain a
function $\phi \circ p$ on $T^* \MM_G$. Hitchin \cite{Hitch} has
shown that for different $\phi$'s these functions Poisson commute with
respect to the natural symplectic structure on $T^* \MM_G$, and
together they define an algebraically completely integrable
system.

Let us generalize the Hitchin systems to the case of Higgs fields with
singularities. Let $x_i, i=1,\ldots,n$ be a collection of distinct
points on $X$. Let us choose a collection of positive integers $m_i,
i=1,\ldots,n$. Denote by $\MM_{G,(x_i)}^{(m_i)}$ the moduli space of
semi-stable $G$-bundles on $X$ with level structures of order $m_i$ at
$x_i, i=1,\ldots,n$. Recall that a level structure on a $G$-bundle
${\mc F}$ on a curve $X$ of order $m$ at a point $x \in X$ is a
trivialization of ${\mc F}$ on the $(m-1)$st infinitesimal
neighborhood of $x$. Then a cotangent vector to a point of
$\MM_{G,(x_i)}^{(m_i)}$ is a Higgs field with poles of orders at most
$m_i$ at the points $x_i$,
$$
\eta \in H^0(X,\g^*_{\mc F} \otimes \Omega(m_1 x_1 + \ldots + m_n
x_n)).
$$
In the same way as above one constructs a generalized Hitchin map
$$
p: T^* \MM_{G,(x_i)}^{(m_i)} \to {\mc H}_{G,(x_i)}^{(m_i)} =
\bigoplus_{j=1}^\ell H^0(X,\Omega(m_1 x_1 + \ldots + m_n
x_n)^{\otimes(d_j+1)}),
$$
and shows that it defines an algebraically completely integrable
system. These systems have been previously studied in
\cite{Beauville,Bot,Markman,DM,ER}.

We are interested in these integrable systems in the case when $X =
\pone$, with the marked points $z_1,\ldots,z_N,\infty$ (with respect
to a global coordinate $t$, as before). In this case a $G$-bundle is
semi-stable if and only if it is trivial. Thus,
$\MM_{G,(x_i)}^{(m_i)}$ may be identified in this case with the
quotient
$$
\left. \prod_{i=1}^N G_{m_i} \times G_{m_\infty} \right/ G_{\on{diag}}
\simeq \prod_{i=1,\ldots,N} G_{m_i} \times \ol{G}_{m_\infty}.
$$
Here $G_m, \ol{G}_m$ are the Lie groups corresponding to the Lie
algebras $\g_m = \g[[t]]/t^m\g[[t]], \ol\g_m =
t\g[[t]]/t^m\g[[t]]$. Hence $T^*\MM_{G,(x_i)}^{(m_i)}$ is
isomorphic to the Hamiltonian reduction of
$$
T^*(\prod_{i=1}^N G_{m_i} \times G_{m_\infty})
$$
with respect to the diagonal action of $G$ (and the $0$ orbit in
$\g^*$). We can an will identify it with
\begin{equation}    \label{open substack}
\prod_{i=1,\ldots,N} G_{m_i} \times \ol{G}_{m_\infty} \times
\prod_{i=1,\ldots,N} (\g_{m_i})^* \times (\ol\g_{m_\infty})^*,
\end{equation}
A point of \eqref{open substack} is then a collection
$(g_\al,A_\al(t)dt)_{\al = 1,\ldots,N,\infty}$, where $g_i \in
G_{m_i}, i=1,\ldots,N; g_\infty \in \ol{G}_{m_\infty}$,
$$
A_i(t) dt= \sum_{k=-m_i}^{-1} A_{i,k} t^k dt, \quad A_{i,k} \in \g^*,
\qquad A_\infty(t) dt = \sum_{k=-m_\infty}^{-2} A_{\infty,k} t^k dt,
\quad A_{\infty,k} \in \g^*.
$$
We associate to it a
$\g^*$-valued one-form on $\pone$ with singularities at the marked
points (also known as an ``$L$-operator''):
\begin{equation}    \label{eta}
\eta = \sum_{i = 1,\ldots,N} \on{Ad}_{g_i}(A_i(t-z_i)) d(t-z_i) +
\on{Ad}_{g_\infty}(A_{\infty}(t^{-1})) d(t^{-1}),
\end{equation}
whose polar parts are given by the one-forms $A_i(t_\al) dt_\al$,
conjugated by $g_\al$'s. The Hitchin map $p$ takes this one-form to
\begin{multline}    \label{Hitchin map}
\eta \mapsto (\ol{P}_1(\eta),\ldots,\ol{P}_\ell(\eta)) \in {\mc
H}_{G,(z_i),\infty}^{(m_i),m_\infty} = \\ \bigoplus_{j=1}^\ell
H^0(\pone,\Omega(m_1 z_1 + \ldots + m_N z_N + m_\infty
\infty)^{\otimes(d_j+1)}).
\end{multline}

Now any linear functional $\phi: {\mc
H}_{G,(z_i),\infty}^{(m_i),m_\infty} \to \C$ gives rise to a function
$H_\phi = p \circ \phi$, and, according to the above general result of
\cite{Hitch,Beauville,Bot,Markman,DM,ER} these functions Poisson
commute with each other.

By definition, the symplectic manifold \eqref{open substack} fibers
over $(\prod_\al G_{m_\al})/G$. The fiber over the identify coset is
the Poisson submanifold
\begin{equation}    \label{fiber}
\prod_{i=1,\ldots,N} (\g_{m_i})^* \times (\ol\g_{m_\infty})^*
\end{equation}
with its natural Kirillov-Kostant structure. Let $\wt{p}$ be the
restriction of the Hitchin map to this subspace. Then we obtain a
system of Poisson commuting Hamiltonians $\wt{p} \circ \phi$ on
\eqref{fiber} for $\phi \in ({\mc
  H}_{G,(z_i),\infty}^{(m_i),m_\infty})^*$.

The Poisson algebra of polynomial functions on \eqref{fiber} is
isomorphic to
\begin{equation}    \label{assoc gr}
\bigotimes_{i=1,\ldots,N} S(\g_{m_i}) \otimes S(\ol\g_{m_\infty}).
\end{equation}
Therefore $\wt{p}$ gives rise to a homomorphism of commutative
algebras\footnote{Recall that for an affine algebraic variety $Z$ we
denote by $\on{Fun} Z$ the algebra of regular functions on $Z$.}
\begin{equation}    \label{olPsi}
\ol\Psi_{(z_i),\infty}^{(m_i),m_\infty}: \on{Fun} {\mc
H}_{G,(z_i),\infty}^{(m_i),m_\infty} \to \bigotimes_{i=1,\ldots,N}
S(\g_{m_i}) \otimes S(\ol\g_{m_\infty})
\end{equation}
(actually, it is easy to see that the image belongs to the subalgebra
of $G$-invariants).

\begin{lem}    \label{is inj}
The homomorphism $\ol\Psi_{(z_i),\infty}^{(m_i),m_\infty}$ is injective.
\end{lem}

\proof It is known that (in the general case) the Hitchin map is
surjective, see \cite{Hitch,Fal,Bot,Markman,DM}. This implies the
statement of the lemma.\qed

Let $\ol\ZZ_{(z_i),\infty}^{(m_i),m_\infty}(\g)$ be the image of
$\ol\Psi_{(z_i),\infty}^{(m_i),m_\infty}$. This is a Poisson commutative
subalgebra of the Poisson algebra \eqref{assoc gr}.

Now we can explain the connection between the ramified Hitchin systems
and the generalized Gaudin algebras
$\ZZ_{(z_i),\infty}^{(m_i),m_\infty}(\g)$ introduced in
\secref{universal}.  Note that the algebra
$\ZZ_{(z_i),\infty}^{(m_i),m_\infty}(\g)$ is defined as a quotient of a
homomorphism of universal enveloping algebra and therefore inherits a
filtration that is compatible with the PBW filtration on
$\bigotimes_{i=1,\ldots,N} U(\g_{m_i}) \otimes
U(\ol\g_{m_\infty})$. The associated graded algebra to the latter is
precisely $\bigotimes_{i=1,\ldots,N} S(\g_{m_i}) \otimes
S(\ol\g_{m_\infty})$.

\begin{thm}    \label{quantized}
The associated graded of the commutative subalgebra
$$
\ZZ_{(z_i),\infty}^{(m_i),m_\infty}(\g) \subset
\bigotimes_{i=1,\ldots,N} U(\g_{m_i}) \otimes U(\ol\g_{m_\infty})
$$
is the Poisson commutative subalgebra
$$
\ol\ZZ_{(z_i),\infty}^{(m_i),m_\infty}(\g) \subset
\bigotimes_{i=1,\ldots,N} S(\g_{m_i}) \otimes S(\ol\g_{m_\infty}).
$$
\end{thm}

Thus, we obtain that $\ZZ_{(z_i),\infty}^{(m_i),m_\infty}(\g)$ is a
quantization of the Poisson commutative algebra
$\ol\ZZ_{(z_i),\infty}^{(m_i),m_\infty}(\g)$ of Hitchin's Hamiltonians.

The proof of this result given below in \secref{back} follows from the
fact that the algebra $\zz(\ghat)$, whose quotient is
$\ZZ_{(z_i),\infty}^{(m_i),m_\infty}(\g)$, is the quantization of the
algebra of invariant functions on $\g^*[[t]]$. It is this local result
that is responsible for the global quantization results such as
\thmref{quantized} or \thmref{th:quantisation} below (and the results
of Beilinson--Drinfeld \cite{BD} for general curves in the unramified
case). We will explain this local result in the next section.

\subsection{The associated graded algebra of $\zz(\hatg)$}
\label{grad}

The PBW filtration on $U\hatg_\ka$ induces one on $\IV_ {0,\ka}$, with
associated graded
$$\gr(\IV_{0,\ka})=S(\hatg_\ka\ppart/\hatg_+)=S(\g\ppart/\g[[t]]).$$
The action of $\g[[t]]$ on $\IV_{0,\ka}$ is readily seen to preserve
this filtration and therefore descends to one on $S(\g\ppart/\g[[t]])$.
This latter action is independent of the level $\ka$ and is given by
derivations induced by the adjoint action of $\g[[t]]$ on $\g\ppart/
\g[[t]]$.

We therefore obtain an embedding
$$\gr(\IV_{0,\ka}^{\g[[t]]})\subset S(\g\ppart/\g[[t]])^{\g[[t]]}$$
which, for $\ka=\ka_c$ is in fact an equality, according to Theorem
\ref{thm:FF} below. We begin by describing the right-hand side.

Split $\g\ppart$ (considered as a vector space) as
\begin{equation}\label{eq:decomposition}
\g\ppart=\g[[t]] \oplus \hatg_{-},
\end{equation}
where $\hatg_{-} = t^{-1}\g[t^{-1}]$, and identify
$S(\g\ppart/\g[[t]])$ with $S(\hatg_-)$. Under this
identification, the adjoint action of $X\in\g[[t]]$ on $Y\in \ghat_-
\cong\g\ppart/\g[[t]]$ is given by
$$\ad(X)_- Y=[X,Y]_-,$$
where $Z_-$ is the component of $Z\in\g\ppart$ along $\hatg_-$.

Let $T=-\pa_t\in\DerO$ be the translation operator acting on
$\g\ppart$, so that $T X_n=-n X_{n-1}$ where $X_n=X\otimes t^n \in
\hatg_\ka$, $X\in\g$, $n\in\Z$. Then $T$ preserves the decomposition
\eqref{eq:decomposition} and extends to a derivation of $S(\hatg_-)$.

\begin{lem}\label{le:T ad}
$$[T,\ad(X)_-]=\ad(TX)_-$$
\end{lem}
\proof Both sides are derivations and they coincide on $\hatg_-$. \qed

\smallskip

For any $X\in\g$, set
$$\ol{X}(z)=\sum_{n<0}\ol{X}_n z^{-n-1}\in (\hatg_-)[[z]]\subset
S(\hatg_-)[[z]],$$ where $\ol{X}_n, n<0$, denotes $X \otimes t^n$
considered as an element of $S(\hatg_-)$. We extend the assignment
$X\rightarrow \ol{X}(z)$ to an algebra homomorphism
$$S\g\longrightarrow S(\hatg_-)[[z]]$$
and denote the image of $P\in S\g$ by $P(z)=\sum_{n<0} P_n
z^{-n-1}$. The following summarizes the main properties of the
map $P\mapsto P(z)$.

\begin{lem}\label{pr:P(z)}
The following holds for any $P\in S\g$, $X\in\g$ and $k\in\IN$,
\begin{enumerate} 
\item $\displaystyle{TP(z)=\frac{dP(z)}{dz}}$
\item $\ad(X_k)_- P(z)=z^k(\ad(X)P)(z)$
\end{enumerate}
\end{lem}
\proof Both sides of the identities are derivations which are readily
seen to coincide on the elements $Y(z)$, $Y\in\g$. \qed

\smallskip

Note that property (i) above is equivalent to the fact that, for any
$P \in S \g$ and $n>0$,
\begin{equation}\label{eq:recurrence}
P_{-n}=\frac{T^{(n-1)}}{(n-1)!}P_{-1}
\end{equation}
while (ii) implies that the assignment $P\rightarrow P(z)$ maps
$(S\g)^\g$ to $S(\hatg_-)^{\g[[t]]}[[z]]$.

Recall from \lemref{chevalley} that $(S\g)^\g$ is a polynomial algebra
in $\ell=\rk(\g)$ generators and that one may choose a system of
homogeneous generators $\ol{P}_1,\ldots,\ol{P}_\ell$ such that
$\deg\ol{P}_i=d_i+1$, where the $d_i$'s are the exponents of $\g$.

The following theorem is due to A. Beilinson and V. Drinfeld \cite{BD}
(see \cite{F:wak}, Prop. 9.3, for an exposition).

\begin{thm} \label{thm:BD}
The algebra $S(\g\ppart/\g[[t]])^{\g[[t]]}$ is freely generated by the
elements $\ol{P}_{i,n}$, $i=1,\ldots,\ell$, $n<0$.
\end{thm}

The following result, due to \cite{FF:gd} (see \cite{F:wak}, Theorem
9.6), enables us to quantize the Hitchin Hamiltonians.

\begin{thm} \label{thm:FF}
The inclusion
$$
\on{gr} \zz(\ghat) = \gr(\IV_{0,\ka_c}^\g[[t]]) \hookrightarrow
S(\g\ppart/\g[[t]])^{\g[[t]]}
$$
is an isomorphism.
\end{thm}

In other words, all $\g[[t]]$-invariants in $S(\g\ppart/\g[[t]])$ may
be quantized.

In particular, each generator $\ol{P}_{i,-1}$ of
$S(\g\ppart/\g[[t]])^{\g[[t]]}$ may be lifted to an element $S_i \in
\zz(\ghat)$ whose symbol in $S(\g\ppart/\g[[t]])^{\g[[t]]}$ is equal
to $\ol{P}_{i,-1}$. The element
$$
\frac{T^{-n-1}}{(-n-1)!} S_i \in \zz(\ghat)
$$
then gives us a lifting for $\ol{P}_{i,n}, n<-1$. Explicit formulas
for these elements are unknown in general (however, recently some
elegant formulas have been given in \cite{CT} in the case when
$\g=\sw_n$). But for our purposes we do not need explicit formulas for
the $S_i$'s. All necessary information about their structure is in
fact contained in the following lemma.

Let us observe that both $\V_{0,\ka_c}$ and $S(\g\ppart/\g[[t]])$ are
$\Z_+$-graded with the degrees assigned by the formula $\deg A_n =
-n$, and this induces compatible $\Z_+$-gradings on $\zz(\ghat)$ and
$S(\g\ppart/\g[[t]])^{\g[[t]]}$. This grading is induced by the vector
field $L_0 = - t\pa_t$. It follows from the construction that $\deg
\ol{P}_{i,-1} = d_i+1$. Therefore we may, and will, choose $S_i \in
\V_{0,\ka_c}$ to be homogeneous of the same degree.

In $\V_0 = \V_{0,\ka_c}$ we have a basis of lexicographically ordered
monomials of the form
$$
J^{a_1}_{-n_1} \ldots J^{a_m}_{-n_m} v_0, \qquad n_1 \geq
\ldots \geq n_m > 0,
$$
where $\{ J^a \}$ is a basis of $\g$. The element $\ol{P}_{i,-1}$ of
$S(\g\ppart/\g[[t]])^{\g[[t]]}$ is a linear combination of monomials
in $\ol{J}^a_{-1}$ of degree $d_i+1$. Let $P_i$ be the element of
$\V_0$ obtained by replacing each of these monomials in the
$\ol{J}^a_{-1}$'s by the corresponding lexicographically ordered
monomial in the $J^a_{-1}$'s. Then $P_i$ is the leading term of the
sought-after element $S_i \in \zz(\ghat)$ (for $i=1$ we actually have
$S_1 = P_1$, but for $i>1$ there are lower order terms). By taking
into account the requirement that $\deg S_i = d_i+1$, we obtain the
following useful

\begin{lem}    \label{lexico}
The element $S_i \in \zz(\ghat)$ is equal to $P_i$ plus the sum of
lexicographically ordered monomials of orders less than $d_i+1$. Each
of these lower order terms contains at least one factor $J^a_n$ with
$n<-1$.
\end{lem}

\subsection{Back to the Hitchin systems}    \label{back}

Now we are ready to prove \thmref{quantized}.

\smallskip

\noindent{\em Proof of \thmref{quantized}.}
Observe that the homomorphism $\Phi_{u,(z_i)}^{(m_i),m_\infty}$
given by formula \eqref{Phi} preserves filtrations. Therefore it
gives rise to a homomorphism of associated graded algebras
$$
\ol\Phi_{u,(z_i)}^{(m_i),m_\infty}: S(\hatg_-) \to
\bigotimes_{i=1,\ldots,N} S(\g_{m_i}) \otimes S(\ol\g_{m_\infty}).
$$
Let us consider the image of $S(\hatg_-)^{\g[[t]]}$ under this map. As
in the quantum case, it may also be defined as the span of the
coefficients of the series $\ol\Phi_{u,(z_i)}^{(m_i),m_\infty}(A), A
\in S(\ghat_-)^{\g[[t]]}$, appearing in front of the monomials of the
form $\prod_{i=1}^N (u-z_i)^{-n_i} u^{n_\infty}$.

It follows from the description of $S(\ghat_-)^{\g[[t]]}$ given in
\thmref{thm:BD} and the construction of the Hitchin map \eqref{Hitchin
map} that the image of $S(\hatg_-)^{\g[[t]]}$ under
$\ol\Phi_{u,(z_i)}^{(m_i),m_\infty}$ is precisely the subalgebra
$$
\ol\ZZ_{(z_i),\infty}^{(m_i),m_\infty}(\g) \subset
\bigotimes_{i=1,\ldots,N} S(\g_{m_i}) \otimes S(\ol\g_{m_\infty}).
$$
On the other hand, we know that
$\ol\ZZ_{(z_i),\infty}^{(m_i),m_\infty}(\g)$ is the free polynomial
algebra of regular functions on the graded vector space ${\mc
H}_{G,(z_i),\infty}^{(m_i),m_\infty}$. From its description in formula
\eqref{Hitchin map} we obtain that we may choose as homogeneous
generators of this polynomial algebra the images under
$\ol\Phi_{u,(z_i)}^{(m_i),m_\infty}$ of the following generators of
$S(\ghat_-)^{\g[[t]]}$:
$$
\ol{P}_{i,-n_i-1}, \qquad
i=1,\ldots,\ell; \; 0 \leq n_i \leq (d_i+1)(\sum_j m_j + m_\infty -2)
$$
(note that the number on the right is the degree of the line bundle
$\Omega(m_1 z_1 + \ldots + m_N z_N + m_\infty
\infty)^{\otimes(d_j+1)}$). Each of them lifts to a generator of
$\zz(\ghat)$ (same notation, but without a bar). It is clear from
the definition of $\Phi_{u,(z_i)}^{(m_i),m_\infty}$ that the images of
these generators of $\zz(\ghat)$ in
$\ZZ_{(z_i),\infty}^{(m_i),m_\infty}(\g)$ under the homomorphism
$\Phi_{u,(z_i)}^{(m_i),m_\infty}$ generate
$\ZZ_{(z_i),\infty}^{(m_i),m_\infty}(\g)$.

Thus, we obtain two sets of generators of
$\ZZ_{(z_i),\infty}^{(m_i),m_\infty}(\g)$ and
$\ol\ZZ_{(z_i),\infty}^{(m_i),m_\infty}(\g)$, such that the latter are
the symbols of the former. In addition, the latter are algebraically
independent. Therefore we find that the symbol of any homogeneous
polynomial in the generators of
$\ZZ_{(z_i),\infty}^{(m_i),m_\infty}(\g)$ is equal to the corresponding
polynomial in their symbols. Thus, we obtain the desired isomorphism
$$
\on{gr} \ZZ_{(z_i),\infty}^{(m_i),m_\infty}(\g) \simeq
\ol\ZZ_{(z_i),\infty}^{(m_i),m_\infty}(\g).
$$
Furthermore, it fits into a commutative diagram
$$
\begin{CD}
\on{gr} \zz(\ghat) @>\Phi_{u,(z_i)}^{(m_i),m_\infty}>> \on{gr}
\ZZ_{(z_i),\infty}^{(m_i),m_\infty}(\g)
\\ @VVV @VVV \\ S(\ghat_-)^{\g[[t]]}
@>{\ol\Phi_{u,(z_i)}^{(m_i),m_\infty}}>>
\ol\ZZ_{(z_i),\infty}^{(m_i),m_\infty}(\g)
\end{CD}
$$
where the vertical maps are isomorphisms.\qed

\medskip

\begin{remark}
The problem we had to deal with in the above proof is that {\em a
priori} we do not have a well-defined map
\begin{equation}    \label{well-defined map}
\on{gr} \ZZ_{(z_i),\infty}^{(m_i),m_\infty}(\g) \to
\ol\ZZ_{(z_i),\infty}^{(m_i),m_\infty}(\g).
\end{equation}
However, such a map may be constructed following \cite{BD}.

We have discussed above the Hitchin map on the moduli space of Higgs
bundles on $\pone$ corresponding to the {\em trivial} $G$-bundle (as
this is the only $G$-bundle on $\pone$ that is semi-stable). This
moduli space is in fact an open substack in the moduli stack
$\Bun_{G,(z_i),\infty}^{(m_i),m_\infty}$ of all Higgs bundles
(corresponding to arbitrary $G$-bundles). The Hitchin map extends to a
proper morphism from the entire stack of Higgs bundles to ${\mc H}_G$
(see, e.g., \cite{BD}). This implies that the algebra of global
functions on the stack of Higgs bundles is equal the algebra of
polynomial functions on the Hitchin space ${\mc H}_{G,(x_i)}^{(m_i)}$;
that is, the algebra $\ol\ZZ_{(z_i),\infty}^{(m_i),m_\infty}(\g)$. On
the other hand, as in \cite{BD}, the algebra
$\ZZ_{(z_i),\infty}^{(m_i),m_\infty}(\g)$ may be identified with a
subalgebra (and, {\em a posteriori}, the entire algebra) of the
algebra of (critically twisted) global differential operators on the
stack $\Bun_{G,(z_i),\infty}^{(m_i),m_\infty}$ of $G$-bundles on
$\pone$ with the level structures. Now we obtain a well-defined
injective map \eqref{well-defined map}: it corresponds to taking the
symbol of a differential operator. (Recall that the symbol of a global
differential operator on a variety, or an algebraic stack, $M$ is a
function on the cotangent bundle to $M$. In our case, $M =
\Bun_{G,(z_i),\infty}^{(m_i),m_\infty}$ and $T^*M$ is isomorphic to
the stack of Higgs bundles.) The surjectivity of \eqref{well-defined
map} follows in the same way as above. Thus, we can obtain another
proof of \thmref{quantized} this way.\qed
\end{remark}

\subsection{Hitchin systems with non-trivial characters}

In \secref{multi-point} we have generalized the definition of
$\ZZ_{(z_i),\infty}^{(m_i),m_\infty}(\g)$ to allow for non-trivial
characters $\chi_i: t^{m_i} \g[[t]] \to \C$ and $\chi_\infty:
t^{m_\infty} \g[[t]] \to \C$. A similar generalization is also
possible classically. In order to define it, we apply a Hamiltonian
reduction.

Let us recall that any character $\chi: t^{m} \g[[t]] \to \C$ is
necessarily trivial on the Lie subalgebra $t^{2m} \g[[t]]$ and hence
is determined by a character $t^{m} \g[[t]]/t^{2m} \g[[t]] \to \C$,
which is just an arbitrary linear functional on $t^{m} \g[[t]]/t^{2m}
\g[[t]]$ as this Lie algebra is abelian. We will write it in the form
$$
\chi = \sum_{k=-m-1}^{-2m} \chi_k t^k dt, \qquad \chi_k \in \g^*
$$
(with respect to the residue pairing).

Now consider the cotangent bundle $T^*
\MM_{G,(z_i),\infty}^{(2m_i),2m_\infty}$ which is isomorphic to
\eqref{open substack} with $m_i \mapsto 2m_i, m_\infty \mapsto
2m_\infty$. This is a Poisson manifold equipped with a Poisson action
of the commutative Lie group $\prod_{\al=1,\ldots,N,\infty}
G_{m_\al,2m_\al}$, where $G_{m,2m}$ is the Lie group of $t^{m}
\g[[t]]/t^{2m} \g[[t]]$. Now we apply the Hamiltonian reduction with
respect to this action and the one-point orbit in the dual space to
the Lie algebra of $\prod_{\al=1,\ldots,N,\infty} G_{m_\al,2m_\al}$
corresponding to $((\chi_i),\chi_\infty)$.

Let us denote the resulting Poisson manifold by ${\mc
M}_{G,(z_i),\infty;(\chi_i),\chi_\infty}^{(m_i),m_\infty}$. Its points
may be identified with collections $((g_i),g_\infty;\eta)$, where $g_i
\in G_{m_i}, g_\infty \in \ol{G}_{m_\infty}$ (where $G_m$ is the Lie
group of $\g_m$ and $\ol{G}_m$ is the Lie group of $\ol{\g}_m$) and
$\eta$ is a one-form \eqref{eta}, where $A_i(t-z_i)$ now has the form
$$
A_i(t-z_i) d(t-z_i) = \left( \sum_{k=-m_i}^{-1} A_{i,k} (t-z_i)^{k} +
\sum_{k=-m_i-1}^{-2m_i} \eta_k (t-z_i)^{k} \right) d(t-z_i), \qquad
A_{i,k} \in \g^*,
$$
and similarly for $A_\infty(t^{-1}) d(t^{-1})$. In other words, we
allow the polar parts of $\eta$ to have orders $2m_\al$, but fix
the $m_\al$ most singular terms in the expansion to be given by the
character $\on{Ad}_{g_\al}(\chi_\al)$.

Consider now the Hitchin map defined by formula \eqref{Hitchin
map}. Taking pull-backs of polynomial functions via the Hitchin map,
we obtain a Poisson commutative subalgebra
$\ol\A_{(z_i),\infty}^{(m_i),m_\infty}(\g)_{(\chi_i),\chi_\infty}$ in
the algebra of functions on the Poisson manifold ${\mc
M}_{G,(z_i),\infty;(\chi_i),\chi_\infty}^{(m_i),m_\infty}$
generalizing $\ol\ZZ_{(z_i),\infty}^{(m_i),m_\infty}(\g)$. This manifold
fibers over $\prod_{i=1}^N G_{m_i} \times \ol{G}_{m_\infty}$. In the
case when all $\chi_\al \equiv 0$ the fiber over the identify element
in $\prod_{i=1}^N G_{m_i} \times \ol{G}_{m_\infty}$ is a Poisson
submanifold isomorphic to $\prod_{i=1}^N \g_{m_i}^* \times
\ol\g_{m_\infty}^*$ with its Kirillov-Kostant structure. This means
that the restriction of the algebra of Hitchin Hamiltonians to this
fiber gives rise to a Poisson commutative subalgebra of
$\bigotimes_{i=1,\ldots,N} S(\g_{m_i}) \otimes S(\ol\g_{m_\infty})$.

But for general characters $(\chi_i),\chi_\infty$ the fiber over the
identity element is no longer a Poisson submanifold (in other words,
its defining ideal is not a Poisson ideal). This is due to the fact
that $\g_m$ acts non-trivially on non-zero characters $\chi:
t^m\g[[t]] \to \C$. Therefore we cannot restrict the Hitchin system to
this fiber. The best we can do is to consider the commutative Poisson
subalgebra in the algebra of functions on the entire manifold ${\mc
M}_{G,(z_i),\infty;(\chi_i),\chi_\infty}^{(m_i),m_\infty}$. This is
analogous to the fact that we cannot realize the corresponding quantum
algebras
$\A_{(z_i),\infty}^{(m_i),m_\infty}(\g)_{(\chi_i),\chi_\infty}$ as
subalgebras of universal enveloping algebras, as we pointed out in
\secref{multi-point}.\footnote{This discussion suggests these algebras
may instead be realized as commutative subalgebras of algebras of
twisted differential operators on
$\MM_{G,(z_i),\infty}^{(m_i),m_\infty}$, with the twisting determined
by the characters $(\chi_i),\chi_\infty$.}

However, as in the quantum case, there one exception, namely, when the
characters $\chi_i \equiv 0$ for all $i=1,\ldots,N$ and the character
$\chi_\infty$ becomes zero when restricted to $t^{m+1}\g[[t]]$. It may
then take non-zero values on $\g \otimes t^m \simeq \g$, determined by
some $\chi \in \g^*$. In this case the fiber of our Poisson manifold
over the identity is Poisson and is isomorphic to $\prod_{i=1}^N
\g_{m_i}^* \times \ol\g_{m_\infty-1}^*$ with its Kirillov-Kostant
structure. This means that the restriction of the algebra of Hitchin
Hamiltonians to this fiber gives rise to a Poisson commutative
subalgebra of $\bigotimes_{i=1,\ldots,N} S(\g_{m_i}) \otimes
S(\ol\g_{m_\infty-1})$.

This Poisson algebra may be described in more concrete terms as
follows. Note that the fiber of ${\mc
M}_{G,(z_i),\infty;(0),\chi}^{(m_i),m_\infty}$ at the identity is the
space of one-forms
\begin{equation}    \label{eta1}
\eta = \left( \sum_{i = 1,\ldots,N} \sum_{k=-m_i}^{-1} A_{i,k}
(t-z_i)^{k} - \sum_{k=0}^{m_\infty-2} A_{\infty,-k-2} t^k - \chi
t^{m_\infty-1} \right) dt,
\end{equation}
where $A_{\al,k} \in \g^*$. In other words, we fix the leading
singular term of $\eta$ at the point $\infty$ to be equal to $\chi \in
\g^*$.

This space is thus isomorphic to $\prod_{i=1}^N \g_{m_i}^* \times
\ol\g_{m_\infty-1}^*$, and the algebra of polynomial functions on it
is isomorphic, as a Poisson algebra, to $\bigotimes_{i=1,\ldots,N}
S(\g_{m_i}) \otimes S(\ol\g_{m_\infty-1})$, as we discussed above. Now
we define $\ol\A_{(z_i),\infty}^{(m_i),m_\infty}(\g)_{(0),\chi}$ as the
Poisson subalgebra of $\bigotimes_{i=1,\ldots,N} S(\g_{m_i}) \otimes
S(\ol\g_{m_\infty-1})$ generated by the coefficients of the invariant
polynomials $\ol{P}_\ell(\eta)$ in front of the monomials in
$(t-z_i)^{-n_i} t^{n_\infty}$. One shows in the same way as above that
this is a Poisson commutative subalgebra of $\bigotimes_{i=1,\ldots,N}
S(\g_{m_i}) \otimes S(\ol\g_{m_\infty-1})$.

Now, we claim that the algebra
$\A_{(z_i),\infty}^{(m_i),m_\infty}(\g)_{(0),\chi}$, discussed at the
very end of \secref{multi-point}, is a quantization of the Poisson
algebra $\ol\A_{(z_i),\infty}^{(m_i),m_\infty}(\g)_{(0),\chi}$ for any
$\chi \in \g^*$. We will prove this in \secref{sec:quant} in the
special case when $N=1$, $m_1=1$ and $m_2=1$ (we will also set
$z_1=0$, but this is not essential). The proof in general is very
similar and will be omitted.

\medskip

We thus set out to prove that the algebra
$\A^{1,1}_{0,\infty}(\g)_{0,\chi}$, denoted by ${\mc A}_\chi$ in
\secref{Achi}, is a quantization of the algebra
$\ol\A^{1,1}_{0,\infty}(\g)_{0,\chi}$, just defined. Let us denote
the latter by $\ol{\A}_\chi$. By definition, ${\mc A}_\chi$ is a
commutative subalgebra of $U(\g)$, and $\ol{\A}_\chi$ is a Poisson
commutative subalgebra of $S(\g)$. Saying that ${\mc A}_\chi$ is a
quantization of $\ol{\A}_\chi$ simply means that
\begin{equation}    \label{quant of A}
\gr {\mc A}_\chi = \ol{\A}_\chi,
\end{equation}
where $\gr {\mc A}_\chi$ is the associated graded algebra of ${\mc
A}_\chi$ with respect to the PBW filtration on $U(\g)$.

First, we take a closer look at $\ol{\A}_\chi$ and show that it is
nothing but the shift of argument subalgebra introduced in
\cite{MF}.

\subsection{The shift of argument subalgebra $\ol{\A}_\chi$}

Let us fix $\chi\in\g^*$ and consider the space of one-forms
\eqref{eta1} in the special case when $N=1, z_1=0$, $m_1=1$ and
$m_2=1$. This space consists of the one-forms
$$
\eta = \frac{A}{t} - \chi, \qquad A \in \g^*,
$$
and is therefore isomorphic to $\g^*$. Following the above general
definition, we define $\ol{\A}_\chi =
\ol\A^{1,1}_{0,\infty}(\g)_{0,\chi}$ as the subalgebra of $S(\g) =
\on{Fun} \g^*$ generated by the coefficients of the polynomials in
$t$,
$$
\ol{P}_i\left( \frac{A}{t} - \chi \right), \qquad i=1,\ldots,\ell.
$$

Equivalently, $\ol{\A}_\chi$ is the subalgebra of $\on{Fun} \g^*$
generated by the iterated directional derivatives $D^{i}_\chi P$ of
invariant polynomials $P\in(S\g)^\g = \on{Fun} \g^*$ in the direction
$\chi$, where
\begin{equation}    \label{dir der}
D_\chi P(x)=\left.\frac{d}{du}\right|_{u=0}P(x+u\chi).
\end{equation}
This definition makes it clear that $\ol{\mc A}_\chi = \ol{\mc
A}_{c\chi}$ for any non-zero $c \in \C$. We also have
$\on{Ad}_g(\ol{\mc A}_\chi) = \ol{\mc A}_{\on{Ad}_g(\chi)}$ for any $g
\in G$.

Note that for any $P\in S\g$, $x\in\g^*$ and $u\in\IC$, we have
$$
P(x+u\chi)=\sum_{m\geq 0}\frac{u^m}{m!}D^{m}_\chi P(x).
$$
Therefore we find that $\ol{\A}_\chi$ may also be defined as the
subalgebra of $S\g$ generated by the shifted polynomials
\begin{equation}\label{eq:shift}
P_{u\chi}(x)=P(x+u\chi)
\end{equation}
where $P$ varies in $(S\g)^\g$ and $u\in\IC$.

The algebra $\ol{\A}_\chi$ was introduced in \cite{MF}, where it was
shown to be Poisson commutative and of maximal possible transcendence
degree $\dim\bb$, where $\bb$ is a Borel subalgebra of $\g$, provided
$\chi$ is regular semi-simple element of $\g^* \simeq \g$. Here and
below we identify $\g$ and $\g^*$ using a non-degenerate invariant
inner product on $\g$. These results were recently elucidated and
extended by B. Kostant to regular nilpotent elements $\chi$
\cite{Ko2}. We will show now that a slight generalization of the
arguments used in \cite{MF,Ko2} yields the same result for all regular
elements. We recall that $a \in \g$ is called {\em regular} if its
centralizer in $\g$ has the smallest possible dimension; namely,
$\ell$, the rank of $\g$.

\begin{thm}    \label{poly}
Let $\chi$ be any regular element of $\g^* \simeq \g$. Then
$\ol{\A}_\chi$ is a free polynomial algebra in $\dim\bb$
generators
\begin{equation}    \label{iterated}
D^{n_i}_\chi \ol{P}_i, \qquad i=1,\ldots,\ell; \; n_i=0,\ldots,d_i,
\end{equation}
where $\ol{P}_i$ is a generator of $S(\g)^{\g}$ of degree $d_i+1$ and
$D_\chi$ is the derivative in the direction of $\chi$ given by formula
\eqref{dir der}.
\end{thm}

\begin{proof}
It is clear from the definition that $\ol{\A}_\chi$ is generated by
the iterated derivatives \eqref{iterated}. Therefore we need to prove
that the above polynomial functions on $\g^*$ are algebraically
independent for regular $\chi$.  For that it is sufficient to show
that the differentials of these functions at a particular point $\eta
\in \g^*$ are linearly independent. Note that we have
$$
\ol{P}_i(x+u\chi) = \sum_{m=0}^{d_i+1}\frac{u^m}{m!}D^{m}_\chi
\ol{P}_i(x),
$$
where the last coefficient, $D^{d_i+1}_\chi \ol{P}_i(x)$, is a
constant. Thus, we need to show that the cofficients $C_{i,k} =
C_{i,k}(\eta) \in \g$ appearing in the $u$-expansion of the
differential $d\ol{P}_i(x+u\chi)|_{x=\eta}$ of $\ol{P}_i(x+u\chi)$
(considered as a function of $x \in \g^*$ with fixed $u$) at $\eta \in
\g^*$,
$$
d\ol{P}_i(x+u\chi) = \sum_{k=0}^{d_i} C_{i,k} u^k,
$$
are linearly independent.

Note that for any $\eta \in \g^*$ and any $\g$-invariant function
$\ol{P}$ on $\g^*$ we have $[\eta,d\ol{P}(\eta)] = 0$. Therefore we
have
$$
[\eta + u \chi,d\ol{P}_i(\eta+u\chi)] = 0.
$$
Expanding the last equation in powers of $u$, we obtain the following
system (see Lemma 6.1.1 of \cite{MF}):
\begin{align}
[\eta,C_{i,0}] &= 0, \notag \\
[\eta,C_{i,k}] + [\chi,C_{i,k-1}] &= 0, \qquad 0<k<d_i, \label{recurs}
\\ [\chi,C_{i,d_i}] &= 0. \notag
\end{align}

Let us first prove that the elements $C_{i,k}$ are linearly
independent when $\chi$ is a regular nilpotent element, following
an idea from the proof of this theorem for regular semi-simple $\chi$
given in \cite{MF}.

Let $\{ e,2\crho,f \}$ be an $\sw_2$-triple in $\g$ such that $e =
\chi$. We choose the element $\crho$ as our $\eta$ -- this is the
point at which we evaluate the $C_{i,k}$'s. Under the adjoint action
of $\crho$ the Lie algebra $\g$ decomposes as follows:
$$
\g = \bigoplus_{i=-h}^h \g_i,
$$
where $h$ is the Coxeter number. Here $\g_0$ is the Cartan subalgebra
containing $\crho$ and $e$ is an element of $\g_1$.

Let
$$
W_j = \on{span}\{ C_{i,k} \, | \, i=1,\ldots,\ell; \;
k=0,\ldots,j \},
$$
where we set $C_{i,k}=0$ for $k>d_i$. We will now prove that the
following equality is true for all $j \geq 0$:
\begin{equation}    \label{Wj}
W_j = \bigoplus_{i=0}^j \g_i.
\end{equation}

Let us prove it for $j=0$. The first equation in
\eqref{recurs} implies that $C_{i,0}, i=1,\ldots,\ell$, belong to
$\g_0$. Furthermore, since $\crho$ is regular, they span
$\g_0$, by \cite{Ko}. Therefore $W_0 = \g_0$.

Now suppose that we have proved \eqref{Wj} for $j=0,\ldots,m$. Let us
prove it for $j=m+1$.

According to the equations \eqref{recurs}, we have
$$
\on{ad} \crho \cdot W_{m+1} = \on{ad} e \cdot W_m.
$$
It follows from general results on representations of $\sw_2$ that the
map
$$
\on{ad} e: \g_i \to \g_{i+1}
$$
is surjective for all $i \geq 0$. Therefore, using our inductive
assumption, we obtain that $\on{ad} e \cdot W_m =
\bigoplus_{i=1}^{m+1} \g_i$. Since $\on{ad} \crho$ is invertible on
this space and the kernel of $\on{ad} \crho$ on the entire $\g$ is
equal to $\g_0$, we find that $W_{m+1}$ is necessarily contained in
$\bigoplus_{i=0}^{m+1} \g_i$ and its projection onto
$\bigoplus_{i=1}^{m+1} \g_i$ along $\g_0$ is surjective. On the other
hand, $W_{m+1}$ contains $C_{i,0}, i=1,\ldots,\ell$, and therefore
contains $\g_0$. Hence we obtain the equality \eqref{Wj} for
$j=m+1$. This completes the inductive step and hence proves \eqref{Wj}
for all $j \geq 0$.

Setting $j=h$, we obtain that the span of the elements $C_{i,n_i},
i=1,\ldots,\ell; n_i=0,\ldots,d_i$, is equal to $\bigoplus_{i=0}^h
\g_h$, which is a Borel subalgebra $\bb$ of $\g$. Its dimension is
equal to the number of these elements, which implies that they are
linearly independent. This proves that the polynomials
\eqref{iterated} are algebraically independent, and so $\ol{\mc
A}_\chi$ is a free polynomial algebra in $\dim \bb$ generators, if
$\chi$ is a regular nilpotent element.

Now we derive that the same property holds for $\ol{\mc A}_\chi$,
where $\chi$ is an arbitrary regular element, following an argument
suggested to us by L. Rybnikov (see also \cite{Ko2}). Consider the
subset $S \subset \g^* \simeq \g$ of all elements satisfying the
property that the polynomials \eqref{iterated} are algebraically
independent. It is clear that this is a Zariski open subset of $\g$,
and it is non-empty since it contains regular nilpotent elements. Let
$\chi$ be a regular element of $\g^* \simeq \g$. Recall that we have
$\ol{\mc A}_\chi = \ol{\mc A}_{c\chi}$ for any non-zero $c \in \C$ and
$\on{Ad}_g(\ol{\mc A}_\chi) = \ol{\mc A}_{\on{Ad}_g(\chi)}$ for any $g
\in G$, where $G$ is the adjoint group of $\g$. This implies that $S$
is a conic subset of $\g$ that is invariant under the adjoint action
of $G$. According to \cite{Ko}, the adjoint orbit of any regular
element of $\g$ contains an element of the form $e + A$, where $e$ is
a regular nilpotent element of an $\sw_2$-triple $\{ e,2\crho,f \}$,
and $A \in \bigoplus_{i\leq 0} \g_i$. Hence we may restrict ourselves
to the elements $\chi$ of this form.

Suppose that some $\chi$ of this form does not belong to the set
$S$. Then neither do the elements
$$
\chi_c = c \on{Ad} \crho(c)^{-1}(\chi), \qquad c \in \C^\times,
$$
where $\crho: \C^\times \to G$ is the one-parameter subgroup of $G$
corresponding to $\crho \in \g$. But then the limit of $\chi_c$ as $c
\to 0$ should not be in $S$. However, this limit is equal to
$e$, which is a regular nilpotent element, and hence belongs to
$S$. This is a contradiction, which implies that all regular elements
of $\g$ belong to $S$. This completes the proof.
\end{proof}

The Poisson commutativity of $\ol{\A}_\chi$ follows from the
above general results about the commutativity of the Hitchin
Hamiltonians. It also follows from Theorem \ref{th:quantisation}
below.

\subsection{The quantization theorem}    \label{sec:quant}

In \cite{Vin} the problem of the existence of a quantization of
$\ol{\A}_\chi$ was posed: does there exist a commutative subalgebra
${\mc A}_\chi$ of $U(\g)$ which satisfies \eqref{quant of A}?

Such a quantization has been constructed for $\g$ of classical types
in \cite{NO}, using twisted Yangians, and for $\g=\sw_n$ in
\cite{Tar}, using the symmetrization map, and in \cite{CT}, using
explicit formulas. In \cite{Ryb1} it was shown that, if it exists, a
quantization of $\ol{\mc A}_\chi$ is unique for generic $\chi$.

Recently, the quantization problem was solved in \cite{Ryb2} for any
simple Lie algebra $\g$ and any regular semi-simple $\chi \in
\g^*$. More precisely, it was shown in \cite{Ryb2} that the algebra
${\mc A}_\chi$, constructed in essentially the same way as in
\secref{Achi}, is a quantization of $\ol{\mc A}_\chi$ for any regular
semi-simple $\chi$.

We will now prove that ${\mc A}_\chi$ is a quantization of $\ol{\mc
A}_\chi$ for any regular $\chi \in \g^*$. First, we prove the
following statement, which is also implicit in \cite{Ryb2}.

\begin{prop}    \label{contains}
The algebra $\ol{\mc A}_\chi$ is contained in $\on{gr} {\mc
  A}_\chi$ for any $\chi \in \g^*$.
\end{prop}

\proof
By definition (see \secref{Achi}) and Proposition \ref{pr:modes}, the
algebra $\A_\chi$ is generated by the coefficients of the Laurent
expansion of $\id\otimes\chi(\Phi _{u,(0)}^{1,2}(A))$ at $u=0$, where
$A$ ranges over a system of generators of $\zz(\ghat)$. Therefore we
need to show that the symbol of each of these coefficients belongs to
$\ol{\mc A}_\chi$.

Since $\Phi_{u,(0)}^{1,2}(TA)=\partial_u\Phi_{u,(0)}^{1,2}(A)$, where
$T\in\DerO$ is the translation operator given by \eqref{translation},
it suffices in fact to consider a system of generators of $\zz(\ghat)$
as an algebra endowed with the derivation $T$ (equivalently, as a
commutative vertex algebra).

By Theorem \ref{thm:FF}, $\gr(\zz(\ghat))=S(\hatg_-)^{\g[[t]]}$ and,
by Theorem \ref{thm:BD}, $S(\hatg_-)^{\g[[t]]}$ is generated, as an
algebra endowed with the derivation $T$, by the polynomials $P_{-1}$,
as $P$ varies in $(S\g)^\g$. It follows that $\zz(\ghat)$ is
generated, as an algebra with a derivation, by elements $\Pi_{-1}$,
where $\Pi_{-1}$ is such that its symbol is $P_{-1}$, and $P$ varies
in $(S\g)^\g$.

The computation of the symbol of
$\id\otimes\chi(\Phi_{u,(0)}^{1,2}(\Pi_{-1}))$ will be carried out in
Lemma \ref{pr:shift of argument} below. To that end, it will be
necessary to choose the lifts $\Pi_{-1}$ in the following way.  First,
we will only consider homogeneous generators of $S\g$. If $P\in
S^d\g$, then $P_{-1}$ is of degree $d$ \wrt the $\Z_+$-grading on
$S(\hatg_-)$ introduced in \secref{grad}. We may, and will,
assume that the lift $\Pi_{-1}$ to $\zz(\ghat)$ is also homogeneous of
degree $d$.

The proof of \propref{contains} is now completed by the following
calculation. Let $\Pi\in U(\hatg_-)$ be an element of order $d\in\IN$
and degree $d$ with respect to the $\Z_+$-grading on
$U(\hatg_-)$. Consider its symbol $\sigma(\Pi)\in S^d(\g\otimes
t^{-1})$ as an element $P\in S^d\g$ via the linear isomorphism
$\g\otimes t^{-1}\cong\g$. For each $n \in \Z$ denote by
$\left((\id\otimes\chi) \circ\Phi_{u,(0)}^{1,2}(\Pi)\right)_{n}$ the
coefficient of $u^{n}$ in $(\id\otimes\chi) \circ
\Phi_{u,(0)}^{1,2}(\Pi)$.

\begin{lem}\label{pr:shift of argument}
For any $k \geq 0$ the following holds: if $D_\chi^k P \neq 0$,
then the symbol of $(\id\otimes\chi) \circ
\Phi_{u,(0)}^{1,2}(\Pi)$ is equal to
$$
\sigma\left((\id\otimes\chi)
\circ\Phi_{u,(0)}^{1,2}(\Pi)\right)_{-d+k} = \frac{(-1)^k}{k!}
D_\chi^k P.
$$
\end{lem}

\proof By assumption, $\Pi$ is a linear combination of lexicographically
ordered monomials of the form
$$\Pi_{(a_i,n_i,b_j)}=
J^{a_1}_{-n_1}\cdots J^{a_p}_{-n_p}J^{b_1}_{-1}\cdots J^{b_q}_{-1},$$
where $n_1\odots{\geq}n_p\geq 2$, 
$$\ord(\Pi_{(a_i,n_i,b_j)})=p+q\leq d\aand
\deg(\Pi_{(a_i,n_i,b_j)})=\sum_i n_i+q=d.$$
Note that if $p\neq 0$, then
$$\ord(\Pi_{(a_i,n_i,b_j)})=p+q<\sum_i n_i+q=d,$$
since $n_i\geq 2$. Thus, only the monomials with $p=0$ contribute to
$\sigma(\Pi)$, so that
$$\sigma(\Pi)=
\sum_{1\leq b_1\odots{\leq}b_d\leq\dim\g}
\alpha_{b_1,\ldots,b_d} \ol{J}^{b_1}_{-1}\cdots \ol{J}^{b_d}_{-1}$$
for some constants $\alpha_{b_1,\ldots,b_d}\in\IC$, where the product
is now that of $S(\hatg_-)$ and
$$P=
\sum_{1\leq b_1\odots{\leq}b_d\leq\dim\g}
\alpha_{b_1,\ldots,b_d}J^{b_1}\cdots J^{b_d}.
$$
For $P = \ol{P}_i$ this is in fact the statement of \lemref{lexico}.

By \eqref{Psi},
\begin{equation*}
\begin{split}
(\id\otimes\chi) \circ\Phi_{u,(0)}^{1,2}(\Pi_{(a_i,n_i,b_j)})
&=
\left(\frac{J^{b_q}}{u}-\chi(J^{b_q})\right)\cdots
\left(\frac{J^{b_1}}{u}-\chi(J^{b_1})\right)
\frac{J^{a_p}}{u^{n_p}}\cdots
\frac{J^{a_1}}{u^{n_1}}\\
&=
u^{-d}
\left(J^{b_q}-u\chi(J^{b_q})\right)\cdots
\left(J^{b_1}-u\chi(J^{b_1})\right)
J^{a_p}\cdots J^{a_1}.
\end{split}
\end{equation*}
The coefficient of $u^{-d+k}$ in $(\id\otimes\chi)
\circ\Phi_{u,(0)}^{1,2}(\Pi_{(a_i,n_i,b_j)})$ is therefore an element
of $U\g$ of order $\leq p+q-k$. If $p\neq 0$, this is strictly less
than $\sum_i n_i+q-k=d-k$ since $n_i\geq 2$. If, on the other hand,
$p=0$, then the coefficient of $u^{-d+k}$ is proportional to a
lexicographically ordered monomial in $U\g$ of order $d$. It therefore
follows that $\sigma\left((\id\otimes\chi)
\circ\Phi_{u,(0)}^{1,2}(\Pi)\right)_{-d+k}$ is the coefficient of
$u^k$ in
\begin{multline*}
\sum_{1\leq b_1\odots{\leq}b_d\leq\dim\g}\alpha_{b_1,\ldots,b_d}
\left(J^{b_d}-u\chi(J^{b_d})\right)\cdots
\left(J^{b_1}-u\chi(J^{b_1})\right)
=
u^{-d}P(\cdot-u\chi),
\end{multline*}
provided that this coefficient is non-zero (otherwise, we would obtain
the symbol of the first non-zero lower order term). This implies the
statement of the lemma.\qed

\smallskip

\begin{thm}\label{th:quantisation}
For regular $\chi \in \g^*$ we have $\on{gr} {\mc A}_\chi = \ol{\mc
A}_\chi$, and so the commutative algebra $\A_\chi\subset U\g$ is a
quantization of the shift of argument subalgebra $\ol{\A}_\chi\subset
S\g$.
\end{thm}

\begin{proof}
We know that $\ol{\A}_\chi$ is generated by
$$
D_\chi^{k_i} \ol{P}_i, \qquad i=1,\ldots,\ell; \; 0 \leq k_i \leq d_i.
$$
On the other hand, it follows from the definition of ${\mc A}_\chi$
that it is generated by
$$
\left((\id\otimes\chi)
\circ\Phi_{u,(0)}^{1,2}(S_{i})\right)_{-d_i-1+k_i}, \qquad
i=1,\ldots,\ell; \; 0 \leq k_i \leq d_i.
$$
If $\chi$ is regular, then each $D_\chi^{k_i} \ol{P}_i \neq 0$, by
\thmref{poly}, and therefore \lemref{pr:shift of argument} implies that
the generators of $\ol{\A}_\chi$ are equal to the symbols of the
generators of ${\mc A}_\chi$, up to non-zero scalars. In addition, the
generators $D_\chi^{k_i} \ol{P}_i$ of $\ol{\A}_\chi$ are algebraically
independent for regular $\chi$, by \thmref{poly}. Therefore, again
applying \lemref{pr:shift of argument}, we obtain that the symbol of
any non-zero element of ${\mc A}_\chi$ is a non-zero element of
$\ol{\mc A}_\chi$. Therefore $\on{gr} {\mc A}_\chi \subset \ol{\mc
A}_\chi$. Combining this with \propref{contains}, we obtain the
assertion of the theorem.
\end{proof}

For non-regular $\chi$, \propref{contains} implies that all elements
of $\ol{\mc A}_\chi$ may be quantized, i.e., lifted to commuting
elements of ${\mc A}_\chi$, but this still leaves open the possibility
that the quantum algebra ${\mc A}_\chi$ is larger than its classical
counterpart $\ol{\mc A}_\chi$ (we thank L. Rybnikov for pointing this
out). However, we conjecture that this never happens:

\begin{conj}
We have $\on{gr} {\mc A}_\chi = \ol{\mc A}_\chi$ for all $\chi \in
\g^*$.
\end{conj}

By \thmref{th:quantisation}, this conjecture holds in the regular
case, and it also holds in the most degenerate case when $\chi=0$. In
that case ${\mc A}_\chi$ is the center $Z(\g)$ of $U(\g)$ and $\ol{\mc
A}_\chi = S(\g)^{\g}$.

\subsection{The DMT Hamiltonians}    \label{DMT}

As we saw in \secref{gaudin model}, one can write down explicit
formulas for quadratic generators of the algebra
$\ZZ^{(1),1}_{(z_i),\infty}(\g)$; these are the original Gaudin
Hamiltonians $\Xi_i, i=1,\ldots,N$.

In this section we determine quadratic generators of the algebra ${\mc
A}_\chi$ for a regular semi-simple $\chi \in \g$. They turn out to be
none other than the DMT Hamiltonians discussed in the Introduction.

Let $\h$ be the Cartan subalgebra of $\g$ containing $\chi$ and
$\Delta\subset\h^*$ the root system of $\g$. For each
$\alpha\in\Delta$, let $\sl{2}^\alpha=\<e_\alpha,f_
\alpha,h_\alpha\>\subset\g$ be the corresponding three-dimensional
subalgebra and denote by
$$C_\alpha=\frac{(\alpha,\alpha)}{2} \left(e_\alpha
f_\alpha+f_\alpha e_\alpha\right)$$ its truncated Casimir operator
\wrt the restriction to $\sl{2}^\alpha$ of a fixed non-degenerate
invariant inner product $(\cdot,\cdot)$ on $\g$.  Note that
$C_\alpha$ is independent of the choice of the root vectors
$e_\alpha,f_\alpha$ and satisfies
$C_{-\alpha}=C_{\alpha}$. Let
$$\hreg=\h\setminus\bigcup_{\alpha\in\Delta}\Ker(\alpha)$$ be the
variety of regular elements in $\h$, $V$ a $\g$-module and
$\VV=\hreg\times V$ the holomorphically trivial vector bundle over
$\hreg$ with fiber $V$. Millson and Toledano Laredo \cite{MTL,TL},
and, independently, De Concini (unpublished; a closely related
connection was also considered in \cite{FMTV}) introduced the
following holomorphic connection on $\VV$:
\begin{equation*}
\nablak
=
d-\frac{h}{2}\sum_{\alpha\in\Delta}
\frac{d\alpha}{\alpha}\cdot C_\alpha
=
d-h\sum_{\alpha\in\Delta_+}
\frac{d\alpha}{\alpha}\cdot C_\alpha,
\end{equation*}
where $\Delta=\Delta_+\sqcup\Delta_-$ is the partition of $\Delta$
into positive and negative roots determined by a choice of simple
roots $\alpha_1, \ldots,\alpha_n$ of $\g$, and proved that $\nablak$
is flat for any value of the complex parameter $h$. Define, for every
$\gamma\in\h$, the function $T_\gamma:\hreg\rightarrow\End(V)$ by
\begin{equation}\label{eq:T Ham}
T_\gamma(\chi)=
\sum_{\alpha\in\Delta_+}\frac{\alpha(\gamma)}{\alpha(\chi)}C_\alpha
\end{equation}
and note that $T_{a\gamma+a'\gamma'}=a T_\gamma+a' T_\gamma'$.
Then, the flatness of $\nablak$ is equivalent to the equations
$$[D_\gamma-hT_\gamma,D_{\gamma'}-hT_{\gamma'}]=0$$
for any $\gamma,\gamma'\in\h$, where $D_\gamma f(x)=\left.\frac{d}
{dt}\right|_{t=0}f(x+t\gamma)$. Dividing by $h$ and letting $h$ tend to
infinity implies that, for a fixed $\chi\in\hreg$ and any $\gamma,\gamma'
\in\h$,
$$[T_\gamma(\chi),T_{\gamma'}(\chi)]=0.$$

\smallskip

In \cite{Vin} it was shown that the algebra $\A_\chi$ contains the
Hamiltonians $T_\gamma(\chi)$. We now give a proof of this result for
completeness.

Let us identify $\g$ and $\g^*$ using the invariant inner product
$(\cdot,\cdot)$. For each homogeneous polynomial $p\in S\g^*$ of
degree $d$ and $i=1, \ldots,d$, let $p^{(i)}_\chi\in S^i\g^*$ be
defined by
$$p^{(i)}_\chi=\frac{D^{d-i}_\chi p}{(d-i)!},$$ where $D_\chi
f(x)=\left.\frac{d}{dt}\right|_{t=0}f(x+t\chi)$. Note that if $p$ is
invariant under $\g$, then $p^{(i)}_\chi$ is invariant under the
centralizer $\g^\chi\subseteq\g$ of $\chi$. The following result shows
that, when $\chi \in\hreg$, the algebra $\ol{\mc A}_\chi$ contains the
symbols $\ol{T}_\gamma(\chi)$ of the Hamiltonians $T_\gamma(\chi)$
defined by \eqref{eq:T Ham}.

\begin{prop}\label{th:generation}
Assume that $\chi\in\hreg$. Then,
\begin{enumerate}
\item Let $p\in(S\g^*)^{\g}$ be homogeneous and $q\in(S\h^*)^W$ its
restriction to $\h$. Then
$$p_\chi^{(2)}=q_\chi^{(2)}+\ol{T}_{\gamma_q}(\chi),$$
where $q_\chi^{(2)}=D^{d-2}_\chi q/(d-2)!$ and $\gamma_q=dq(\chi)
\in T_\chi^*\h\cong\h$ is the differential of $q$ at $\chi$ and the
cotangent space $T_\chi^*\h$ is identified with $\h$ by means of the
form $(\cdot,\cdot)$.
\item As $p$ varies over the homogeneous elements of $(S\g^*)^{\g}$,
$\gamma_q$ ranges over the whole of $\h$.
\end{enumerate}
\end{prop}

We will need the following

\begin{lem}\label{le:D=d}
If $p\in S\g^*$ is $\g$-invariant and $\chi\in\h$, then, for any $\alpha
\in\Delta$,
$$\alpha(\chi)D_{e_\alpha}D_{f_\alpha}\medspace p^{(2)}_\chi=
D_{h_\alpha}p(\chi).$$
\end{lem}
\proof
Let $P\in(\g^*)^{\otimes d}$ be the symmetric, multilinear function on
$\g$ defined by $p(x)=P(x,\ldots,x)$. The $\g$-invariance of $p$ implies
that
\begin{equation*}
\begin{split}
0
&=
\ad^*(e_\alpha)P(f_\alpha,\underbrace{\chi,\ldots,\chi}_{d-1})\\
&=
-P(h_\alpha,\underbrace{\chi,\ldots,\chi}_{d-1})+
(d-1)\alpha(\chi)P(f_\alpha,e_\alpha,
\underbrace{\chi,\ldots,\chi}_{d-2})\\ &=
-\frac{1}{d}D_{h_\alpha}p(\chi)+
\alpha(\chi)\frac{d-1}{d!}D_{e_\alpha}D_{f_\alpha}D_\chi^{d-2}p,
\end{split}
\end{equation*}
whence the claimed result. \qed\\

\noindent {\em Proof of \propref{th:generation}}.\\ (i) Let
$\nu:\g\rightarrow\g^*$ be the isomorphism of $\g$-modules induced by
$(\cdot,\cdot)$. Since $p_\chi^{(2)}\in S^2\g^*$ is invariant under
$\h=\g^\chi$, and $(e_\alpha,f_\alpha)=2/(\alpha,\alpha)$, we have
\begin{equation*}
\begin{split}
p_\chi^{(2)}
&=
\left.p_\chi^{(2)}\right|_\h+
\sum_{\alpha\in\Delta_+}\left(\frac{(\alpha,\alpha)}{2}\right)^2
\nu(e_\alpha)\nu(f_\alpha)D_{e_\alpha}D_{f_\alpha}p_{\chi}^{(2)}\\
&=
\left.p_\chi^{(2)}\right|_\h+
\sum_{\alpha\in\Delta_+}
\left(\frac{(\alpha,\alpha)}{2}\right)^2
\nu(e_\alpha)\nu(f_\alpha)\frac{D_{h_\alpha}p(\chi)}{\alpha(\chi)}\\
&=
\left.p_\chi^{(2)}\right|_\h+
\sum_{\alpha\in\Delta_+}
\frac{(\alpha,\alpha)}{2}
\nu(e_\alpha)\nu(f_\alpha)\frac{dp(\chi)(\nu^{-1}(\alpha))}
   {\alpha(\chi)},
\end{split}
\end{equation*}
where the second equality follows from lemma \ref{le:D=d}.\\ (ii) As
$p$ ranges over the homogeneous elements of $(S\g^*)^ \g$, $q$ ranges
over those of $(S\h^*)^W=\IC[c_1,\ldots,c_n]$. The differential of the
latter span $T_\chi^*\h$ since the Jacobian of $c_1,\ldots,c_n$ at
$\chi\in\h$ is proportional to $\prod_{\alpha\in\Delta_+}
\alpha(\chi)$ \cite[V.5.4]{Bo} and is therefore non-zero since $\chi$
is regular. \qed

\smallskip

Now, according to \thmref{th:quantisation}, any element of $\ol{\mc
A}_\chi$ may be lifted to ${\mc A}_\chi$. We already know that $\h
\subset {\mc A}_\chi$. An arbitrary lifting of $\ol{T}_{\gamma}(\chi)
\in \ol{\mc A}_\chi$ to $U(\g)$ is equal to $T_\gamma(\chi) + J$,
where $J \in \g$. But the lifting to ${\mc A}_\chi$ has to commute
with $\h$. Therefore $J \in \h$ and so $T_\gamma(\chi)$ itself belongs
to ${\mc A}_\chi$ for all $\gamma \in \h$.

\section{Recollections on opers}    \label{opers}

In order to describe the universal Gaudin algebra and its various
quotients introduced in \secref{construction} we need to recall the
description of $\zz(\G)$. According to \cite{FF:gd,F:wak}, $\zz(\G)$
is identified with the algebra $\on{Fun} \on{Op}_{^L G}(D)$ of
(regular) functions on the space $\on{Op}_{^L G}(D)$ of $^L G$-opers
on the disc $D = \on{Spec} \C[[t]]$. Here $^L G$ is the Lie group of
adjoint type corresponding to the Lie algebra $^L \g$ whose Cartan
matrix is the transpose of that of $\g$. Note that $^L G$ is the
Langlands dual group of the connected simply-connected Lie group $G$
with Lie algebra $\g$.

In this section we collect results on opers that we will need (for a
more detailed exposition, see \cite{F:book}). Then in the next section
we describe $\zz(\G)$ and the Gaudin algebras in terms of opers.

\subsection{Definition of opers}    \label{defn of opers}

Let $G$ be a simple algebraic group of adjoint type, $B$ a Borel
subgroup and $N = [B,B]$ its unipotent radical, with the corresponding
Lie algebras $\n \subset \bb\subset \g$. The quotient $H = B/N$ is a
torus. Choose a splitting $H \to B$ of the homomorphism $B \to H$ and
the corresponding splitting $\h \to \bb$ at the level of Lie
algebras. Then we have a Cartan decomposition $\g = \n_- \oplus
\h \oplus \n$. We choose generators $\{ e_i \}, i=1,\ldots,\ell$,
of $\n$ and generators $\{ f_i \}, i=1,\ldots,\ell$ of $\n_-$
corresponding to simple roots, and denote by $\crho \in \h$ the sum of
the fundamental coweights of $\g$. Then we have the following
relations:
$$
[\crho,e_i] = e_i, \quad [\crho,f_i] = -f_i, \qquad i=1,\ldots,\ell.
$$

A $G$-oper on a smooth curve $X$ (or a disc $D \simeq \on{Spec}
\C[[t]]$ or a punctured disc $D^\times = \on{Spec} \C\ppart$) is by
definition a triple $(\F,\nabla,\F_B)$, where $\F$ is a principal
$G$-bundle $\F$ on $X$, $\nabla$ is a connection on $\F$ and $\F_B$ is
a $B$-reduction of $\F$ such that locally on $X$, in the analytic or
\'etale topology, it has the following form. Choose a coordinate $t$
and a trivialization of $\F_B$ on a sufficiently small open subset of
$X$ (in the analytic or \'etale topology, respectively) over which
$\F_B$ may be trivialized. Then the connection operator
$\nabla_{\pa_t}$ corresponding to the vector field $\pa_t$ has the
form
\begin{equation} \label{form of nabla}
\nabla_{\pa_t} = \pa_t + \sum_{i=1}^\ell \psi_i(t) f_i + {\mb v}(t),
\end{equation}
where each $\psi_i(t)$ is a nowhere vanishing function, and ${\mb
v}(t)$ is a $\bb$-valued function. The space of $G$-opers on $X$ is
denoted by $\on{Op}_G(X)$.

This definition is due to A. Beilinson and V. Drinfeld \cite{BD} (in
the case when $X$ is the punctured disc opers were first introduced in
\cite{DS}).

In particular, suppose that $U = \on{Spec} R$ and $t$ is a coordinate
on $U$. It is clear what this means if $U = D$ or $D^\times$, and if
$U$ is an affine curve with the ring of functions $R$, then $t$ is an
\'etale morphism $U \to {\mathbb A}^1$ (for example, if $U = \C^\times
= \on{Spec} \C[t,t^{-1}]$, then $t^n$ is a coordinate for any non-zero
integer $n$). Then $\on{Op}_G(U)$ has the following explicit
realization: it is isomorphic to the quotient of the space of
operators of the form\footnote{In order to simplify notation, from now
on we will write $\nabla$ for $\nabla_{\pa_t}$.}
\begin{equation} \label{another form of nabla1}
\nabla = \pa_t + \sum_{i=1}^\ell \psi_i(t) f_i + {\mb v}(t), \qquad
{\mb v}(t) \in \bb(R),
\end{equation}
where each $\psi_i(t) \in R$ is a nowhere vanishing function, by the
action of the group $B(R)$. Recall that the gauge transformation of an
operator $\pa_t+A(t)$, where $A(t) \in \g(R)$ by $g(t) \in G(R)$ is
given by the formula
$$
g \cdot (\pa_t + A(t)) = \pa_t + g A(t) g^{-1} - \pa_t g \cdot g^{-1}.
$$

There is a unique element $g(t) \in H(R)$ such that the gauge
transformation of the operator \eqref{another form of nabla1} by
$H(R)$ has the form
\begin{equation} \label{another form of nabla}
\nabla = \pa_t + \sum_{i=1}^\ell f_i + {\mb v}(t), \qquad
{\mb v}(t) \in \bb(R).
\end{equation}
This implies that $\on{Op}_G(U)$ is isomorphic to the quotient of the
space of operators of the form \eqref{another form of nabla} by the
action of the group $N(R)$.

\subsection{Canonical representatives}    \label{can repr}

Set
$$
p_{-1} = \sum_{i=1}^\ell f_i.
$$
The operator $\on{ad} \crho$ defines the principal gradation on $\bb$,
with respect to which we have a direct sum decomposition $\bb =
\bigoplus_{i\geq 0} \bb_i$.  Let $p_1$ be the unique element of degree
1 in $\n$ such that $\{ p_{-1},2\rv,p_1 \}$ is an $\sw_2$-triple. Let
$$
V_{\can} = \oplus_{i \in E} V_{\can,i}
$$
be the space of $\on{ad}
p_1$-invariants in $\n$, decomposed according to the principal
gradation. Here
$$
E = \{ d_1,\ldots,d_\ell \}
$$
is the set of exponents of $\g$ (see \cite{Ko}). Then $p_1$ spans
$V_{\on{can},1}$. Choose a linear generator $p_j$ of $V_{\can,d_j}$
(if the multiplicity of $d_j$ is greater than one, which happens only
in the case $\g=D_{2n}, d_j=2n$, then we choose linearly independent
vectors in $V_{\on{can},d_j}$). The following result is due to
Drinfeld and Sokolov \cite{DS} (the proof is reproduced in the proof
of Lemma 2.1 of \cite{F:opers}).

\begin{lem} \label{free}
The gauge action of $N(R)$ on the space of operators of the form
\eqref{another form of nabla} is free, and each gauge equivalence
class contains a unique operator of the form $\nabla = \pa_t + p_{-1}
+ {\mathbf v}(t)$, where ${\mathbf v}(t) \in V_{\can}(R)$, so that we
can write
\begin{equation} \label{coeff fun}
{\mathbf v}(t) = \sum_{j=1}^\ell v_j(t) \cdot p_j.
\end{equation}
\end{lem}

Let $x$ be a point of a smooth curve $X$ and $D_x = \on{Spec} \OO_x,
D^\times_x = \on{Spec} \K_x$, where $\OO_x$ is the completion of the
local ring of $x$ and $\K_x$ is the field of fractions of $\OO_x$.
Choose a formal coordinate $t$ at $x$, so that $\OO_x \simeq \C[[t]]$
and $\K_x = \C\ppart$. Then the space $\on{Op}_G(D_x)$ (resp.,
$\on{Op}_G(D_x^\times)$) of $G$-opers on $D_x$ (resp., $D_x^\times$)
is the quotient of the space of operators of the form \eqref{form of
nabla} where $\psi_i(t) \neq 0$ take values in $\OO_x$ (resp., in
$\K_x$) and ${\mb v}(t)$ takes values in $\bb(\OO_x)$ (resp.,
$\bb(\K_x)$) by the action of $B(\OO_x)$ (resp., $B(\K_x)$).

To make this definition coordinate-independent, we need to specify the
action of the group of changes of coordinates on this space. Suppose
that $s$ is another coordinate on the disc $D_x$ such that $t =
\varphi(s)$. In terms of this new coordinate the operator
\eqref{another form of nabla} becomes
$$
\nabla_{\pa_t} = \nabla_{\varphi'(s)^{-1} \pa_s} = \varphi'(s)^{-1}
\pa_s + \sum_{i=1}^\ell f_i + \varphi'(s) \cdot {\mb v}(\varphi(s)).
$$
Hence we find that
$$
\nabla_{\pa_s} = \pa_s + \varphi'(s) \sum_{i=1}^\ell f_i + \varphi'(s)
\cdot {\mb v}(\varphi(s)).
$$

\subsection{Opers with singularities}    \label{reg sing}

A $G$-oper on $D_x$ with singularity of order $m$ at $x$ is by
definition (see \cite{BD}, Sect. 3.8.8) a $B(\OO_x)$-conjugacy class
of operators of the form
\begin{equation} \label{oper with RS1}
\nabla = \pa_t + \frac{1}{t^m} \left( \sum_{i=1}^\ell \psi_i(t) f_i +
{\mb v}(t) \right),
\end{equation}
where $\psi_i(t) \in \OO_x, \psi_i(0) \neq 0$, and ${\mb v}(t) \in
\bb(\OO_x)$. Equivalently, it is an $N(\OO_x)$-equivalence class of
operators
\begin{equation} \label{oper with RS}
\nabla = \pa_t + \frac{1}{t^m} \left( p_{-1} + {\mb v}(t) \right),
\qquad {\mb v}(t) \in \bb(\OO_x).
\end{equation}
We denote the space of such opers by $\on{Op}_G^{\leq m}(D_x)$.

Let $\mm_x = t\C[[t]]$ be the maximal ideal of $\OO_x = \C[[t]]$. Then
$(\mm_x/\mm_x^2)^*$ is naturally interpreted as the tangent space
$T_x$ to $x \in X$. Let $T_x^\times$ be the corresponding
$\C^\times$-torsor. Now consider the affine space
$$
\g/G := \on{Spec} (\on{Fun} \g)^G = \on{Spec} (\on{Fun} \h)^W =: \h/W.
$$
It carries a natural $\C^\times$-action. We denote by $(\g/G)_{x,n} =
(\h/W)_{x,n}$ its twist by the $\C^\times$-torsor
$(T_x^\times)^{\otimes n}$:
$$
(\g/G)_{x,n} = (T_x^\times)^{\otimes n} \underset{\C^\times}\times
\g/G.
$$

Define the $m$-{\em residue map}
\begin{equation}    \label{resm}
\on{res}_m: \on{Op}_G^{\leq m}(D_x) \to (\g/G)_{x,m-1} =
(\h/W)_{x,m-1}
\end{equation}
sending $\nabla$ of the form \eqref{oper with RS} to the image of
$p_{-1} + {\mb v}(0)$ in $\h/W$. It is clear that this map is
independent of the choice of coordinate $t$ (which was the reason for
twisting by $(T_x^\times)^{\otimes (m-1)}$). This generalizes the
definition of $1$-residue given in \cite{BD}, Sect 3.8.11. The
$1$-residue which takes values in $\g/G = \h/W$ and no twisting is
needed. The definition of $m$-residue with $m>1$ requires twisting by
$(T_x^\times)^{\otimes (m-1)}$. Alternatively, we may view the
$m$-residue map as a morphism from $\on{Op}_G^{\leq m}(D_x)$ to the
algebraic stack $\g/(G \times \C^\times) \simeq \h/(W \times
\C^\times)$.

Introduce the following affine subspace $\g_{\on{can}}$ of $\g$:
\begin{equation}    \label{gcan}
\g_{\on{can}} = \left\{ p_{-1} + \sum_{j=1}^\ell y_j p_j, \quad y \in
\C \right\}.
\end{equation}
Recall from \cite{Ko} that the adjoint orbit of any regular element in
the Lie algebra $\g$ contains a unique element which belongs to
$\g_{\on{can}}$. Thus, the corresponding morphism $\g_{\on{can}} \to
\g/G = \h/W$ is an isomorphism.

\begin{prop}[\cite{DS}, Prop. 3.8.9]    \label{can real RS}
The natural morphism $\on{Op}_G^{\leq m}(D_x) \to \on{Op}_G(D_x^\times)$
is injective. Its image consists of those $G$-opers on $D_x^\times$
whose canonical representatives have the form
\begin{equation}    \label{can form2}
\nabla = \pa_t + p_{-1} + \sum_{j=1}^\ell t^{-m(d_j+1)} u_j(t) p_j,
\qquad u_j(t) \in \C[[t]].
\end{equation}
Moreover,
\begin{equation}    \label{res2}
\on{res}_m(\nabla) = p_{-1} + \left( u_1(0) + \frac{1}{4} \delta_{m,1}
  \right) p_1 + \sum_{j>1} u_j(0) p_j,
\end{equation}
which is an element of $\g_{\on{can}} = \g/G$ (here we use the
trivialization of $T_x$ induced by the coordinate $t$ in the
definition of $\on{res}_m$).
\end{prop}

\begin{proof}
First, we bring an oper \eqref{oper with RS1} to the form
$$
\pa_t + \frac{1}{t^m} \left( p_{-1} + \sum_{j=1}^\ell c_j(t) p_j
\right), \qquad c_j(t) \in \C[[t]]
$$
(in the same way as in the proof of \lemref{free}). Next, we apply the
gauge transformation by $\crho(t)^{-m}$ and obtain
$$
\pa_t + p_{-1} + m\crho t^{-1}
+ \sum_{j=1}^\ell t^{-m(d_j+1)} c_j(t) p_j, \qquad c_j(t)
\in \C[[t]].
$$
Finally, applying the gauge transformation by $\exp(-mp_1/2t)$, we
obtain the operator
\begin{equation}    \label{with res rho}
\pa_t + p_{-1} + \left( t^{-m-1} c_1(t) - \frac{m^2-2}{4} t^{-2}
\right) p_1 + \sum_{j>1} t^{-m(d_j+1)} c_j(t) p_j, \qquad c_j(t)
\in \C[[t]].
\end{equation}

Thus, we obtain an isomorphism between the space $\on{Op}_G^{\leq
m}(D_x)$ and the space of opers on $D_x^\times$ of the form \eqref{can
form2}. Moreover, comparing formula \eqref{with res rho} with formula
\eqref{can form2} we find that $u_1(t) = c_1(t) - \frac{m^2-2}{4}
t^{m-1}$ and $u_j(t) = c_j(t)$ for $j>1$. Therefore the $m$-residue of
the oper \eqref{can form2} is equal to \eqref{res2}.
\end{proof}

If $m=1$, then the corresponding opers are called opers with {\em
regular singularity}. In $m>1$, then they are called opers with {\em
irregular singularity}. In what follows we will often refer to the
$1$-residue of an oper with regular singularity simply as residue.

Given $\nu \in \h/W$, we denote by $\on{Op}_G^{\leq 1}(D_x)_{\nu}$ the
subvariety of $\on{Op}_G^{\leq 1}(D_x)$ which consists of those opers
that have residue $\nu \in \h/W$.

In particular, the residue of a regular oper $\pa_t + p_{-1} + {\mb
v}(t)$, where ${\mb v}(t) \in \bb(\OO_x)$, is equal to
$\varpi(-\crho)$, where $\varpi$ is the projection $\h \to \h/W$ (see
\cite{BD}). Indeed, a regular oper may be brought to the form
\eqref{oper with RS}, using the gauge transformation with $\crho(t)
\in B(\K_x)$, after which it takes the form
$$
\pa_t + \frac{1}{t} \left( p_{-1} - \crho + t \cdot
\crho(t) ({\mb v}(t)) \crho(t)^{-1} \right).
$$
If ${\mb v}(t)$ is regular, then so is $\crho(t) ({\mb v}(t))
\crho(t)^{-1}$. Therefore the residue of this oper in $\h/W$ is
equal to $\varpi(-\crho)$, and so $\on{Op}_G(D_x) =
\on{Op}_G^{\leq 1}(D_x)_{\varpi(-\crho)}$.

Next, we consider opers with irregular singularities in more
detail. Let $m > 1$. Denote by $\pi$ the projection $\g \to \g/G
\simeq \h/W$. Any point in $\g/G$ may be represented uniquely in the
form $\pi(y)$, where $y$ is a regular element of $\g$ of the form
\begin{equation}    \label{olchi}
y = p_{-1} + \ol{y} = p_{-1} + \sum_{j=1}^\ell y_j p_j.
\end{equation}
Let $N^{(1)}[[t]]$ be the first congruence subgroup of $N[[t]]$, which
consists of all elements of $N[[t]]$ congruent to the identity modulo
$t$. The next lemma follows from the definition.

\begin{lem}    \label{order two}
Let $y$ be a regular element of $\g$ of the form \eqref{olchi}. Then,
for $m > 1$, the space $\on{Op}^{\leq m}_G(D_x)_{\pi(y)}$ of $G$-opers
with singularity of order $m$ and $m$-residue $\pi(y)$ is isomorphic
to the quotient of the space of operators of the form
$$
\pa_t + \frac{1}{t^m}(p_{-1} + \ol{y}) + {\mb v}(t), \qquad {\mb
  v}(t) \in t^{-1} \bb[[t]],
$$
by  $N^{(1)}[[t]]$. Equivalently, it is isomorphic to the space of
operators of the form
$$
\nabla = \pa_t + p_{-1} + \sum_{j=1}^\ell (y_j t^{-m(d_j+1)} +
t^{-m(d_j+1)+1} w_j(t)) p_j, \qquad w_j(t) \in \C[[t]].
$$
\end{lem}

\subsection{Opers without monodromy}

Suppose that $\cla$ is a dominant integral coweight of $\g$. Let
$\on{Op}_G(D_x)_{\cla}$ be the quotient of the space of operators of
the form
\begin{equation} \label{psi la}
\nabla = \pa_t + \sum_{i=1}^\ell \psi_i(t) f_i + {\mb v}(t),
\end{equation}
where $$\psi_i(t) = t^{\langle \al_i,\cla \rangle}(\kappa_i +
t(\ldots)) \in \OO_x, \qquad \kappa_i \neq 0,$$ and ${\mb v}(t) \in
\bb(\OO_x)$, by the gauge action of $B(\OO_x)$. Equivalently,
$\on{Op}_G(D_x)_{\cla}$ is the quotient of the space of operators of
the form
\begin{equation} \label{psi la1}
\nabla = \pa_t + \sum_{i=1}^\ell t^{\langle \al_i,\cla \rangle} f_i +
{\mb v}(t),
\end{equation}
where ${\mb v}(t) \in \bb(\OO_x)$, by the gauge action of
$N(\OO_x)$. Considering the $N(\K_x)$-class of such an operator, we
obtain an oper on $D_x^\times$. Thus, we have a map
$\on{Op}_G(D_x)_{\cla} \to \on{Op}_G(D_x^\times)$.

To understand better the image of $\on{Op}_G(D_x)_{\cla}$ in
$\on{Op}_{^L G}(D_x^\times)$, we introduce, following \cite{FG:local},
Sect. 2.9, a larger space $\on{Op}_G(D_x)_{\cla}^{\on{nilp}}$ as the
quotient of the space of operators of the form \eqref{psi la1}, where
now
$$
{\mb v}(t) \in \h(\OO_x) \oplus t^{-1}\n(\OO_x),
$$
by the gauge action of $N(\OO_x)$.

Consider an operator of the form \eqref{psi la1} with ${\mb v}(t) \in
\h(\OO_x) \oplus t^{-1}\n(\OO_x)$. Denote by
$$
{\mb v}_{-1} = \sum_{\al \in \De_+} {\mb v}_{\al,-1} e_\al \in \n
$$
the coefficient of ${\mb v}(t)$ in front of $t^{-1}$. Then, according
to the definition, $\on{Op}_G(D_x)_{\cla}$ is the subvariety of
$\on{Op}_G(D_x)^{\on{nilp}}_{\cla}$ defined by the equations ${\mb
v}_{\al,-1} = 0, \al \in \De_+$.

It is clear that the monodromy conjugacy class of an oper of the form
\eqref{psi la1} in $\on{Op}_G(D_x)^{\on{nilp}}_{\cla}$ is equal to
$\exp(2 \pi i {\mb v}_{-1})$. Therefore $\on{Op}_G(D_x)_{\cla}$ is the
locus of monodromy-free opers in $\on{Op}_G(D_x)^{\on{nilp}}_{\cla}$.

We have the following alternative descriptions of
$\on{Op}_G(D_x)_{\cla}$ and $\on{Op}_G(D_x)^{\on{nilp}}_{\cla}$.

\begin{prop}[\cite{FG:local},\cite{F:opers}] \label{no mon}
For any dominant integral coweight $\cla$ of $G$ the map
$\on{Op}_G(D_x)^{\on{nilp}}_{\cla} \to \on{Op}_G(D_x^\times)$ is
injective and its image is equal to $\on{Op}_G^{\leq
1}(D_x)_{\varpi(-\cla-\crho)}$. Moreover, the points of
$\on{Op}_G(D_x)_{\cla} \subset \on{Op}_G(D_x)^{\on{nilp}}_{\cla}$ are
precisely those $G$-opers with regular singularity and residue
$\varpi(-\cla-\crho)$ which have no monodromy around $x$.
\end{prop}

In particular, we find that $\on{Op}_G(D_x)_{\cla}$ is a subvariety in
$\on{Op}_G^{\leq 1}(D_x)_{\varpi(-\cla-\crho)}$ defined by
$|\Delta_+|$ equations corresponding to ${\mb v}_{\al,-1} = 0, \al \in
\De_+$.

One may rewrite the elements ${\mb v}_{\al,-1} = 0, \al \in \De_+$,
which are the generators of the defining ideal of
$\on{Op}_G(D_x)_{\cla}$ in $\on{Fun} \on{Op}_G^{\leq
1}(D_x)_{\varpi(-\cla-\crho)}$, as polynomials in the canonical
coordinates $u_{j,n}, j=1,\ldots,\ell; n > 0$, on $\on{Op}_G^{\leq
1}(D_x)_{\varpi(-\cla-\crho)}$ obtained via \propref{can real RS}
(here $u_{j,n}$ is the $t^{n}$-coefficient of $u_j(t)$ in \eqref{can
form2}). It is easy to see that these generators are homogeneous of
degrees $\langle \al,\cla+\crho \rangle, \al \in \De_+$, with respect
to the grading on $\on{Fun} \on{Op}_{^L G}^{\leq
1}(D_x)_{\varpi(-\cla-\crho)}$ defined by the assignment $\deg u_{j,n}
= n$. For instance, for $\cla=0$ these generators are $u_{j,n_j}, j
= 1,\ldots,\ell; n_j = 1,\ldots,d_j$. Other examples are discussed in
\cite{F:icmp}, Sect. 3.9.

\section{Spectra of generalized Gaudin algebras and opers with
  irregular singularities}    \label{spectrum}

In this section we describe the algebra of endomorphisms of
$\V_{0,\ka_c}$, the center of the completed enveloping algebra
$U_{\ka_c}(\ghat)$ and the action of the center on various
$\ghat$-modules of critical level in terms of $^L G$-opers. We then
derive from this description and some general results on coinvariants
\cite{FB,F:faro} that the spectrum of the universal Gaudin algebra
$\ZZ_{(z_i),\infty}(\g)$ is identified with the space of all $^L
G$-opers on $\pone$ with singularities of arbitrary orders at
$z_1,\ldots,z_N$ and $\infty$. The spectrum of the quotient
$\ZZ_{(z_i),\infty}^{(m_i),m_\infty}(\g)$ of $\ZZ_{(z_i),\infty}(\g)$
is identified with the subspace of those $^L G$-opers which have
singularities of orders at most $m_i$ at $z_i$ and $m_\infty$ at
$\infty$. We also describe the spectrum of the algebra ${\mc A}_\chi$
for regular semi-simple and regular nilpotent $\chi$, and the joint
eigenvalues of ${\mc A}_\chi$, and its multi-point generalizations, on
tensor products of finite-dimensional $\g$-modules.

\subsection{The algebra of endomorphisms of $\V_{0,\ka}$}
\label{alg end}

Let $\g$ be a simple Lie algebra. Recall that, by definition, the
Langlands dual Lie algebra $^L \g$ is the Lie algebra whose Cartan
matrix is the transpose to the Cartan matrix of $\g$. In what follows
we will choose Cartan decompositions of $\g$ and $^L \g$ and use the
canonical identification
$$
\h^* = {}^L \h
$$
between $\h^*$ and the Cartan subalgebra $^L\h$ of the Langlands dual
Lie algebra $^L\g$. In particular, we will identify the weights and
roots of $\g$ with the coweights and coroots of $^L \g$, respectively,
and vice versa.

Recall that we denote by $\zz_\ka(\ghat)$ the algebra of endomorphisms
of the vacuum module $\V_{0,\ka}$ of level $\ka$ (see
\secref{action}). We will now describe $\zz_\ka(\G)$ following
\cite{FF:gd,F:wak} (for a detailed exposition, see \cite{F:book}).

Let $^L G$ be the adjoint group of $^L \g$ and $\on{Op}_{^L G}(D)$ the
space of $^L G$-opers on the disc $D = \on{Spec} \C[[t]]$ (see
\secref{defn of opers}).  Denote by $\on{Fun} \on{Op}_{^L G}(D)$ the
algebra of regular functions on $\on{Op}_{^L G}(D)$. In view of
\lemref{free}, it is isomorphic to the algebra of functions on the
space of $\ell$-tuples $(v_1(t),\ldots,v_\ell(t))$ of formal Taylor
series, i.e., the space $\C[[t]]^\ell$. If we write
$$
v_i(t) = \sum_{n<0} v_{i,n} t^{-n-1},
$$
then we obtain
\begin{equation} \label{descr of opers}
\on{Fun} \on{Op}_{^L G}(D) \simeq \C[v_{i,n}]_{i=1,\ldots,\ell; n<0}.
\end{equation}

Let $\DerO = \C[[t]] \pa_t$ be the Lie algebra of continuous
derivations of the topological algebra $\OO = \C[[t]]$. The action of
its Lie subalgebra $\on{Der}_0 \OO = t \C[[t]] \pa_t$ on $\OO$
exponentiates to an action of the group $\AutO$ of formal changes of
variables. Both $\DerO$ and $\AutO$ naturally act on $\V_{0,\ka}$ in a
compatible way, and these actions preserve $\zz_\ka(\G)$. They also
act on the space $\on{Op}_{^L G}(D)$. One can check that the vector
field $-t \pa_t$ defines a $\Z_+$-grading on $\on{Fun} \on{Op}_{^L
G}(D)$ such that $\deg v_{i,n} = d_i-n$, and the vector field $-\pa_t$
acts as a derivation such that $-\pa_t \cdot v_{i,n} = -n v_{i,n-1}$.

Recall the critical invariant inner product $\ka_c$ introduced in
\secref{universal}.

\begin{thm}[\cite{FF:gd,F:wak}] \label{center} \

{\em (1)} $\zz_\ka(\G) = \C$ if $\ka \neq \ka_c$.

{\em (2)} There is a canonical isomorphism $$\zz_{\ka_c}(\G) \simeq
\on{Fun} \on{Op}_{^L G}(D)$$ of algebras which is compatible with the
actions of $\DerO$ and $\AutO$.
\end{thm}

Since $\zz_\ka(\G)$ is trivial for $\ka\neq \ka_c$, we will set
$\ka=\ka_c$ and omit $\ka$ from our notation.

\medskip

Let again $x$ be a smooth point of a curve $x$ and $\OO_x \subset
\K_x$ be as in \secref{can repr}. We have the Lie algebra $\g(\K_x)$
and its central extension $\ghat_{\ka_c,x}$ defined by the commutation
relations \eqref{comm rel}. Since the residue is
coordinate-independent, the Lie algebra $\ghat_{\ka_c,x}$ is
independent on the choice of an isomorphism $\K_x \simeq \C\ppart$. We
have a Lie subalgebra $\g(\OO_x) \subset \ghat_{\ka_c,x}$ and we
define the corresponding vacuum module of level $\ka_c$ as
$$
\V_{0,x} = \on{Ind}_{\g(\OO_x) \oplus \C {\mb 1}}^{\ghat_{\ka_c,x}}
\C.
$$
Set
$$
\zz(\ghat)_x = (\V_{0,x})^{\g(\OO_x)} = \on{End}_{\ghat_{\ka_c,x}}
\V_{0,x}.
$$
Then the compatibility of the isomorphism of \thmref{center} with the
action of $\AutO$ implies the existence of the following canonical
(i.e., coordinate-independent) isomorphism:
\begin{equation} \label{isom x}
\zz(\G)_x \simeq \on{Fun} \on{Op}_{^L G}(D_x).
\end{equation}

Recall from \secref{grad} that the action of $L_0 = -t\pa_t \in \DerO$
on the module $\V_0$ defines a $\Z$-grading on it such that $\deg v_0
= 0, \deg J^a_n = -n$. The action of $L_{-1} = -\pa_t \in \DerO$ is
given by the translation operator $T$ defined by formula
\eqref{translation}. \thmref{center} and the isomorphism \eqref{descr
of opers} imply that there exist non-zero vectors $S_i \in
\V_0^{\g[[t]]}$ of degrees $d_i+1, i = 1,\ldots,\ell$, such that
$$
\zz(\G) = \C[S^{(n)}_i]_{i = 1,\ldots,\ell; n\geq 0} v_0,
$$
where $S_i^{(n)} = T^n S_i$, and under the isomorphism of
\thmref{center} we have
\begin{equation}    \label{S and v1}
S^{(n)}_i \mapsto n! v_{i,-n-1}, \quad n \geq 0.
\end{equation}
The $\Z$-gradings on both algebras get identified and the
action of $T$ on $\zz(\G)$ becomes the action of $-\pa_t$ on $\on{Fun}
\on{Op}_{^L G}(D)$.

Now recall from \secref{grad} that the PBW filtration on
$U(\ghat_{\ka_c})$ induces a filtration on $\V_0$ such that the
associated graded is identified with
$$
S(\g\ppart/\g[[t]]) = \on{Fun} (\g\ppart/\g[[t]])^* =
\on{Fun} (\g^*[[t]]dt) \simeq \on{Fun} \g^*[[t]],
$$
where we use the canonical non-degenerate pairing
$$
\phi(t)dt \in \g^*[[t]]dt,A(t) \in \g\ppart/\g[[t]] \mapsto
\on{Res}_{t=0} \langle \phi(t),A(t) \rangle dt,
$$
and a coordinate $t$ on $D$. Let
$$
\on{Inv} \g^*[[t]] = (\on{Fun} \g^*[[t]])^{\g[[t]]}
$$
be the algebra of $\g[[t]]$-invariant functions on
$\g^*[[t]]$. According to \thmref{thm:FF}, the map
$$
\on{gr} \zz(\ghat) \to \on{gr} \V_0 = \on{Fun} \g^*[[t]]
$$
gives rise to an isomorphism
$$
\on{gr} \zz(\ghat) \simeq \on{Inv} \g^*[[t]].
$$
In particular, the symbols of the generators $S^{(n)}_i$ of
$\zz(\ghat)$ have the following simple description.

Let $\on{Inv} \g^*$ be the algebra of $\g$-invariant functions on
$\g^*$. By \thmref{chevalley},
$$
\on{Inv} \g^* = \C[\ol{P}_i]_{i=1,\ldots\ell},
$$
where the generators $\ol{P}_i$ may be chosen in such a way that
they are homogeneous and $\on{deg} \ol{P}_i = d_i + 1$.

As in \secref{grad}, we use the elements $\ol{P}_i$ to construct
generators of the algebra $\on{Inv} \g^*[[t]]$. We will use the
generators $\ol{J}^a_n, n<0$, of $S (\g\ppart/\g[[t]]) = \on{Fun}
\g^*[[t]]$, which are the symbols of $J^a_n \vac \in \V_0$. These are
linear functions on $\g^*[[t]]$ defined by the formula
\begin{equation}    \label{ol J}
\ol{J}^a_n(\phi(t)) = \on{Res}_{t=0} \langle \phi(t),J^a \rangle t^n
dt.
\end{equation}
We will also write
$$\ol J^a(z) = \sum_{n<0} \ol J^a_n z^{-n-1}.$$

Let us write $\ol{P}_i$ as a polynomial in the $\ol{J}^a$'s, $\ol{P}_i
= \ol{P}_i(\ol J^a)$. Define a set of elements $\ol{P}_{i,n} \in
\on{Fun} \g^*[[t]]$ by the formula
\begin{equation}    \label{ol P}
\ol{P}_i(\ol J^a(z)) = \sum_{n<0} \ol{P}_{i,n} z^{-n-1}.
\end{equation}
Note that each of the elements $\ol{P}_{i,n}$ is a finite polynomial
in the $\ol{J}^a_n$'s. Now Theorems \ref{thm:BD} and \ref{thm:FF}
imply the following:

\begin{lem}    \label{free pol}
The generators $\ol{P}_i$ of $\on{Inv} \g^*$ may be chosen
in such a way that the symbol of $S_i^{(n)} \in \zz(\ghat)$ is equal
to $n!  \ol{P}_{i,-n-1}$.
\end{lem}

For example, we may choose the degree $2$ vector $S_1$ to be the
Segal-Sugawara vector \eqref{sugawara} (it is unique up to a non-zero
scalar). Its symbol is equal to
$$
\ol{P}_{1,-1} = \frac{1}{2} \sum_a \ol{J}_{a,-1} \ol{J}^a_{-1},
$$
where $\ol{P}_1 = \frac{1}{2} \sum_a J_a J^a$ is the quadratic Casimir
generator of $\on{Inv} \g^*$.

The algebra $\on{Fun} \on{Op}_{^L G}(D)$ has a canonical filtration
such that the associated graded algebra is isomorphic to
$$
\on{Inv} {}^L\g[[t]] = (\on{Fun} {}^L \g[[t]])^{^L G[[t]]}
$$
(see \cite{F:wak}, Sect. 11.3). The spectrum of this algebra is the
jet scheme of $$^L \g/{}^L G = \on{Spec} (\on{Fun} {}^L \g)^{^L G}.$$
Using the canonical isomorphisms
$$
^L \g/{}^L G = {}^L \h/W = \h^*/W = \g^*/G,
$$
we identify the jet scheme of $^L \g/{}^L G$ with the jet scheme of
$\g^*/G$. Therefore we have a canonical isomorphism
\begin{equation}    \label{isom of gr}
\on{Inv} {}^L\g[[t]] = \on{Inv} \g^*[[t]].
\end{equation}
On the other hand, we know that
$$
\on{gr} \zz(\ghat) = \on{Inv} \g^*[[t]].
$$

The following result is proved in \cite{F:wak}, Theorem 11.4.

\begin{prop}    \label{sign1}
The isomorphism $\zz(\ghat) \simeq \on{Fun} \on{Op}_{^L G}(D)$ of
\thmref{center},(2) preserves the filtrations on both algebras, and
the corresponding isomorphism of the associated graded algebras is the
isomorphism \eqref{isom of gr} multiplied by $(-1)^n$ on the subspaces
of degree $n$.
\end{prop}

\subsection{The center of the completed enveloping algebra}

Recall from \cite{FB} that $\V_0 = \V_{0,\ka_c}$ is a vertex algebra,
and $\zz(\G)$ is its commutative vertex subalgebra; in fact, it is the
center of $\V_0$. We will also need the center of the completed
universal enveloping algebra of $\ghat$ of critical level. This
algebra is defined as follows.

Let $U_{\ka_c}(\ghat)$ be the quotient of the universal enveloping
algebra $U(\ghat_{\ka_c})$ of $\ghat_{\ka_c}$ by the ideal generated
by $({\mb 1}-1)$. Define its completion $\wt{U}_{\ka_c}(\ghat)$ as
follows:
$$
\wt{U}_{\ka_c}(\ghat) = \underset{\longleftarrow}\lim \;
U_{\ka_c}(\ghat)/U_{\ka_c}(\ghat) \cdot (\g \otimes t^N\C[[t]]).
$$
It is clear that $\wt{U}_{\ka_c}(\ghat)$ is a topological algebra
which acts on all smooth $\ghat_{\ka_c}$-module, i.e. such that
any vector is annihilated by $\g \otimes t^N\C[[t]]$ for sufficiently
large $N$, and the central ${\mb 1}$ acts as the identity. Let
$Z(\ghat)$ be the center of $\wt{U}_{\ka_c}(\ghat)$.

Denote by $\on{Fun} \on{Op}_{^L G}(D^\times)$ the algebra of regular
functions on the space $\on{Op}_{^L G}(D^\times)$ of $^L G$-opers on
the punctured disc $D^\times = \on{Spec} \C\ppart$. In view of
\lemref{free}, it is isomorphic to the algebra of functions on the
space of $\ell$-tuples $(v_1(t),\ldots,v_\ell(t))$ of formal Laurent
series, i.e., the ind-affine space $\C\ppart^\ell$. If we write $v_i(t)
= \sum_{n \in \Z} v_{i,n} t^{-n-1}$, then we obtain that $\on{Fun}
\on{Op}_{^L G}(D^\times)$ is isomorphic to the completion of the
polynomial algebra $\C[v_{i,n}]_{i = 1,\ldots,\ell; n \in \Z}$ with
respect to the topology in which the basis of open neighborhoods of
zero is formed by the ideals generated by $v_{i,n}, i = 1,\ldots,\ell;
n \geq N$, for $N \geq 0$.

\begin{thm}[\cite{FF:gd,F:wak}] \label{center1}
There is a canonical isomorphism $$Z(\G) \simeq \on{Fun} \on{Op}_{^L
G}(D^\times)$$ of complete topological algebras which is compatible
with the actions of $\DerO$ and $\AutO$.
\end{thm}

For a smooth point $x \in X$ as above, we have the Lie algebra
$\ghat_{\ka_c,x}$. We define the completed enveloping algebra
$\wt{U}_{\ka_c}(\ghat_x)$ of $\ghat_{\ka_c,x}$ in the same way as
above. Let $\Z(\ghat)_x$ be its center. Then \thmref{center1} implies
the following isomorphism:
$$
Z(\G)_x \simeq \on{Fun} \on{Op}_{^L G}(D_x^\times).
$$

Each element $A \in \V_0$ gives rise to a ``vertex operator'' which is
a formal power series
$$
Y[A,z] = \sum_{n \in \Z} A_{[n]} z^{-n-1}, \qquad A_{[n]} \in
\wt{U}_{\ka_c}(\ghat)
$$
(see \cite{FB}, Sect. 4.2). In particular, we have the elements
$S_{i,[n]}$ attached to the generators $S_i \in \zz(\ghat) \subset
\V_0$. Under the isomorphism of \thmref{center1} we have
\begin{equation}    \label{S and v}
S_{i,[n]} \to v_{i,n}.
\end{equation}

The algebra $\wt{U}_{\ka_c}(\ghat)$ has a PBW filtration, and its
associated graded algebra is the completed symmetric algebra
$\wt{S}(\g\ppart)$ of $\g\ppart$, which we identify with the
topological algebra $\on{Fun} \g^*\ppart$. Let $\ol{P}_{i,n}$ be the
symbol of the element $S_{i,[n]}$ in $\on{Fun} \g^*\ppart$. These
elements are given by the formula
\begin{equation}    \label{ol P1}
\ol{P}_i(\ol{J}^a(z)) = \sum_{n \in \Z} \ol{P}_{i,n} z^{-n-1},
\end{equation}
where
$$
\ol{J}^a(z) = \sum_{n \in \Z} \ol{J}^a_n z^{-n-1},
$$
and the $\ol{J}^a_n$'s are the generators of $\on{Fun} \g^*\ppart$,
defined by the formula
$$
\ol{J}^a_n(\phi(t)) = \on{Res}_{t=0} \langle \phi(t),J^a \rangle t^n
dt, \qquad \phi(t) \in \g^*\ppart, n \in \Z,
$$
(compare with formula \eqref{ol J}).

It is easy to see that the elements $\ol{P}_{i,n}$ are
$G\ppart$-invariant elements of $\on{Fun} \g^*\ppart$. Moreover, they
are topological generators of the algebra $(\on{Fun}
\g^*\ppart)^{G\ppart}$. More precisely, $(\on{Fun}
\g^*\ppart)^{G\ppart}$ is isomorphic to a completion of the free
polynomial algebra in $\ol{P}_{i,n}, i=1,\ldots, \ell; n \in \Z$, see
\cite{BD}, Theorem 3.7.5.

On the other hand, the algebra $\on{Fun} \on{Op}_{^L G}(D^\times)$
also has a canonical filtration such that the associated graded
algebra is canonically isomorphic to $(\on{Fun} {}^L \g\ppart)^{^L
G\ppart}$. In the same way as at the end of \secref{spectrum} we
obtain a canonical isomorphism
$$
(\on{Fun} {}^L \g\ppart)^{^L G\ppart} \simeq (\on{Fun}
\g^*\ppart)^{G\ppart}.
$$
We can now take the symbol of $v_{i,n}, n \in \Z$, in $(\on{Fun}
\g^*\ppart)^{G\ppart}$ and compare it with $\ol{P}_{i,n} =
\sigma(S_{i,[n]})$, where $v_{i,n}$ and $S_{i,[n]}$ are related by
formula \eqref{S and v}. Using the commutative vertex algebra
structures on $\zz(\ghat)$ and $\on{Fun} \on{Op}_{^L G}(D)$ and
\propref{sign1}, we obtain the following:

\begin{lem}    \label{sign}
The symbol of $v_{i,n}$ is equal to $(-1)^{d_i+1} \ol{P}_{i,n}$.
\end{lem}

If $M$ is a smooth $\ghat_{\ka_c}$-module, then the action of
$Z(\ghat)$ on $M$ gives rise to a homomorphism
$$
Z(\ghat) \to \on{End}_{\ghat_{\ka_c}} M.
$$
For example, if $M=\V_0$, then using Theorems \ref{center} and
\ref{center1} we identify this homomorphism with the surjection
$$
\on{Fun} \on{Op}_{^L G}(D^\times) \twoheadrightarrow \on{Fun}
\on{Op}_{^L G}(D)
$$
induced by the natural embedding
$$
\on{Op}_{^L G}(D) \hookrightarrow \on{Op}_{^L G}(D^\times).
$$

Recall that the Harish-Chandra homomorphism identifies the center
$Z(\g)$ of $U(\g)$ with the algebra $(\on{Fun} \h^*)^W$ of polynomials
on $\h^*$ which are invariant with respect to the action of the Weyl
group $W$. Therefore a character $Z(\g) \to \C$ is the same as a
point in $\on{Spec} (\on{Fun} \h^*)^W$ which is the quotient $\h^*/W$.
For $\la \in \h^*$ we denote by $\varpi(\la)$ its projection onto
$\h^*/W$. In particular, $Z(\g)$ acts on $M_\la$ and $V_\la$ via its
character $\varpi(\la+\rho)$. We also denote by $I_\la$ the maximal
ideal of $Z(\g)$ equal to the kernel of the homomorphism $Z(\g) \to
\C$ corresponding to the character $\varpi(\la+\rho)$.

Recall that we have an isomorphism
$$
\g^*/G \simeq \h^*/W = {}^L \h/W.
$$
Let us denote by $\pi$ the corresponding map $\g^* \to {}^L \h/W$.

\begin{thm}\label{factor RS}\ %

{\em (1)} Let
$$
\UU_m = \on{Ind}_{t^m\g[[t]] \oplus \C {\mb 1}}^{\ghat_{\ka_c}} \C.
$$
The homomorphism $Z(\ghat) \to \on{End}_{\ghat_{\ka_c}} \UU_m$ factors
as
$$
Z(\ghat) \simeq \on{Fun} \on{Op}_{^L G}(D^\times)
\twoheadrightarrow \on{Fun} \on{Op}_{^L G}^{\leq m}(D) \hookrightarrow
\on{End}_{\ghat_{\ka_c}} \UU_m.
$$

{\em (2)} Let
$$
\I_{m,\chi} = \on{Ind}_{t^{m}\g[[t]] \oplus \C {\mb
    1}}^{\ghat_{\ka_c}} \C_\chi,
$$
where $\chi \in \g^*$ and $t^{m}\g[[t]]$ acts
on $\C_\chi$ via
$$
t^{m}\g[[t]] \to \g \otimes t^{m} \overset{\chi}\longrightarrow
\C.
$$
The homomorphism $Z(\ghat) \to \on{End}_{\ghat_{\ka_c}}
\I_{m,\chi}$ factors as
\begin{equation}    \label{inj for Ichi}
Z(\ghat) \simeq \on{Fun} \on{Op}_{^L G}(D^\times)
\twoheadrightarrow \on{Fun} \on{Op}_{^L G}^{\leq
  (m+1)}(D)_{\pi(-\chi)} \rightarrow \on{End}_{\ghat_{\ka_c}}
\I_{m,\chi},
\end{equation}
where $\on{Op}_{^L G}^{\leq (m+1)}(D)_{\pi(-\chi)}$ is the space of
opers with singularity of order $m+1$ and the $(m+1)$-residue
$\pi(-\chi)$.

In addition, if $\chi$ is a regular element of $\g^*$, then the last
map in \eqref{inj for Ichi} is injective.

{\em (3)}
Let $M$ be a $\g$-module on which the center $Z(\g)$ acts
via its character $\varpi(\la+\rho)$, and let $\M$ be the induced
$\ghat_{\ka_c}$-module. Then the homomorphism $Z(\ghat) \to
\on{End}_{\ghat_{\ka_c}} \M$ factors as
$$
Z(\ghat) \simeq \on{Fun} \on{Op}_{^L G}(D^\times) \twoheadrightarrow
\on{Fun} \on{Op}_{^L G}^{\leq 1}(D)_{\varpi(-\la-\rho)} \to
\on{End}_{\ghat_{\ka_c}} \M.
$$

{\em (4)}
For an integral dominant weight $\la \in \h^*$ the
homomorphism
$$
\on{Fun} \on{Op}_{^L G}(D^\times) \to
\on{End}_{\ghat_{\ka_c}} \V_\la
$$
factors as
$$
\on{Fun} \on{Op}_{^L G}(D^\times) \to \on{Fun} \on{Op}_{^L
G}(D)_\la \to \on{End}_{\G_{\ka_c}} \V_\la.
$$
\end{thm}

\begin{proof}

According to \propref{can real RS}, the ideal of $\on{Op}_{^L G}^{\leq
m}(D)$ in $\on{Fun} \on{Op}_{^L G}(D^\times)$ is generated by
$v_{i,n_i}, i=1,\ldots,\ell; n_i \geq m(d_i+1)$. By formula \eqref{S
and v}, these correspond to $S_{i,[n_i]}, i=1,\ldots,\ell; n_i \geq
m(d_i+1)$. It follows from the definition of vertex operators that for
any $A \in \V_0$ of degree $N$, the operators $A_{[n]}, n\geq mN$, act
by zero on any vector that is annihilated by $\g \otimes t^m
\C[[t]]$. Since $\deg S_i = d_i+1$, we obtain that the homomorphism
$Z(\ghat) \to \on{End}_{\ghat_{\ka_c}} \UU_m$ factors as
$$
Z(\ghat) \simeq \on{Fun} \on{Op}_{^L G}(D^\times)
\twoheadrightarrow \on{Fun} \on{Op}_{^L G}^{\leq m}(D) \to
\on{End}_{\ghat_{\ka_c}} \UU_m.
$$
To complete the proof of part (1), we need to show that the last
homomorphism is injective. It suffices to prove that the
map
$$
\on{Fun} \on{Op}_{^L G}^{\leq m}(D) \to \UU_m,
$$
applied by acting on the generating vector of $\UU_m$, is
injective. This map preserves natural filtrations on both spaces,
and it is sufficient to show that the corresponding map of the
associated graded is injective.

The PBW filtration on $U(\ghat_{\ka_c})$ induces a filtration on
$\UU_m$, and the associated graded space is identified with
the symmetric algebra $S(\g\ppart/t^m\g[[t]])$. On the other hand, we
have, by \propref{can real RS},
$$
\on{Fun} \on{Op}_{^L G}^{\leq m}(D) \simeq
\C[v_{i,n_i}]_{i=1,\ldots,\ell;n<m(d_i+1)}.
$$
Now \lemref{sign} implies that
$$
\on{gr} \on{Fun} \on{Op}_{^L G}^{\leq m}(D) \simeq
\C[\ol{P}_{i,n_i}]_{i=1,\ldots,\ell;n<m(d_i+1)}
$$
(see formula \eqref{ol P1}). Thus, we need to show that the map
$$
\C[\ol{P}_{i,n_i}]_{i=1,\ldots,\ell;n<m(d_i+1)} \to
S(\g\ppart/t^m\g[[t]]),
$$
is injective. Let us apply the automorphism $\ol{J}^a_n
\mapsto \ol{J}_{n-m}$ of $\g\ppart$ (considered as a vector space) to
both sides. Then we have $\ol{P}_{i,n} \mapsto \ol{P}_{i,n-m(d_i+1)}$,
and so the above map becomes
$$
\C[\ol{P}_{i,n_i}]_{i=1,\ldots,\ell;n<0} \to S(\g\ppart/\g[[t]]),
$$
which is injective by \thmref{thm:BD}. This completes the proof of
part (1).

Next, we prove part (2). We will work with $\I_{m-1,\chi}$ instead of
$\I_{m,\chi}$. According to the result of part (1), we have an
injective homomorphism from the algebra
$$
\on{Fun} \on{Op}_{^L G}^{\leq m}(D^\times) \simeq
\C[v_{i,n_i}]_{i=1,\ldots,\ell; n_i < m(d_i+1)},
$$
to $\on{End}_{\ghat_{\ka_c}} \UU_m$, sending
$$
v_{i,n_i} \mapsto S_{i,[n_i]}.
$$
In addition, the commutative algebra
$$
S(\g) = S(\g \otimes t^{m-1}) = \C[\ol{J}^a_{m-1}]_{a=1\ldots,\dim \g}
$$
also acts on $\UU_m$ and commutes with the action $\on{Fun}
\on{Op}_{^L G}^{\leq m}(D^\times)$. Let us denote the generating
vector of $\UU_m$ by $v_m$. Since the elements $J^a_n, n \geq m$, of
$\ghat_{\ka_c}$ annihilate $v_m$, we obtain from \lemref{lexico} that
the generator $S_{i,[m(d_i+1)-1]}$ acts on the generating vector of
$\UU_m$ by multiplication by $\ol{P}_i(\ol{J}^a_{m-1})$.

It follows from the definition that $\I_{m-1,\chi}$ is the quotient of
$\UU_m$ by the maximal ideal in $S(\g) = \on{Fun} \g^*$ corresponding
to $\chi \in \g^*$. Therefore we find that $S_{i,[m(d_i+1)-1]} \in
Z(\ghat)$ acts on the generating vector of $\I_{m-1,\chi}$, and hence
on the entire module $\I_{m-1,\chi}$, by multiplication by the value
of $\ol{P}_i \in \on{Fun} \g^*$ at $\chi \in \g^*$. By \lemref{sign}
and \propref{can real RS}, this means that the action of $\on{Fun}
\on{Op}_{^L G}^{\leq m}(D^\times)$ on $\I_{m-1,\chi}$ factors through
the algebra of functions on opers with singularity of order $m$ and
the $m$-residue $\pi(-\chi)$.

Thus, we obtain that the homomorphism $Z(\ghat)
\to \on{End}_{\ghat_{\ka_c}} \I_{m-1,\chi}$ factors as
$$
Z(\ghat) \simeq \on{Fun} \on{Op}_{^L G}(D^\times)
\twoheadrightarrow \on{Fun} \on{Op}_{^L G}^{\leq m}(D)_{\pi(-\chi)}
\to \on{End}_{\ghat_{\ka_c}} \I_{m-1,\chi}.
$$
To complete the proof of part (2), we need to show that the last
homomorphism is injective if $\chi$ is regular. It suffices to show
that the natural map
\begin{equation}    \label{into Ichi}
\C[S_{i,[n_i]}]_{i=1,\ldots,\ell; n_i < m(d_i+1)-1} \to \I_{m-1,\chi}
\end{equation}
obtained by acting on the generating vector of $\I_{m-1,\chi}$, is
injective. As in the proof of part (1), we will derive this from the
injectivity of the corresponding maps of associated graded spaces.

The associated graded space of $\I_{m-1,\chi}$ with respect to the PBW
filtration on the universal enveloping algebra of $\ghat_{\ka_c}$ is
naturally identified with
$$
S(\g\ppart/t^{m-1}\g[[t]]) \simeq \on{Fun} (t^{-m+1} \g^*[[t]]
dt) \simeq \on{Fun} (t^{-m+1} \g^*[[t]]).
$$
Now, using \lemref{lexico}, we obtain that the symbol of the image of
a polynomial $R$ in $S_{i,[n_i]}$ under the map \eqref{into Ichi},
considered as a function on $t^{-m+1} \g^*[[t]]$, is equal
to the same polynomial $R$ in which we make a replacement
$$
S_{i,[n]} \mapsto \ol{P}_{i,n}(\cdot + \chi t^{-m})
$$
followed by the shift of argument by $\chi t^{-m}$. Thus, injectivity
of the associated graded map of \eqref{into Ichi} is equivalent to the
algebraic independence of the restrictions of $$\ol{P}_{i,n_i}, \qquad
i=1,\ldots,\ell; \; n_i < m(d_i+1)-1,$$ which are functions on
$\g^*\ppart$, to $\chi t^{-m} + t^{-m+1} \g^*[[t]]$. Let us prove this
algebraic independence.

To simplify the argument, let us apply the automorphism of $\g\ppart$
sending $\ol{J}^a_n$ to $\ol{J}^a_{n-m}$. Then $\chi t^{-m} + t^{-m+1}
\g^*[[t]]$ becomes $\chi + t\g^*[[t]]$ and $\ol{P}_{i,n} \mapsto
\ol{P}_{i,n-m(d_i+1)}$. We therefore need to prove that the
restrictions of the polynomials $\ol{P}_{i,n}, i=1,\ldots,\ell; n<-1$,
to $\chi + t\g^*[[t]] \subset \g^*[[t]]$ are algebraically
independent. This would follow if we show that their differentials at
the point $\chi$ are linearly independent.

Let us identify the cotangent space to $\chi \in \g^*$ with $(\g^*)^*
= \g$ and the cotangent space to $\chi \in \chi + t\g^*[[t]]$ with
$(t\g^*[[t]])^* = \g\ppart/t^{-1}\g[[t]]$. Since $\chi$ is regular, we
obtain from \cite{Ko} that the values $d\ol{P}_i|_\chi$ of the
differentials $d\ol{P}_i$ of the generators $\ol{P}_i$ of the algebra
$\on{Inv} \g^*$ at $\chi \in \g^*$ are linearly independent vectors in
the centralizer $\g_\chi$ of $\chi$, considered as a subspace of
$\g$. Using the explicit formula \eqref{ol P} for $\ol{P}_{i,n}$, we
find that the value of $d\ol{P}_{i,n}, n<-1$, at $\chi \in \g^*[[t]]$
is equal to
$$
(d\ol{P}_i|_\chi) \otimes t^n \in \g\ppart/t^{-1}\g[[t]].
$$
These vectors are linearly independent for $i=1,\ldots,\ell$ and
$n<-1$. Therefore the restrictions of the polynomials $\ol{P}_{i,n},
i=1,\ldots,\ell; n<-1$, to $\chi + t\g^*[[t]] \subset \g^*[[t]]$ are
algebraically independent, and so the map \eqref{into Ichi} is
injective. This completes the proof of part (2).

Part (3) of the Theorem is obtained by combining Theorem 12.4, Lemma
9.4 and Proposition 12.8 of \cite{F:wak}.

Part (4) is established in \cite{FG:fusion}, Lemma 1.7.
\end{proof}

We note that parts (1) and (2) of the theorem may be interpreted as
saying that the supports of the $Z(\ghat)$-modules $\UU_m$ and
$\I_{m,\chi}$, considered as subvarieties in $\on{Op}_{^L
G}(D^\times)$, are equal to $\on{Op}_{^L G}^{\leq (m+1)}(D)$ and
$\on{Op}_{^L G}^{\leq (m+1)}(D)_{\pi(-\chi)}$, respectively. If $\chi$
were not regular, then the support of $\I_{m,\chi}$ would still be
contained in $\on{Op}_{^L G}^{\leq (m+1)}(D)_{\pi(-\chi)}$, but it
would not be equal to it (equivalently, the last map in \eqref{inj for
Ichi} would not be injective). For example, if $\chi=0$, we have
$\I_{m,0} = \UU_m$, and so the support of $\I_0$ is equal to
$\on{Op}^{\leq m}_{^L G}(D)$. This explains the special role played by
regular characters $\chi$.

\subsection{Description of the Gaudin algebras}

Now we are ready to show that the Gaudin algebras
$\ZZ_{(z_i),\infty}^{(m_i),m_\infty}(\g)$, which were introduced in
\secref{universal}, are isomorphic to algebras of functions on opers
on $\pone \bs \{ z_1,\ldots,z_N,\infty \}$ with appropriate
singularities at the points $z_1,\ldots,z_N$ and $\infty$. Here we
follow the analysis of \cite{F:faro}, Sect. 2.5, where the Gaudin
algebras were described in the case regular singularities (when all
$m_i$ and $m_\infty$ are equal to $1$).

Let us denote by
$$
\on{Op}_{^L G}(\pone)^{(m_i),m_\infty}_{(z_i),\infty}
$$
the space of $^L G$-opers on $\pone \bs \{
z_1,\ldots,z_N,\infty \}$ whose restriction to $D_{z_i}^\times$
belongs to
$$
\on{Op}^{\leq m_i}_{^L G}(D_{z_i}^\times) \subset \on{Op}_{^L
  G}(D_{z_i}^\times), \qquad i=1,\ldots,N,
$$
and whose restriction to $D_{\infty}^\times$ belongs to
$$
\on{Op}^{\leq m_\infty}_{^L G}(D_{\infty}^\times)  \subset \on{Op}_{^L
  G}(D_{\infty}^\times).
$$
Thus, points of $\on{Op}_{^L
  G}(\pone)^{(m_i),m_\infty}_{(z_i),\infty}$ are $^L G$-opers on
$\pone$ with singularities at $z_i, i=1,\ldots,N$, and $\infty$ of
orders $m_i, i=1,\ldots,N$, and $m_\infty$, respectively (and regular
elsewhere).

\begin{thm}    \label{ZZ descr} \

{\em (1)} There is an isomorphism of algebras
$$
\ZZ_{(z_i),\infty}^{(m_i),m_\infty}(\g) \simeq \on{Fun} \on{Op}_{^L
G}(\pone)^{(m_i),m_\infty}_{(z_i),\infty}.
$$

\medskip

Let us fix a $\g_{m_i}$-module $M_i$ for each $i=1,\ldots,N$, and a
$\ol\g_{m_\infty}$-module $M_\infty$. Then the following holds:

\medskip

{\em (2)} Suppose that we have $m_j=1$ and let $M_j$ be a $\g$-module
on which the center $Z(\g)$ acts via its character
$\varpi(\la+\rho)$. Then the action of
$\ZZ_{(z_i),\infty}^{(m_i),m_\infty}(\g)$ on $\bigotimes_{i=1}^N M_i
\otimes M_\infty$ factors through the algebra of functions on a subset
of $\on{Op}_{^L G}(\pone)^{(m_i),m_\infty}_{(z_i),\infty}$, which
consists of the opers with regular singularity and $1$-residue
$\varpi(-\la_j-\rho)$ at $z_j$.

{\em (3)} Under the assumptions of part (2), suppose in addition that
$\la_j$ is an integral dominant weight and $M_j = V_{\la_j}$ is the
finite-dimensional irreducible $\g$-module with highest weight
$\la_j$. Then the action of $\ZZ_{(z_i),\infty}^{(m_i),m_\infty}(\g)$
on $\bigotimes_{i=1}^N M_i \otimes M_\infty$ factors through the algebra
of functions on a subset of $\on{Op}_{^L
G}(\pone)^{(m_i),m_\infty}_{(z_i),\infty}$, which consists of the
opers with regular singularity at $z_j$, $1$-residue
$\varpi(-\la_j-\rho)$, and trivial monodromy around $z_j$.

{\em (4)} Now let $\chi \in \g^*$ be a regular element. Set $M_j =
\I_{m_j,\chi}$. Then the action of
$\ZZ_{(z_i),\infty}^{(m_i),m_\infty}(\g)$ on $\bigotimes_{i=1}^N M_i
\otimes M_\infty$ factors through the algebra of functions on a subset
of $\on{Op}_{^L G}(\pone)^{(m_i),m_\infty}_{(z_i),\infty}$, which
consists of the opers with singularity of order $m_j+1$ at $z_j$ and
$(m_j+1)$-residue $\pi(-\chi)$ (here we trivialize the tangent space
to $z_j$ using the global coordinate $t$ on $\pone$).

An analogous result also holds with $z_j$ replaced by $\infty$.
\end{thm}

\noindent {\em Proof} is a word-for-word repetition of the argument
used in the proof of Theorem 2.7 of \cite{F:faro} (which corresponds
to the special case of parts (1)--(3) of the above theorem when all
the $m_i$'s and $m_\infty$ are equal to $1$). Using the general
results of \cite{FB} about the action of commutative vertex algebras
on coinvariants, we show the following. Let $\M_1,\ldots,\M_N$ and
$\M_\infty$ be $\ghat_{\ka_c}$-modules. Then, as explained in
\secref{universal}, the universal Gaudin algebra
$\ZZ_{(z_i),\infty}(\g)$ acts on the corresponding space
$H(\M_1,\ldots,\M_N,\M_\infty)$ of coinvariants.  Suppose that the
action of $Z(\ghat) \simeq \on{Fun} \on{Op}_{^L G}(D^\times)$ on
$\M_i$ factors through $\on{Fun} \on{Op}^{M_i}_{^L G}(D^\times)$,
where $\on{Op}^{M_i}_{^L G}(D^\times) \subset \on{Op}_{^L
G}(D^\times)$. Then, in the same way as in the proof of Theorem 2.7 of
\cite{F:faro} we obtain that the action of $\ZZ_{(z_i),\infty}(\g)$
factors through the space of $^L G$-opers on $\pone \bs \{
z_1,\ldots,z_N,\infty \}$ whose restriction to the punctured disc
$D_{z_i}^\times, i=1,\ldots,N$ (resp., $D_\infty^\times$) belongs to
$\on{Op}^{M_i}_{^L G}(D_{z_i}^\times)$ (resp., $\on{Op}^{M_i}_{^L
G}(D_\infty^\times)$).

Now we combine this result with the local statements of \thmref{factor
RS} describing the action of the center on
$\ghat_{\ka_c}$-modules. This gives us the sought-after assertions
about the factorization of the action of the universal Gaudin algebra
$\ZZ_{(z_i),\infty}(\g)$ on particular modules.

To complete the proof of part (1), we need to show that the
homomorphism
\begin{equation}    \label{inject}
\on{Fun} \on{Op}_{^L
  G}(\pone)^{(m_i),m_\infty}_{(z_i),\infty} \to \bigotimes_{i=1}^N
U(\g_{m_i}) \otimes U(\ol\g_{m_\infty})
\end{equation}
obtained this way is injective. But its image is, by definition, the
algebra $\ZZ_{(z_i),\infty}^{(m_i),m_\infty}(\g)$, and so the
injectivity of this homomorphism implies that we have an isomorphism
\begin{equation}    \label{isomor}
\on{Fun} \on{Op}_{^L G}(\pone)^{(m_i),m_\infty}_{(z_i),\infty}
\simeq \ZZ_{(z_i),\infty}^{(m_i),m_\infty}(\g),
\end{equation}
as stated in part (1) of the theorem.

In order to prove the injectivity we pass to the associated graded
algebras. According to \thmref{quantized}, at the level of associated
graded algebras the homomorphism \eqref{inject} becomes the
homomorphism
$$
\ol\Psi_{(z_i),\infty}^{(m_i),m_\infty}: \on{Fun} {\mc
H}_{G,(z_i),\infty}^{(m_i),m_\infty} \to \bigotimes_{i=1,\ldots,N}
S(\g_{m_i}) \otimes S(\ol\g_{m_\infty})
$$
defined in formula \eqref{olPsi}. By \lemref{is inj},
$\ol\Psi_{(z_i),\infty}^{(m_i),m_\infty}$ is injective. This implies
the injectivity of \eqref{inject}. Therefore we obtain an isomorphism
\eqref{isomor}.
\qed

\smallskip

Note that the embedding
$$
\on{Op}_{^L G}(\pone)^{(m_i),m_\infty}_{(z_i),\infty} \hookrightarrow 
\on{Op}_{^L G}(D_u)
$$
obtained by restricting an oper to the disc $D_u$ around a point $u
\in \pone \{ z_1,\ldots,z_N,\infty \}$ gives rise to a surjective
homomorphism
$$
\on{Fun} \on{Op}_{^L G}(D_u) \twoheadrightarrow \on{Fun} \on{Op}_{^L
G}(\pone)^{(m_i),m_\infty}_{(z_i),\infty}.
$$
The corresponding homomorphism
$$
\zz(\ghat)_u \twoheadrightarrow
\ZZ_{(z_i),\infty}^{(m_i),m_\infty}(\g)
$$
is nothing but the homomorphism $\Phi_{(z_i),\infty}^{(m_i),m_\infty}$
from \secref{universal}.

Observe also that part (1) of the theorem implies that the universal
Gaudin algebra $\ZZ_{(z_i),\infty}(\g)$ is isomorphic to the
(topological) algebra functions on the ind-affine space of all
(meromorphic) $^L G$-opers on $\pone \bs \{ z_1,\ldots,z_N,\infty \}$.

\subsection{The case of singularity of order $2$}

Consider a special case of the theorem when we have two points $z_1=0$
and $\infty$, set $m_1=1$ and $m_\infty=1$ and choose a regular $\chi
\in \g^*$ corresponding to the point $\infty$. The corresponding
Gaudin algebra $\A_{0,\infty}^{1,1}(\g)_{0,\chi}$ is the algebra
${\mc A}_\chi$ introduced in \secref{Achi}.

Let $\on{Op}_{^L G}(\pone)_{\pi(-\chi)}$ be the space of $^L G$-opers
on $\pone$ with regular singularity at the point $0$ and with
singularity of order $2$ at $\infty$, with $2$-residue
$\pi(-\chi)$. This space has the following concrete realization for
regular $\chi$.

Let us pick an element of the form
$$
-p_{-1}-\ol\chi = -p_{-1} - \sum_{j=1}^\ell \ol\chi_j p_j \in
{}^L \g_{\on{can}}
$$
(see formula \eqref{gcan} for the definition of $^L \g_{\on{can}}$) in
the conjugacy class $\pi(-\chi) \in \g^*/G = {}^L \g/{}^L G$. Then it
follows from \lemref{order two} that on the punctured disc
$D_\infty^\times$ at $\infty$ (with coordinate $s=t^{-1}$) each
element of $\on{Op}_{^L G}(\pone)_{\pi(-\chi)}$ may be uniquely
represented by a connection operator of the form
$$
\pa_s - p_{-1} - \sum_{j=1}^\ell (\ol\chi_j s^{-2d_j-2} + s^{-2d_j-1}
u_j(s)) p_j, \qquad u_j(s) = \sum_{n\geq 0} u_{j,n} s^n \in \C[[s]]
$$
(the sign in front of $p_{-1}$ may be eliminated by a gauge
transformation with $\crho(-1)$, which would result in multiplication
of $u_j(s)$ by $(-1)^j$, but we prefer not to do this). This oper
extends to $\pone \bs 0$ if and only if each $u_j(s)$ belongs to
$\C[s]$. To understand its behavior at $0$, we apply the change of
variables $s=t^{-1}$. After applying the gauge transformation with
$\crho(t^{-2})$, we find that the restriction of this oper to the
punctured disc $D_0^\times$ at $0 \in \pone$ is equal to
$$
\pa_t + p_{-1} + \frac{2\crho}{t} + \sum_{j=1}^\ell (\ol\chi_j +
\wt{u}_j(t)) p_j, \qquad \wt{u}_j(t) = t^{-1} u_j(t^{-1}).
$$
Next, we apply the gauge transformation with $\exp(-p_1/t)$ and obtain
$$
\pa_t + p_{-1} + \sum_{j=1}^\ell (\ol\chi_j + \wt{u}_j(t)) p_j, \qquad
\wt{u}_j(t) = t^{-1} u_j(t^{-1}).
$$
This oper has regular singularity at $0$ if and only if
$$
u_j(s) = \sum_{n=0}^{j} u_{j,n} s^n.
$$

Thus, we find that each element of $\on{Op}_{^L
  G}(\pone)_{\pi(-\chi)}$ may be represented uniquely by an operator
of the form
\begin{equation}    \label{order 2}
\pa_t + p_{-1} + \sum_{j=1}^\ell \left(\ol\chi_j + \sum_{n=0}^{d_j}
u_{j,n} t^{-n-1} \right) p_j.
\end{equation}
Note that according to \propref{can real RS}, its $1$-residue at $0$
is equal to
\begin{equation}    \label{1-res}
p_{-1} + \sum_{j=1}^\ell (u_{j,d_j} + \frac{1}{4} \delta_{j,1}) p_j \in
{}^L \g_{\on{can}} \simeq {}^L \g/{}^L G = \g/G.
\end{equation}
In particular, we obtain that
\begin{equation}    \label{free pol1}
\on{Fun} \on{Op}_{^L G}(\pone)_{\pi(-\chi)} \simeq
\C[u_{j,n_j}]_{j=1,\ldots,\ell;n_j=0,\ldots,d_j}.
\end{equation}

\begin{thm}    \label{isom for Achi}
If $\chi \in \g^*$ is regular, then the algebra ${\mc A}_\chi$ is
isomorphic to the algebra of functions on $\on{Op}_{^L
G}(\pone)_{\pi(-\chi)}$.
\end{thm}

\proof According to \thmref{ZZ descr},(4), we have a surjective
homomorphism
\begin{equation}    \label{hom to A}
\on{Fun} \on{Op}_{^L G}(\pone)_{\pi(-\chi)} \to {\mc A}_\chi.
\end{equation}
To show that it is an isomorphism, it is sufficient to prove that the
corresponding homomorphism of the associated graded algebras is an
isomorphism. According to \cite{Ryb2} and \thmref{th:quantisation},
$\on{gr} {\mc A}_\chi = \ol{\mc A}_\chi$. But for regular $\chi$ the
algebra $\ol{\mc A}_\chi$ is a free polynomial algebra with generators
$D_\chi^{n_i} \ol{P}_i, i=1,\ldots,\ell; n_i=0,\ldots,d_i+1$, by
\cite{MF,Ko2} and \thmref{poly}.

On the other hand, it follows from formula \eqref{free pol1},
\lemref{sign} and the discussion in the proof of
\thmref{th:quantisation} that $\on{gr} \on{Fun} \on{Op}_{^L
G}(\pone)_{\pi(-\chi)}$ is isomorphic to the same free polynomial
algebra. Hence the map \eqref{hom to A} is an isomorphism. \qed

\subsection{Joint eigenvalues on finite-dimensional modules}

Let us first recall the results of \cite{F:opers,F:faro} (see also
\cite{FFR,F:icmp}) on the joint eigenvalues of the ordinary Gaudin
algebra $\ZZ^{(1),1}_{(z_i),\infty}(\g) \subset U(\g)^{\otimes N}$ on
the tensor products $\bigotimes_{i=1}^N V_{\la_i}$ of irreducible
finite-dimensional $\g$-modules $V_{\la_i}$. Let $\on{Op}_{^L
G}(\pone)_{(z_i),\infty;(\la_i),\la_\infty}$ be the set of $^L
G$-opers on $\pone$ with regular singularities at $z_i, i=1,\ldots,N$
and $\infty$, with residues $\varpi(-\la_i-\rho)$ and
$\varpi(-\la_\infty-\rho)$, and with trivial monodromy
representation. Then according to \cite{F:opers}, Corollary 4.8 (see
also \cite{F:faro}, Theorem 2.7,(3), and \thmref{ZZ descr},(3) above),
we have the following description of the joint eigenvalues of
$\ZZ^{(1),1}_{(z_i),\infty}(\g) \subset U(\g)^{\otimes N}$ on
$\bigotimes_{i=1}^N V_{\la_i}$.

\begin{thm}    \label{reg sing eig}
The set of joint eigenvalues of $\ZZ^{(1),1}_{(z_i),\infty}(\g)$ on
$\bigotimes_{i=1}^N V_{\la_i}$ (without multiplicities) is a subset of
$\on{Op}_{^L G}(\pone)_{(z_i),\infty;(\la_i),\la_\infty}$.
\end{thm}

Further, Conjecture 1 of \cite{F:faro} states that this inclusion is
actually a bijection. Now we discuss analogous results and conjectures
for the generalized Gaudin algebras corresponding to irregular
singularity of order $2$ at $\infty$.

We start with the simplest such Gaudin algebra, namely, the algebra
${\mc A}_\chi \subset U(\g)$. Consider its action on the irreducible
finite-dimensional $\g$-module $V_\la$, where $\la$ is a dominant
integral weight. Note that from the point of view of the general
construction of \secref{multi-point} this action comes about through
the action of ${\mc A}_\chi = \A_{0,\infty}^{1,1}(\g)_{0,\chi}$ on
the space of coinvariants
$$
H(\V_\la \otimes \I_{1,\chi}) \simeq (V_\la \otimes I_\chi)/\g \simeq
V_\la.
$$
Suppose that $\chi$ is regular semi-simple. Then ${\mc A}_\chi$
contains a Cartan subalgebra $\h$ of $\g$, which is the centralizer of
$\chi$ in $\g$. Therefore the action of ${\mc A}_\chi$ preserves the
weight decomposition of $V_\la$ with respect to the $\h$-action. It is
natural to ask what are the joint generalized eigenvalues of ${\mc
A}_\chi$ on these components.

Let
$$
\on{Op}_{^L G}(\pone)^\la_{\pi(-\chi)} \subset \on{Op}_{^L
  G}(\pone)_{\pi(-\chi)}
$$
be the set of $^L G$-opers on $\pone$ with regular singularity at the
point $0$ with the $1$-residue $\varpi(-\la-\rho)$, singularity of
order $2$ at the point $\infty$ with the $2$-residue $\pi(-\chi)$ and
trivial monodromy. Then, according to \thmref{ZZ descr},(3), the
action of ${\mc A}_\chi$ on $U(\g)$ factors through the homomorphism
$$
{\mc A}_\chi \simeq \on{Fun} \on{Op}_{^L G}(\pone)_{\pi(-\chi)}
\twoheadrightarrow \on{Fun} \on{Op}_{^L G}(\pone)^\la_{\pi(-\chi)}.
$$
In other words, we obtain the following description of the joint
generalized eigenvalues of the commutative algebra ${\mc A}_\chi$
(which are the points of $\on{Op}_{^L G}(\pone)_{\pi(-\chi)}$,
according to \thmref{isom for Achi}) on its generalized eigenvectors
in $V_\la$.

\begin{thm}    \label{joint}
For regular semi-simple $\chi \in \g^*$ the set of joint generalized
eigenvalues of ${\mc A}_\chi$ on $V_\la$ (without multiplicities) is a
subset of $\on{Op}_{^L G}(\pone)^\la_{\pi(-\chi)}$.
\end{thm}

Concretely, elements of $\on{Op}_{^L G}(\pone)^\la_{\pi(-\chi)}$,
whose $1$-residue at $0$ is equal to $\varpi(-\la-\rho)$ are
represented by the connections of the form \eqref{order 2} with the
expression \eqref{1-res} equal to $\varpi(-\la-\rho)$. As we explained
after \propref{no mon}, the condition that this connection has trivial
monodromy around $0$ imposes a set of $\on{dim} {}^L \n$ algebraic
equations on the oper, and $\on{Op}_{^L G}(\pone)^\la_{\pi(-\chi)}$ is
just the set of solutions of these equations plus $\dim {}^L \h$
equations corresponding to the $1$-residue condition at $0$. Note that
the dimension of $\on{Op}_{^L G}(\pone)_{\pi(-\chi)}$ is equal to
$\on{dim} {}^L \bb$, and so it is reasonable to expect that
$\on{Op}_{^L G}(\pone)^\la_{\pi(-\chi)}$ is a finite set.

\begin{conj}    \label{conj for Achi}
The injective map of \thmref{joint} is a bijection for any regular
semi-simple element $\chi \in \g^*$ and dominant integral weight
$\la$.
\end{conj}

If the action of the algebra ${\mc A}_\chi$ on $V_\la$ is
diagonalizable (which we expect to happen for generic regular
semi-simple $\chi$), then \propref{joint} would give us a labeling of
an eigenbasis of ${\mc A}_\chi$ in $V_\la$ by monodromy-free $^L
G$-opers on $\pone$ with prescribed singularities at $0$ and $\infty$.

\medskip

\propref{joint} has the following multi-point generalization. Consider
the algebra
$$
\A_{(z_i),\infty}^{(1),1}(\g)_{(0),\chi} \subset U(\g)^{\otimes
  N}
$$
(see \secref{multi-point} and \cite{Ryb2}). It may be obtained as the
quotient of
$$
\ZZ_{(z_i),\infty}^{(1),2}(\g) \subset U(\g)^{\otimes N} \otimes
U(\ol\g_2)
$$
obtained by applying the character $\ol\g_2 \to \g
\overset{\chi}\longrightarrow \C$ along the last factor.

Let
$$
\on{Op}_{^L G}(\pone)^{(1),2}_{(z_i);\pi(-\chi)} \subset
\on{Op}_{^L G}(\pone)^{(1),2}_{(z_i),\infty}
$$
be the space of $^L G$-opers on $\pone$ with regular
singularities at the points $z_i, i=1,\ldots,N$, and with singularity
of order $2$ at the point $\infty$ with the $2$-residue $\pi(-\chi)$.

\begin{conj}
For any regular $\chi \in \g^*$
$$
\A_{(z_i),\infty}^{(1),1}(\g)_{(0),\chi} \simeq \on{Fun}
\on{Op}_{^L G}(\pone)^{(1),2}_{(z_i);\pi(-\chi)}.
$$
\end{conj}

Now let $\la_1,\ldots,\la_N$ be a collection of dominant integral
weights of $\g$. Consider the action of
$\A_{(z_i),\infty}^{(1),1}(\g)_{(0),\chi}$ on the tensor product
$\bigotimes_{i=1}^N V_{\la_i}$. We will now assume that $\chi$ is
regular semi-simple.

Let 
$$
\on{Op}_{^L G}(\pone)^{(\la_i)}_{(z_i);\pi(-\chi)} \subset
\on{Op}_{^L G}(\pone)^{(1),2}_{(z_i);\pi(-\chi)}
$$
be the set of $^L G$-opers on $\pone$ with regular singularities at
the points $z_i, i=1,\ldots,N$, with the $1$-residues
$\varpi(-\la_i-\rho)$ and trivial monodromy around these points, and
with singularity of order $2$ at the point $\infty$ with the
$2$-residue $\pi(-\chi)$. Then \thmref{ZZ descr} implies the following
result.

\begin{thm}    \label{subset}
There is an injective map from the set of joint generalized
eigenvalues of the commutative algebra
$\A_{(z_i),\infty}^{(1),1}(\g)_{(0),\chi}$ on $\bigotimes_{i=1}^N
V_{\la_i}$ (without multiplicities) to $\on{Op}_{^L
G}(\pone)^{(\la_i)}_{(z_i);\pi(-\chi)}$.
\end{thm}

We propose the following analogue of \conjref{conj for Achi}:

\begin{conj}    \label{bijection}
The injective map of \thmref{subset} is a bijection.
\end{conj}

This should be viewed as an analogue of Conjecture 1 of
\cite{F:faro}. The motivation for both conjectures comes from the
geometric Langlands correspondence (see the discussion in
\cite{F:faro} after Conjecture 1).

\section{Bethe Ansatz in Gaudin models with irregular singularities}
\label{bethe1}

In this section we develop an analogue of the Bethe Ansatz method for
constructing eigenvectors of the Gaudin algebra in the case of
irregular singularities. For definiteness, we will consider here the
case of the Gaudin algebras
$\A^{(1),1}_{(z_i),\infty}(\g)_{(0),\chi}$, but our methods may be
generalized to yield eigenvectors in other generalized Gaudin models.

The construction of the Bethe Ansatz for Gaudin models with regular
singularities is explained in detail in \cite{F:faro}, Sects. 4--5,
following \cite{FFR} (for another approach, see \cite{RV}). In this
section we will follow the same approach, using notation and results
of \cite{F:faro}.

\subsection{Wakimoto modules}

The construction of eigenvectors of the Hamiltonians of the ordinary
Gaudin model developed in \cite{FFR,F:faro} utilizes a class of
$\ghat_{\ka_c}$-modules called {\em Wakimoto modules}. These modules
were defined in \cite{Wak} for $\g=\sw_2$ and in
\cite{FF:usp,FF:si,F:wak} for a general simple Lie algebra $\g$ (for a
detailed exposition, see \cite{F:book}). Here we will follow the
notation of \cite{F:wak,F:faro}, where we refer the reader for more
details.

Wakimoto modules over $\ghat_{\ka_c}$ are parameterized by connections
on an $^L H$-bundle $\Omega^{-\rho}$ on the punctured disc
$D^\times = \on{Spec} \C\ppart$. Here $\Omega^{-\rho}$ is defined as
the push-forward of the $\C^\times$-bundle corresponding to the
canonical line bundle $\Omega$ on $D^\times$ under the homomorphism
$\C^\times \to {}^L H$ corresponding to the integral coweight $-\rho$
of $^L H$ (we recall that $^L H$ is a Cartan subgroup of the group $^L
G$ of adjoint type). A choice of coordinate $t$ on the disc $D$ gives
rise to a trivialization of $\Omega$, and hence of $\Omega^{-\rho}$. A
connection on $\Omega^{-\rho}$ may then be written as an operator
$$
\ol\nabla = \pa_t + \nu(t), \qquad \nu(t) \in {}^L \h\ppart =
\h^*\ppart
$$
(see \cite{F:wak}, Sect. 5.5). If $s$ is another coordinate such that
$t=\varphi(s)$, then this connection will be represented by the
operator
\begin{equation} \label{trans for conn1}
\pa_s + \varphi'(s) \nu(\varphi(s)) + \rho \cdot
\frac{\varphi''(s)}{\varphi'(s)}.
\end{equation}

Let $\Con_{D^\times}$ be the space of all connections on the $^L
H$-bundle $\Omega^{-\rho}$ on $D^\times$. Denote by $b_{i,n}, n \in
\Z$, the function on $\Con_{D^\times}$ defined by the formula
$$
\pa_t + \nu(t) \mapsto \on{Res}_{t=0} \langle \chal_i,\nu(t) \rangle t^n
dt.
$$
The algebra $\on{Fun} \Con_{D^\times}$ of functions on
$\Con_{D^\times}$ is a complete topological algebra
$$
\on{Fun} \Con_{D^\times} \simeq \underset{\longleftarrow}\lim \;
\C[b_{i,n}]_{i=1,\ldots,\ell;n\in\Z}/I_N,
$$
where $I_N$ is the ideal generated by $b_{i,n}, i=1,\ldots,\ell; n\geq
N$. A module over this algebra is called smooth if every vector is
annihilated by an ideal $I_N$ for large enough $N$. In particular,
each $\ol\nabla \in \Con_{D^\times}$ gives rise to a one-dimensional
smooth module over $\on{Fun} \Con_{D^\times}$, which we denote by
$\C_{\ol\nabla}$. Equivalently, $\C_{\ol\nabla}$ may be viewed as a
module over the commutative vertex algebra $\pi_0 =
\on{Fun} \Con_{D}$ (see \cite{F:wak}, Sect. 4.2).

Next, we define the Weyl algebra $\cA^{\g}$ with generators
$a_{\al,n}, a^*_{\al,n}$, $\al \in \De_+$, $n \in \Z$, and relations
\begin{equation} \label{commina}
[a_{\al,n},a_{\beta,m}^*] = \delta_{\al,\beta}
\delta_{n,-m}, \hskip.3in
[a_{\al,n},a_{\beta,m}] = [a_{\al,n}^*,a_{\beta,m}^*] = 0.
\end{equation}
Let $M_{\g}$ be the Fock representation of $\cA^{\g}$ generated by a
vector $\vac$ such that
$$
a_{\al,n} \vac = 0, \quad n\geq 0; \qquad a^*_{\al,n} \vac = 0, \quad
n>0.
$$
It carries a vertex algebra structure (see \cite{F:wak}). It follows
from the general theory of \cite{FB}, Ch. 5, that a module over the
vertex algebra $M_\g$ is the same as a smooth module over the Weyl
algebra $\cA^{\g}$, i.e., one such that every vector is annihilated by
$a_{\al,n}, a^*_{\al,n}$ for large enough $n$.

According to \cite{F:wak}, Theorem 4.7, we have a homomorphism of
vertex algebras
$$
w_{\ka_c}: \V_0 \to M_\g \otimes \pi_0,
$$
which sends the center $\zz(\ghat) \subset \V_0$ to $\pi_0$. Moreover,
the corresponding homomorphism of commutative algebras
$$
\on{Fun} \on{Op}_{^L G}(D) \to \on{Fun} \Con_D
$$
is induced by the {\em Miura transformation}
$$
\Con_D \to \on{Op}_{^L G}(D)
$$
defined in \cite{F:wak}, Sect. 10.3, following \cite{DS}. We recall
that this map sends a Cartan connection $\ol\nabla = \pa_t + \nu(t)$
to the $^L G$-oper which is the $N[[t]]$-gauge equivalence class of
$\nabla = \pa_t + p_{-1} + \nu(t)$. There is a similar map over the
punctured disc $D^\times$ or a smooth curve.

This result implies that for any smooth $\cA^\g$-module $L$ and any
smooth module $R$ over $\on{Fun} \Con_{D^\times}$ the tensor product
$L \otimes R$ is a $\ghat_{\ka_c}$-module. In particular, taking $R =
\C_{\ol\nabla}$, we obtain a $\ghat_{\ka_c}$-module $L \otimes
\C_{\ol\nabla}$. These are the Wakimoto modules. The following is
proved in \cite{F:wak}, Theorem 12.6.

\begin{thm}    \label{center and miura}
The center $Z(\ghat) = \on{Fun} \on{Op}_{^L G}(D^\times)$ acts on $L
\otimes \C_{\ol\nabla}$ via the character corresponding to the $^L
G$-oper $\ol{\mb b}^*(\ol\nabla)$, where
$$
\ol{\mb b}^*(\ol\nabla): \Con_{D^\times} \to \on{Op}_{^L G}(D^\times)
$$
is the Miura transformation on $D^\times$. In particular, the action
of $Z(\ghat)$ is independent of the choice of the module $L$.
\end{thm}

\subsection{Coinvariants of Wakimoto modules}

The idea of \cite{FFR,F:faro} is to construct eigenvectors of Gaudin
algebras using spaces of coinvariants of the tensor products of
Wakimoto modules. These coinvariants are defined (following the
general definition of \cite{FB}, Ch. 10) with respect to the vertex
algebra $M_\g \otimes \pi_0$. By functoriality of coinvariants, the
above homomorphism $w_{\ka_c}$ of vertex algebras gives rise to linear
maps from the spaces of coinvariants with respect to $\V_{0} =
\V_{0,\ka_c}$, which are just the spaces of $\g_{(z_i)}$-coinvariants
introduced in \secref{construction}, to the spaces of coinvariants
with respect to $M_\g \otimes \pi_0$. Because $M_\g \otimes \pi_0$ is
a much simpler vertex algebra, its coinvariants are easy to compute,
and in the cases we consider below they turn out to be
one-dimensional. Therefore we obtain linear functionals on the spaces
of $\g_{(z_i)}$-coinvariants that are of interest to us. It then
follows from the construction that these linear functionals are
eigenvectors of the corresponding Gaudin algebra.

Let us choose a set of distinct points $x_1,\ldots,x_p$ on $\C \subset
\pone$. We attach to each of these points, and to the point
$\infty$, an $M_\g$-module and a $\pi_0$-module. As the $M_\g$-modules
attached to $x_i, i=1,\ldots,p$, we choose $M_\g$ itself, and to
$\infty$ we attach another module $M'_\g$ generated by a vector
$\vac'$ such that
$$
a_{\al,n} \vac' = 0, \quad n > 0; \qquad a^*_{\al,n} \vac' = 0, \quad
n \geq 0.
$$
As the $\pi_0$-modules attached to $x_i, i=1,\ldots,p$, we
take the one-dimensional modules $\C_{\ol\nabla_i} = \C_{\nu_i(z)}$,
where
$$
\ol\nabla_i = \pa_z + \nu_i(z), \qquad \nu_i(z) \in \h^*\zpart = {}^L
\h\zpart,
$$
and as the $\pi_0$-module attached to the point $\infty$ we take
$\C_{\ol\nabla_\infty} = \C_{\nu_\infty(z)}$, where $\ol\nabla_\infty
= \pa_z + \nu_\infty(z)$.

The corresponding space of coinvariants for $M_\g \otimes \pi_0$ is
the tensor product of the spaces
$$
H_{M_\g}(\pone;(x_i),\infty;(M_\g),M'_\g) \qquad \on{and} \qquad
H_{\pi_0}(\pone,(x_i),\infty;(\C_{\nu_i(z)}),\C_{\nu_\infty}(z))
$$
of coinvariants for the $M_\g$-modules and the $\pi_0$-modules,
respectively (see \cite{F:faro} for their definition). The following
result is proved in \cite{F:faro}, Prop. 4.9 (see also \cite{FFR},
Prop. 4).

\begin{prop}    \label{analyt} \

{\em (1)} The space $H_{M_\g}(\pone;(x_i),\infty;(M_\g),M'_\g)$ is
one-dimensional and the projection of the vector $\vac^{\otimes N}
\otimes \vac'$ on it is non-zero.

{\em (2)} The space
$H_{\pi_0}(\pone,(x_i),\infty;(\C_{\nu_i(z)}),\C_{\nu_\infty}(z))$ is
one-dimensional if and only if there exists a connection $\ol{\nabla}$
on the $^L H$-bundle $\Omega^{-\rho}$ on $\pone \bs \{
x_1,\ldots,x_p,\infty \}$ whose restriction to the punctured disc at
each $x_i$ is equal to $\pa_t + \nu_i(t-x_i)$, and whose restriction
to the punctured disc at $\infty$ is equal to $\pa_{t^{-1}} +
\nu_\infty(t^{-1})$.

Otherwise,
$H_{\pi_0}(\pone,(x_i),\infty;(\C_{\nu_i(z)}),\C_{\nu_\infty}(z))=0$.
\end{prop}

Formula \eqref{trans for conn1} shows that if we have a connection on
$\Omega^{-\rho}$ over $\pone$ whose restriction to $\pone \bs \infty$
is represented by the operator $\pa_t + \nu(t)$, then its restriction
to the punctured disc $D_\infty^\times$ at $\infty$ reads, with
respect to the coordinate $s = t^{-1}$
$$
\pa_s - s^{-2} \nu(s^{-1}) - 2\rho s^{-1}.
$$

In \cite{FFR,F:faro} we chose $\nu_i(z)$ to be of the form $$\nu_i(z) =
\frac{\nu_i}{z} + \sum_{n\geq 0} \nu_{i,n} z^n,$$ and $\nu_\infty(z)$ to
be of the form
$$\nu_\infty(z) = \frac{\nu_\infty}{z} + \sum_{n\geq 0} \nu_{\infty,n}
z^n.$$
In other words, we consider connections on $\Omega^{-\rho}$ with
regular singularities. The condition of the proposition is then
equivalent to saying that the restriction of $\ol\nabla$ to $\pone \bs
\infty$ is represented by the operator $\pa_t + \nu(t)$, where
$$
\nu(t) = \sum_{i=1}^p \frac{\nu_i}{t-x_i},
$$
and $\pa_t + \nu_i(t-x_i)$ is the expansion of $\nu(t)$ at $x_i,
i=1,\ldots,p$, while $\nu_\infty(t^{-1})$ is the expansion of $- t^2
\nu(t) - 2\rho t$ in powers of $t^{-1}$.

Now we will choose $\nu_i(z)$ to be the same as above, but we will
choose $\nu_\infty(z)$ to be of the form
$$
\nu_\infty(z) = \frac{\chi}{z^2} + \sum_{n\geq -1} \nu_{\infty,n} z^n,
$$
where $\chi \in \h^*$. Then the condition of the
proposition means that the restriction of $\ol\nabla$ to $\pone \bs
\infty$ is represented by the operator
$$
\pa_t - \chi + \sum_{i=1}^p \frac{\nu_i}{t-x_i}.
$$
This connection has irregular singularity at $\infty$. Note also that
we have $\nu_{\infty,-1} = -2\rho - \sum_{i=1}^p \nu_i$.

Using the homomorphism $w_{\ka_c}$, we obtain $\ghat_{\ka_c}$-module
structures on $M_\g \otimes \C_{\nu_i(z)}$ and on $M'_\g \otimes
\C_{\nu_\infty(z)}$. As explained in \cite{F:faro}, Sect. 4.2,
functoriality of coinvariants implies that there is a natural map from
the space of $\g_{(x_i)}$-coinvariants $$H(M_\g \otimes
\C_{\nu_1(z)},\ldots,M_\g \otimes \C_{\nu_p(z)},M'_\g \otimes
\C_{\nu_\infty(z)}),$$ defined in \secref{coinvariants}, to
the corresponding space of coinvariants with respect to $M_\g \otimes
\pi_0$, which is
$$
H_{M_\g}(\pone;(x_i),\infty;(M_\g),M'_\g) \otimes
H_{\pi_0}(\pone,(x_i),\infty;(\C_{\nu_i(z)}),\C_{\nu_\infty}(z)).
$$

Now, if $\nu_i(z)$ and $\nu_\infty(z)$ are as above, the latter space
is one-dimensional, by \propref{analyt}. Hence we obtain a non-zero
linear functional
\begin{equation} \label{non zero f}
\tau_{(x_i)}: H(M_\g \otimes \C_{\nu_1(z)},\ldots,M_\g
\otimes \C_{\nu_p(z)},M'_\g \otimes \C_{\nu_\infty}(z)) \to \C,
\end{equation}
which we normalize so that its value on $\vac^{\otimes N} \otimes
\vac'$ is equal to $1$.

\subsection{Construction of the Bethe vectors}

We are now ready to construct eigenvectors of the Gaudin algebra
$\A^{(1),1}_{(z_i),\infty}(\g)_{(0),\chi}$. The idea is to use the
space of coinvariants of the tensor product of specially selected
Wakimoto modules. We attach them to the points $z_i, i=1,\ldots,\ell$,
and $\infty$, and also to additional points $w_1,\ldots,w_m$. The
modules attached to the points $z_i, i=1,\ldots,\ell$ will be of the
form $M_\g \otimes \C_{\la_i(z)}$, where the connection
$\pa_z+\la_i(z)$ has regular singularity. For such modules we have a
homomorphism $\M^*_{\la_i} \to M_\g \otimes \C_{\la_i(z)}$, where
$\la_i$ is the most singular coefficient of $\la_i(z)$. The module
attached to $\infty$ will be of the form $M'_\g \otimes
\C_{\la_\infty(z)}$, where the connection $\pa_z + \la_\infty(z)$ has
singularity of order $2$. We then have a homomorphism $\I_{1,\chi} \to
M'_\g \otimes \C_{\la_\infty(z)}$, where $\chi$ is the most singular
coefficient of $\la_\infty(z)$. Finally, the module attached to $w_j$
will be of the form $M_\g \otimes \C_{\mu_j(z)}$, where $\mu_j(z) =
-\al_{i_j}/z + \ldots$. Considered as a $\ghat_{\ka_c}$-module, this
module contains a vector annihilated by $\g[[t]]$, provided that a
certain system of equations, called {\em Bethe Ansatz equations}, is
satisfied. If it is satisfied, then we can use these
$\g[[t]]$-invariant vectors to construct eigenvectors of the Gaudin
algebra $\A^{(1),1}_{(z_i),\infty}(\g)_{(0),\chi}$ in
$\bigotimes_{i=1}^N M_{\la_i}$.

Now we explain this in more detail. Let us look more closely at the
Wakimoto modules $M_\g \otimes \C_{\nu_i(z)}$ and $M'_\g \otimes
\C_{\nu_\infty(z)}$. Let $M_\la^*$ be the $\g$-module contragredient
to the Verma module $M_\la$, and $\M^*_\la$ the corresponding induced
$\ghat_{\ka_c}$-module. As shown in \cite{F:faro}, Sect. 4.4, we have
a non-trivial homomorphism of $\ghat_{\ka_c}$-modules
$$
\M^*_\la \to M_\g \otimes \C_{\nu(z)}, \qquad \on{if} \quad \nu(z) =
\frac{\nu}{z} + \sum_{n \geq 0} \nu_n z^n.
$$

Now suppose that $\nu = -\al_i, i = 1,\ldots,\ell$. Consider
the vector
\begin{equation}    \label{null vector}
e^R_{i,-1}\vac \in M_\g \otimes \C_{\nu(z)}, \qquad \nu(z) =
-\frac{\al_i}{z} + \sum_{\n \geq 0} \nu_n z^n
\end{equation}
(see \cite{F:faro}, formula (4.7) for the definition of
$e^R_{i,-1}$).

According to \cite{FFR}, Lemma 2 (see also \cite{F:faro}, Lemma 4.5),
we have

\begin{lem}    \label{sing}
The vector \eqref{null vector} is annihilated by $\g[[t]]$ if and only
if we have
\begin{equation}    \label{cond bethe}
\langle \chal_i,\nu_0 \rangle = 0.
\end{equation}
\end{lem}

Finally, consider the $\ghat_{\ka_c}$-module $M'_\g \otimes
\C_{\nu(z)}$, where
\begin{equation}    \label{uz}
\nu(z) = \frac{\chi}{z^2} + \sum_{n\geq -1} \nu_n z^n, \qquad \chi \in
\h^*.
\end{equation}
Using the explicit formulas for the homomorphism $w_{\ka_c}$ (see
\cite{F:faro}, Theorem 4.1), we obtain that the vector $\vac'
\in M'_\g \otimes \C_{\nu(z)}$ satisfies
$$
t^2 \g[[t]] \cdot \vac' = 0, \qquad (A \otimes t) \cdot \vac' =
\chi(A) \vac', \quad A \in \g,
$$
where we extend $\chi \in \h^*$ to a linear functional on $\g^*$ via
the projection $\g \to \h$ obtained using the Cartan decomposition of
$\g$ (abusing notation, we will denote this extension by the same
symbol $\chi$).

This implies the following:

\begin{lem}
For any $\chi \in \h^*$ there is a homomorphism of
$\ghat_{\ka_c}$-modules
$$
\I_{1,\chi} \to M'_\g \otimes \C_{\nu(z)},
$$
where $\nu(z)$ is given by formula \eqref{uz}, sending the generating
vector of $\I_{1,\chi}$ to $\vac' \in M'_\g \otimes \C_{\nu(z)}$.
\end{lem}

Now let us fix an $N$-tuple of distinct points $z_1,\ldots,z_N$ of $\C
= \pone \bs \{ \infty \}$, an element $\chi \in \h^*$, a set of
weights $\la_i \in \h^*, i=1,\ldots,N$, and a set of simple roots
$\al_{i_j}, j=1,\ldots,m$, of $\g$. Consider a connection
$\ol{\nabla}$ on $\Omega^{-\rho}$ on $\pone$ whose restriction to
$\pone \bs \infty$ is equal to $\pa_t + \la(t)$, where
\begin{equation} \label{sv}
\la(t) = - \chi + \sum_{i=1}^N \frac{\la_i}{t-z_i} - \sum_{j=1}^m
\frac{\al_{i_j}}{t-w_j},
\end{equation}
where $w_1,\ldots,w_m$ is an $m$-tuple of distinct points of $\pone
\bs \{ \infty \}$ such that $w_j \neq z_i$. Denote by $\la_i(t-z_i)$
the expansions of $\la(t)$ at the points $z_i, i=1,\ldots,N$, and by
$\mu_j(t-w_j)$ the expansions of $\la(t)$ at the points $w_j,
j=1,\ldots,m$. We have: $$\la_i(z) = \frac{\la_i}{z} + \cdots, \quad
\quad \mu_j(z) = - \frac{\al_{i_j}}{z} + \mu_{j,0} + \cdots,$$ where
\begin{equation}    \label{muj0}
\mu_{j,0} = - \chi + \sum_{i=1}^N \frac{\la_i}{w_j-z_i} - \sum_{s\neq
  j} \frac{\al_{i_s}}{w_j-w_s}.
\end{equation}

Finally, the expansion of this connection near $\infty$ reads
$\pa_s + \la_\infty(s)$, where $s=t^{-1}$. Then $\la_\infty(s) =
s^{-2} \la(s^{-1}) - 2\rho s^{-1}$. Therefore we have
\begin{equation} \label{la infty}
\la_\infty(z) = \frac{\chi}{z^2} - \frac{\sum_{i=1}^N \la_i -
  \sum_{j=1}^m \al_{i_j} + 2\rho}{z} + \cdots.
\end{equation}
In the previous subsection we constructed a non-zero linear
functional $\tau_{(z_i),(w_j)}$ on the corresponding space of
$\g_{(z_i),(w_j)}$-coinvariants
\begin{equation}    \label{tauzw}
\tau_{(z_i),(w_j)}: H((M_\g \otimes \C_{\la_i(z)}),(M_\g \otimes
\C_{\mu_j(z)}), M'_\g \otimes \C_{\la_\infty(z)}) \to \C.
\end{equation}
(In particular, this implies that this space of coinvariants is itself
non-zero.)

Next, we use \lemref{sing}, according to which the vectors
$e^R_{i_j,-1} \vac \in M_\g \otimes \C_{\mu_j(z)}$ are
$\g[[t]]$-invariant if and only if the equations $\langle
\chal_{i_j},\mu_{j,0} \rangle = 0$ are satisfied, where $\mu_{j,0}$ is
the constant coefficient in the expansion of $\la(t)$ at $w_j$ given
by formula \eqref{muj0}. This yields the following system of
equations:
\begin{equation} \label{bethe}
\sum_{i=1}^N \frac{\langle \chal_{i_j},\la_i \rangle}{w_j-z_i} -
\sum_{s\neq j} \frac{\langle \chal_{i_j},\al_{i_s} \rangle}{w_j-w_s} =
\langle \chal_{i_j},\chi \rangle, \quad j=1,\ldots,m.
\end{equation}
These are the {\em Bethe Ansatz equations} of our Gaudin model.

This is a system of equations on the complex numbers $w_j,
j=1,\ldots,m$, to each of which we attach a simple root $\al_{i_j}$.
We have an obvious action of a product of symmetric groups permuting
the points $w_j$ corresponding to simple roots of the same kind. In
what follows, by a {\em solution} of the Bethe Ansatz equations we
will understand a solution defined up to these permutations.  We will
adjoin to the set of all solutions associated to all possible
collections $\{ \al_{i_j} \}$ of simple roots of $\g$, the unique
``empty'' solution, corresponding to the empty set of simple roots
(when this system of equations is empty).

Suppose that equations \eqref{bethe} are satisfied. Then we obtain a
homomorphism of $\G_{\ka_c}$-modules
$$
\bigotimes_{i=1}^N \M^*_{\la_i} \otimes \V_0^{\otimes m} \otimes
\I_{1,\chi} \to \bigotimes_{i=1}^N M_\g \otimes \C_{\la_i(z)} \otimes
\bigotimes_{j=1}^m M_\g \otimes \C_{\mu_j(z)} \otimes (M'_\g \otimes
\C_{\la_\infty(z)}),
$$
which sends the vacuum vector $v_0$ in the $j$th copy of $\V_0$ to
$e^R_{i_j,-1} \vac \in M_\g \otimes \C_{\mu_j(z)}$. Hence we obtain
the corresponding map of the spaces of $\g_{(z_i),(w_j)}$-coinvariants
$$
H((\M^*_{\la_i}),(\V_0),\I_{1,\chi}) \to H((M_\g \otimes
\C_{\la_i(z)}),(M_\g \otimes \C_{\mu_j(z)}),M'_\g \otimes
\C_{\la_\infty(z)}).
$$
But the insertion of $\V_0$ does not change the space of coinvariants,
according to \propref{pr:propagation}. Hence the first space of
coinvariants is isomorphic to the space of $\g_{(z_i)}$-coinvariants
$$
H((\M^*_{\la_i}),\I_{1,\chi}) \simeq (\bigotimes_{i=1}^N M^*_{\la_i}
\otimes I_\chi)/\g \simeq \bigotimes_{i=1}^N M^*_{\la_i}.
$$
Composing the corresponding map
$$
\bigotimes_{i=1}^N M^*_{\la_i} \to H((M_\g \otimes
\C_{\la_i(z)}),(M_\g \otimes \C_{\mu_j(z)}),M'_\g \otimes
\C_{\la_\infty(z)})
$$
with the linear functional $\tau_{(z_i),(w_j)}$ (see \eqref{tauzw}),
we obtain a linear functional \begin{equation} \label{linfunct}
\psi(w_1^{i_1},\ldots,w_m^{i_m}): \bigotimes_{i=1}^N M^*_{\la_i} \to
\C.
\end{equation}

The functional $\psi(w_1^{i_1},\ldots,w_m^{i_m})$ coincides with the
functional denoted in the same way in \cite{F:faro},
Sect. 4.4. Indeed, it is obtained by computing the same coinvariants
as in the setting of \cite{F:faro} and here. What is different here is
that the numbers $w_j$'s are solutions of the Bethe Ansatz equations
\eqref{bethe}, whereas in \cite{F:faro} they are the solutions of
another set of equations, namely equations (4.18) of \cite{F:faro},
which correspond to the special case $\chi=0$ (note that in this case
$\I_{1,\chi} = \UU_1$, which is the module that was attached to the
point $\infty$ in \cite{F:faro}).

According to formula (4.22) of \cite{F:faro} (based on the
computations performed in \cite{ATY} and in \cite{FFR}, Lemma 3),
$\psi(w_1^{i_1},\ldots,w_m^{i_m})$ corresponds to the vector
$$
\phi(w_1^{i_1},\ldots,w_m^{i_m}) \in \bigotimes_{k=1}^N M_{\la_k}
$$
given by the formula
\begin{equation} \label{genbv}
\phi(w_1^{i_1},\ldots,w_m^{i_m}) = \sum_{p=(I^1,\ldots,I^N)} \;
\bigotimes_{k=1}^N \; \; \prod_{s \in I^k} \;
\frac{f_{i_s}}{(w_s^{i_s}-w_{s+1}^{i_{s+1}})} \cdot v_{\la_k}
\end{equation}
(up to a scalar). Here $f_i$ denotes a generator of the Lie
algebra $\n_- \subset \g$ corresponding to the $i$th simple root, and
$v_{\la_k}$ is a highest weight vector in $M_{\la_k}$. The
summation is taken over {\em all ordered partitions} $I^1 \cup I^2
\cup \ldots \cup I^N$ of the set $\{1,\ldots,m\}$, where $I^k = \{
j^k_1,j^k_2,\ldots,j^k_{a_j} \}$, and the product is taken from
left to right, with the convention that the $w$ with the lower index
$j^k_{a_j}+1$ is $z_k$. Note that we differentiate between partitions
obtained by permuting elements within each subset $I^k$.


Thus, $\psi(w_1^{i_1},\ldots,w_m^{i_m})$ is the linear functional on
$\bigotimes_{i=1}^N M^*_{\la_i}$ equal to the pairing with
$\phi(w_1^{i_1},\ldots,w_m^{i_m})$.

In particular, we find that $\phi(w_1^{i_1},\ldots,w_m^{i_m})$ has
weight
\begin{equation}    \label{weight}
\sum_{i=1}^N \la_i - \sum_{j=1}^m \al_{i_j}.
\end{equation}

We call $\phi(w_1^{i_1},\ldots,w_m^{i_m})$ the {\em Bethe vector}
corresponding to a solution of the Bethe Ansatz equations
\eqref{bethe}.

\begin{thm}    \label{it is eigenvector}
If the Bethe Ansatz equations \eqref{bethe} are satisfied, then the
vector \linebreak $\phi(w_1^{i_1},\ldots,w_m^{i_m})$ given by formula
\eqref{genbv} is either equal to zero or is an eigenvector of the
Gaudin algebra $\A^{(1),1}_{(z_i),\infty}(\g)_{(0),\chi}$ in
$\bigotimes_{i=1}^N M_{\la_i}$.
\end{thm}

\noindent {\em Proof} is identical to the proof of Theorem 4.11 of
\cite{F:faro} (which is based on \cite{FFR}, Theorem 3). Let us look
at the linear functional $\psi(w_1^{i_1},\ldots,w_m^{i_m})$ as a map
of coinvariants
$$
\bigotimes_{i=1}^N M^*_{\la_i} \simeq H((\M^*_{\la_i}),\I_{1,\chi})
\to \\ H((M_\g \otimes
\C_{\la_i(z)}),(M_\g \otimes \C_{\mu_j(z)}),M'_\g \otimes
\C_{\la_\infty(z)}) \simeq \C.
$$
By functoriality of coinvariants, this map intertwines the natural
actions of $$\zz(\ghat)_u = \on{Fun} \on{Op}_{^L G}(D_u)$$ on the left
and right hand sides. By \thmref{ZZ descr}, on the left hand side this
action corresponds to the action of the Gaudin algebra
$\A^{(1),1}_{(z_i),\infty}(\g)_{(0),\chi}$ on $\bigotimes_{i=1}^N
M^*_{\la_i}$. On the other hand, by \thmref{center and miura}, the
action of $\on{Fun} \on{Op}_{^L G}(D_u)$ on the right hand side
factors through the Miura transformation
$$
\on{Fun} \on{Op}_{^L G}(D_u) \to \on{Fun} \Con_{D_u}.
$$
This means that $\psi(w_1^{i_1},\ldots,w_m^{i_m})$ (and hence
$\phi(w_1^{i_1},\ldots,w_m^{i_m})$) is an eigenvector of $\zz(\ghat)_u
= \on{Fun} \on{Op}_{^L G}(D_u)$ whose eigenvalue is the $^L G$-oper on
$\pone$ whose restriction to $D_u$ is the Miura transformation of the
connection $\pa_t + \la(t)$ restricted to $D_u$. \qed

\smallskip

In particular, $\phi(w_1^{i_1},\ldots,w_m^{i_m})$ is an eigenvector of
the operators $\Xi_{i,\chi}, i=1,\ldots,N$, given by formula
\eqref{Gaudin ham1} and of the DMT Hamiltonians $T_\ga(\chi)$ (acting
via the diagonal action of $\g$ on $\bigotimes_{i=1}^N
M_{\la_i}$). The Bethe Ansatz procedure for these quadratic
Hamiltonians was considered, from a different perspective, in
\cite{MTV}, for $\g=\sw_n$. The results of \cite{MTV} are
in agreement with \thmref{it is eigenvector}. We also note that in the
case when $\g=\sw_2$ the diagonalization problem for the operators
$\Xi_{i,\chi}, i=1,\ldots,N$, was studied in \cite{Skl}.

\subsection{Eigenvalues on Bethe vectors}

Next, we compute the joint eigenvalues of the algebra
$\A^{(1),1}_{(z_i),\infty}(\g)_{(0),\chi}$ on the Bethe vector
$\phi(w_1^{i_1},\ldots,w_m^{i_m})$.

It follows from \thmref{ZZ descr}
that the action of $\A^{(1),1}_{(z_i),\infty}(\g)_{(0),\chi}$ on the
space of coinvariants $H((\M^*_{\la_i}),\I_{1,\chi})$ factors through the
algebra
$$
\on{Fun} \on{Op}^{\on{RS}}_{^L G}(\pone)^{(\la_i)}_{(z_i);\pi(-\chi)},
$$
where $\on{Op}^{\on{RS}}_{^L G}(\pone)^{(\la_i)}_{(z_i);\pi(-\chi)}$
is the space of $^L G$-opers on $\pone$ with regular singularities at
the points $z_i, i=1,\ldots,N$, with the $1$-residues
$\varpi(-\la_i-\rho)$, and with singularity of order $2$ at the point
$\infty$ with the $2$-residue $\pi(-\chi)$. Therefore the joint
eigenvalues of $\A^{(1),1}_{(z_i),\infty}(\g)_{(0),\chi}$ on
$\phi(w_1^{i_1},\ldots,w_m^{i_m})$ are recorded by a point in
$\on{Op}^{\on{RS}}_{^L G}(\pone)^{(\la_i)}_{(z_i);\pi(-\chi)}$.

To describe this point, we let
$$
\Con_{\pone \bs \{ (z_i),(w_j),\infty \}}
$$
denote the space of connections on $\Omega^{-\rho}$ over $\pone \bs \{
(z_i),(w_j),\infty \}$. Then we have a Miura transformation
$$
\mu_{(z_i),(w_j),\infty}: \Con_{\pone \bs \{ (z_i),(w_j),\infty
\}} \to \on{Op}_{^L G}(\pone \bs \{ (z_i),(w_j),\infty \}).
$$

\begin{lem}
Let $\ol\nabla = \pa_t + \la(t)$, where $\la(t)$ is
given by formula \eqref{sv}, and the numbers $w_j$ satisfy the Bethe
Ansatz equations \eqref{bethe}. Then
$$
\mu_{(z_i),(w_j),\infty}(\ol\nabla) \in \on{Op}^{\on{RS}}_{^L
G}(\pone)^{(\la_i)}_{(z_i);\pi(-\chi)} \subset \on{Op}_{^L G}(\pone \bs \{
(z_i),(w_j),\infty \}).
$$
\end{lem}

\proof We need to show that the restriction of
$\mu_{(z_i),(w_j),\infty}(\ol\nabla)$ to $D_{w_j}^\times,
j=1,\ldots,m$, belongs to the subspace $\on{Op}_{^L G}(D_{w_j})$ of
regular opers, and the restriction of
$\mu_{(z_i),(w_j),\infty}(\ol\nabla)$ to $D_{\infty}^\times$ belongs
to $\on{Op}^{\leq 2}(D_\infty)_{\pi(-\chi)}$.

To see the former, we recall that the restriction of the connection
$\ol\nabla$ to $D_{w_j}^\times$ has the form
$$
\pa_t - \frac{\al_{i_j}}{t-w_j} + \mu_{j,0} + \cdots,
$$
where $\mu_{j,0}$ is given by formula \eqref{muj0}. The Bethe
Ansatz equations \eqref{bethe} are equivalent to the equations $\langle
\chal_{i_j},\mu_{j,0} \rangle = 0$. According to \cite{F:faro},
Lemma 3.5, this implies that the Miura transformation of this
connection is regular at $w_j$.

To prove the latter, we recall that the restriction of $\ol\nabla$ to
$D_\infty^\times$ has the form \eqref{la infty}:
$$
\pa_s + \frac{\chi}{s^2} + \cdots, \qquad s=t^{-1}.
$$
The Miura transformation of this connection is the $^L G$-oper which
is the $N(\!(s)\!)$-gauge equivalence class of the operator
$$
\nabla = \pa_s + p_{-1} + \frac{\chi}{s^2} + \cdots.
$$
Applying the gauge transformation with $\crho(s)^2$, we obtain the
connection
$$
\pa_s + \frac{1}{s^2}(p_{-1} + \chi) - \frac{2\crho}{s} + \cdots.
$$
This oper has singularity of order $2$, and its $2$-residue is equal to
$\pi(-p_{-1}-\chi) = \pi(-\chi)$, since $\chi \in \h^* \subset \g^*$.
\qed

\smallskip

Now we are ready to describe the joint eigenvalues of
$\A^{(1),1}_{(z_i),\infty}(\g)_{(0),\chi}$ on the Bethe vectors.

\begin{thm}    \label{joint eigenvalue}
The joint eigenvalues of $\A^{(1),1}_{(z_i),\infty}(\g)_{(0),\chi}$
on
$$
\phi(w_1^{i_1},\ldots,w_m^{i_m}) \in \bigotimes_{i=1}^N M_{\la_i},
$$
where $w_1,$ $\ldots,w_m$ satisfy the Bethe Ansatz equations
\eqref{bethe}, are given by the $^L G$-oper
$$
\mu_{(z_i),(w_j),\infty}(\ol\nabla) \in \on{Op}_{^L
  G}^{\on{RS}}(\pone)^{(\la_i)}_{(z_i);\pi(-\chi)}
$$
where $\ol\nabla = \pa_t + \la(t)$ with $\la(t)$ given by formula
\eqref{sv}.
\end{thm}

\proof According to the proof of \thmref{it is eigenvector},
$\phi(w_1^{i_1},\ldots,w_m^{i_m})$ is an eigenvector of $\zz(\ghat)_u
= \on{Fun} \on{Op}_{^L G}(D_u)$ and its eigenvalue is the $^L G$-oper
on $\pone$ whose restriction to $D_u$ is the Miura transformation of
$\ol\nabla|_{D_u}$. But then this oper is nothing but
$\mu_{(z_i),(w_j),\infty}(\ol\nabla)$. \qed

\subsection{Bethe Ansatz for finite-dimensional modules}

Now we specialize to the case when all $\g$-modules $M_i$ are
irreducible and finite-dimensional. Thus, $M_i = V_{\la_i}$ for some
dominant integral highest weight $\la_i$. In this case the oper 
$\mu_{(z_i),(w_j),\infty}(\ol\nabla)$ automatically has no monodromy
around the points $z_1,\ldots,z_N$. Indeed, its restriction to
$D_{z_i}^\times$ is the $^L G$-oper which is the gauge equivalence
class of
$$
\pa_t + p_{-1} - \frac{\la_i}{t-z_i} + {\mb u}(t), \qquad {\mb u}(t)
\in \h^*[[t-z_i]].
$$
Applying the gauge transformation with $\la(t-z_i)^{-1}$, we obtain
the operator
$$
\pa_t + \sum_{j=1}^\ell t^{\langle \la_i,\chal_j \rangle} f_j + {\mb
  u}(t), \qquad {\mb u}(t) \in \h^*[[t-z_i]],
$$
which is regular at $t=z_i$. Hence this $^L G$-oper has trivial
monodromy around $z_i$ for each $i=1,\ldots,N$.

Now consider the Bethe vector $\phi(w_1^{i_1},\ldots,w_m^{i_m}) \in
\bigotimes_{i=1}^N M_{\la_i}$. Let
$\ol\phi(w_1^{i_1},\ldots,w_m^{i_m})$ be its projection onto
$\bigotimes_{i=1}^N V_{\la_i}$. Then we find from \thmref{joint
eigenvalue} that the eigenvalues of the Gaudin algebra on it are given
by the $^L G$-oper
$$
\mu_{(z_i),(w_j),\infty}(\ol\nabla) \in \on{Op}_{^L
  G}(\pone)^{(\la_i)}_{(z_i);\pi(-\chi)} \subset \on{Op}_{^L
  G}^{\on{RS}}(\pone)^{(\la_i)}_{(z_i);\pi(-\chi)},
$$
where $\on{Op}_{^L G}(\pone)^{(\la_i)}_{(z_i);\pi(-\chi)}$ stands for
the monodromy-free locus in $\on{Op}_{^L
  G}^{\on{RS}}(\pone)^{(\la_i)}_{(z_i);\pi(-\chi)}$.

This is compatible with the statement of \thmref{subset} which states
that the eigenvalue on {\em any} eigenvector in $\bigotimes_{i=1}^N
V_{\la_i}$ belongs to this monodromy-free locus.

An interesting problem is the {\em completeness of the Bethe Ansatz}:
is it true that the Bethe eigenvectors gives us a basis of
$\bigotimes_{i=1}^N V_{\la_i}$ for generic $z_1,\ldots,z_N$ and
$\chi$? If this is so, this would mean that there is a basis of
$\bigotimes_{i=1}^N V_{\la_i}$ labeled by solutions of the Bethe
Ansatz equations \eqref{bethe}.

In the case of the ordinary Gaudin model (corresponding to $\chi=0$)
this problem has been investigated in great detail. The completeness
of the Bethe Ansatz has been proved in some special cases in
\cite{SV,MV1,S}. However, examples constructed in \cite{MV2} show that
Bethe Ansatz may be incomplete for some fixed highest weights $\la_i$
and all possible values of $z_i$'s.

Completeness of the Bethe Ansatz in the ordinary Gaudin model is
discussed in detail \cite{F:faro}, Sect. 5. According to Theorem
2.7,(3) \cite{F:faro}, which is recalled in \thmref{reg sing eig}, the
joint eigenvalues of the Gaudin algebra are realized as a subset of
the set
$$
\on{Op}_{^L G}(\pone)_{(z_i),\infty;(\la_i),\la_\infty} = \on{Op}_{^L
  G}(\pone)^{(\la_i)}_{(z_i);\pi(0)}
$$
of monodromy-free opers with $\chi=0$, so that they have regular
singularity at $\infty$. It was conjectured in \cite{F:faro} that this
inclusion is actually a bijection. We can try to prove this by
associating to an oper in $\on{Op}_{^L
G}(\pone)^{(\la_i)}_{(z_i);\pi(0)}$ an eigenvector of the Gaudin
algebra in $\bigotimes_{i=1}^N V_{\la_i}$. This can be done by using
the {\em Miura opers}, as explained in \cite{F:faro}. The problem is
caused by the {\em degenerate} opers in $\on{Op}_{^L
G}(\pone)^{(\la_i)}_{(z_i);\pi(0)}$, in the terminology of
\cite{F:opers,F:faro}.

Each non-degenerate oper gives rise to a Bethe vector of the form
$\ol\phi(w_1^{i_1},\ldots,w_m^{i_m})$. If all opers in $\on{Op}_{^L
G}(\pone)^{(\la_i)}_{(z_i);\pi(0)}$ are non-degenerate and the
corresponding Bethe vectors are non-zero, then the Bethe Ansatz is
complete (see \cite{F:opers}, Prop. 4.10, and \cite{F:faro},
Prop. 5.5). However, a degenerate oper does not give rise to a Bethe
vector. Thus, we may have a bijection between $\on{Op}_{^L
G}(\pone)^{(\la_i)}_{(z_i);\pi(0)}$ and the set of joint eigenvalues
of the Gaudin algebra on $\bigotimes_{i=1}^N V_{\la_i}$ (as
conjectured in \cite{F:faro}), but the Bethe Ansatz may still be
incomplete because of the presence of degenerate opers, to which we
cannot attach Bethe vectors (we believe that this is the reason behind
the counterexample of \cite{MV2}). However, as explained in
\cite{F:faro}, Sect. 5.5, the construction of the Bethe vectors
presented above may be generalized so as to enable us to attach
eigenvectors (of a slightly different form) even to degenerate
opers. This gives us a way to construct an eigenbasis of the Gaudin
algebra in $\bigotimes_{i=1}^N V_{\la_i}$ even if there are degenerate
opers in $\on{Op}_{^L G}(\pone)^{(\la_i)}_{(z_i);\pi(0)}$.\footnote{We
note that for $\g=\sw_n$ and generic $\chi$ it follows from
\cite{Ryb2} that the action of
$\A^{(1),1}_{(z_i),\infty}(\g)_{(0),\chi}$ on $\bigotimes_{i=1}^N
V_{\la_i}$ is diagonalizable and has simple spectrum.}

The approach of \cite{F:faro} may be generalized to the case of the
Gaudin model with regular semi-simple $\chi$. However, the notion of
Miura oper becomes more subtle here because the oper connection has
irregular singularity at $\infty$. We hope to discuss this question
in more detail elsewhere.

\subsection{The case of ${\mc A}_\chi$}

Now we specialize the above results to the case of the Gaudin algebra
${\mc A}_\chi = \A^{1,1}_{0,\infty}(\g)_{0,\chi}$, corresponding
to $N=1$ with $z_1 =0$ and $\la_1 = \la$. We fix a regular semi-simple
$\chi \in \g^*$.

In this case the Bethe vectors in the Verma module $M_\la$ have the
form
\begin{multline}    \label{bethe for one mod}
\phi(w_1^{i_1},\ldots,w_m^{i_m}) = \\ \sum_{\sigma \in S_m}
\frac{f_{i_{\sigma(1)}} f_{i_{\sigma(2)}} \ldots
f_{i_{\sigma(m)}}}{(w_{{\sigma(1)}} -
w_{{\sigma(2)}})(w_{{\sigma(2)}} - w_{{\sigma(3)}}) \ldots
(w_{{\sigma(m-1)}}-w_{{\sigma(m)}}) w_{{\sigma(m)}}} v_{\la},
\end{multline}
where the sum is over all permutations on $m$ letters. This vector has
the weight
$$
\la - \sum_{j=1}^m \al_{i_j}.
$$

The corresponding Bethe Ansatz equations \eqref{bethe} have the form
\begin{equation}    \label{one mod}
\frac{\langle \chal_{i_j},\la \rangle}{w_j} -
\sum_{s\neq j} \frac{\langle \chal_{i_j},\al_{i_s} \rangle}{w_j-w_s} =
\langle \chal_{i_j},\chi \rangle, \quad j=1,\ldots,m.
\end{equation}

If the Bethe Ansatz for an irreducible $\g$-module $V_\la$ is complete
for a particular $\chi \in \g^*$, then the Bethe vectors \eqref{bethe
for one mod} give us a basis of $V_\la$ labeled by the solutions of
the equations \eqref{one mod}. We remark that in the case of
$\g=\sw_n$ it was shown in \cite{Ryb2} that ${\mc A}_\chi$ is closely
related to the Gelfand-Zetlin algebra. This suggests that for a
general simple Lie algebra $\g$ the Gaudin algebra ${\mc A}_\chi$ may
give us a new powerful tool for analyzing finite-dimensional
representations.

\end{document}